\documentclass[preprint,12pt]{elsarticle}  % arXiv / MOX-report
\usepackage[a4paper, total={6.5in, 8.66in}]{geometry}
\pdfoutput=1

\let\today\relax
\makeatletter
\def\ps@pprintTitle{%
    \let\@oddhead\@empty
    \let\@evenhead\@empty
    \def\@oddfoot{\footnotesize\itshape
         {Submitted preprint} \hfill\today}%
    \let\@evenfoot\@oddfoot
    }
\makeatother

\usepackage[english]{babel}
\usepackage[utf8]{inputenc}
\usepackage[printonlyused]{acronym}
\usepackage{amsmath}
\usepackage{amssymb}
\usepackage{mathtools}
\usepackage{mathrsfs}
\usepackage{siunitx}
\usepackage{empheq}
\usepackage{hyperref}
\usepackage[capitalise]{cleveref}
\usepackage{booktabs}
\usepackage{multirow}

\biboptions{sort&compress}

\graphicspath{{figures/}}

\makeatletter
\newcommand{\pushright}[1]{\ifmeasuring@#1\else\omit\hfill$\displaystyle#1$\fi\ignorespaces}
\newcommand{\pushleft}[1]{\ifmeasuring@#1\else\omit$\displaystyle#1$\hfill\fi\ignorespaces}
\makeatother

\usepackage[dvipsnames]{xcolor}
\hypersetup{colorlinks=true}
\usepackage[color=Apricot,colorinlistoftodos]{todonotes}

\everymath{\displaystyle}

\newcommand{\mLperMsquare}{[\si{\milli\liter \per \meter^{2}}]}

\newcommand{\lifex}{\texttt{life\textsuperscript{x}}}

\newcommand{\PhyEP}{\mathscr{E}}
\newcommand{\PhyEPIMEX}{\mathscr{E}_{\mathrm{IMEX}}}
\newcommand{\PhyIon}{\mathscr{I}}
\newcommand{\PhyIonIMEX}{\mathscr{I}_{\mathrm{IMEX}}}
\newcommand{\PhyMec}{\mathscr{M}}
\newcommand{\PhyMecI}{\mathscr{M}_{\mathrm{I}}}
\newcommand{\PhyAct}{\mathscr{A}}
\newcommand{\PhyActIMEX}{\mathscr{A}_{\mathrm{IMEX}}}
\newcommand{\PhyCirc}{\mathscr{C}}
\newcommand{\PhyCircE}{\mathscr{C}_{\mathrm{E}}}
\newcommand{\PhyCoupl}{\mathscr{V}} % volume-consistency

\newcommand{\Nsub}{N_{\mathrm{sub}}}
\newcommand{\h}{h}
\newcommand{\generica}{\boldsymbol{a}}

\newcommand{\VLVthreedimH}{V_{\LV}^{\mathrm{3D}}}
\newcommand{\VRVthreedimH}{V_{\RV}^{\mathrm{3D}}}
\newcommand{\VLAthreedimH}{V_{\LA}^{\mathrm{3D}}}
\newcommand{\VRAthreedimH}{V_{\RA}^{\mathrm{3D}}}

\newcommand{\Vthreedimk}{V_{k}^{\mathrm{3D}}}

\newcommand{\x}{\mathbf{x}}
\newcommand{\xtilde}{\widetilde{\mathbf{x}}}
\newcommand{\RA}{\mathrm{RA}}
\newcommand{\LA}{\mathrm{LA}}
\newcommand{\RV}{\mathrm{RV}}
\newcommand{\LV}{\mathrm{LV}}
\newcommand{\PT}{\mathrm{PT}}
\newcommand{\AO}{\mathrm{AO}}
\newcommand{\TV}{\mathrm{TV}}
\newcommand{\MV}{\mathrm{MV}}
\newcommand{\AV}{\mathrm{AV}}
\newcommand{\PV}{\mathrm{PV}}

\newcommand{\KinCH}{k \in \{ \RA,\LA,\RV,\LV \}}

\newcommand{\IinCH}{i\in \{ \RA,\LA,\RV,\LV \}}
\newcommand{\IinCHV}{i\in \{ \RA,\LA,\mathrm{V} \}}
\newcommand{\IinCHAR}{i\in \{ \RA,\LA,\RV,\LV,\PT,\AO \}}
\newcommand{\IinValve}{i\in \{ \TV,\MV,\PV,\AV \}}

\newcommand{\Ventr}{\mathrm{V}}
\newcommand{\Atria}{\mathrm{A}}
\newcommand{\Arter}{\mathrm{AR}}
\newcommand{\Valve}{\mathrm{valve}}
\newcommand{\Caps}{\mathrm{caps}}

\newcommand{\Myo}{\mathrm{myo}}
\newcommand{\Epi}{\mathrm{epi}}
\newcommand{\Endo}{\mathrm{endo}}
\newcommand{\Fat}{\mathrm{EAT}}

\newcommand{\OmegaTilde}{\widetilde{\Omega}}

\newcommand{\OmegaMyo}{\Omega^{\Myo}_0}
\newcommand{\OmegaVentr}{\Omega^{\Ventr}_0}
\newcommand{\OmegaRA}{\Omega^{\RA}_0}
\newcommand{\OmegaLA}{\Omega^{\LA}_0}
\newcommand{\OmegaPT}{\Omega^{\PT}_0}
\newcommand{\OmegaAO}{\Omega^{\AO}_0}
\newcommand{\OmegaValve}{\Omega^{\Valve}_0}
\newcommand{\OmegaCaps}{\Omega^{\Caps}_0}

\newcommand{\GammaEpi}{\Gamma^{\Epi}}

\newcommand{\GammaEpiFat}{\Gamma^{\Epi,\Fat}}
\newcommand{\GammaEpiPF}{\Gamma^{\Epi,\mathrm{PF}}}
\newcommand{\GammaEpiArter}{\Gamma^{\Epi,\Arter}}

\newcommand{\GammaEndok}{\Gamma^{\Endo,k}}
\newcommand{\GammaEndoI}{\Gamma^{\Endo,i}}
\newcommand{\GammaEndoLV}{\Gamma^{\Endo,\LV}}
\newcommand{\GammaEndoRV}{\Gamma^{\Endo,\RV}}
\newcommand{\GammaEndoLA}{\Gamma^{\Endo,\LA}}
\newcommand{\GammaEndoRA}{\Gamma^{\Endo,\RA}}
\newcommand{\GammaEndoAo}{\Gamma^{\Endo,\AO}}
\newcommand{\GammaEndoPT}{\Gamma^{\Endo,\PT}}

\newcommand{\GammaRings}{\Gamma^{\mathrm{rings}}}

\newcommand{\w}{{\boldsymbol{w}}}
\newcommand{\wAtria}{\w_1}%{\w_\Atria}
\newcommand{\wVentr}{\w_2}%{\w_\Ventr}
\newcommand{\wCa}{w_{\mathrm{Ca}}}
\newcommand{\wCaAtria}{w_{1,\mathrm{Ca}}}
\newcommand{\wCaVentr}{w_{2,\mathrm{Ca}}}
\newcommand{\Disp}{{\mathbf{d}}}
\newcommand{\DispTilde}{{\widetilde{\Disp}}}
\newcommand{\Circ}{{\boldsymbol{c}}}

\newcommand{\ZeroDim}{\mathrm{0D}}
\newcommand{\ThreeDim}{\mathrm{3D}}

\newcommand{\f}{\mathbf{f}}
\newcommand{\s}{\mathbf{s}}
\newcommand{\n}{\mathbf{n}}

\newcommand{\GammaA}{\Gamma_0^{\mathrm{a}}}
\newcommand{\GammaB}{\Gamma_0^{\mathrm{b}}}
\newcommand{\GammaN}{\Gamma_0^{\mathrm{n}}}

\newcommand{\GammaClosedA}{\overline{\Gamma}_0^{\mathrm{a}}}
\newcommand{\GammaClosedB}{\overline{\Gamma}_0^{\mathrm{b}}}
\newcommand{\GammaClosedN}{\overline{\Gamma}_0^{\mathrm{n}}}

\newcommand{\xiA}{\xi_{\mathrm{a}}}
\newcommand{\xiB}{\xi_{\mathrm{b}}}

\newcommand{\eL}{\widehat{\boldsymbol{e}}_{\mathrm{\ell}}}
\newcommand{\eT}{\widehat{\boldsymbol{e}}_{\mathrm{t}}}
\newcommand{\eN}{\widehat{\boldsymbol{e}}_{\mathrm{n}}}

\newcommand{\thetaEndoi}{\theta_{\mathrm{endo,}j}}
\newcommand{\thetaEpii}{\theta_{\mathrm{epi,}j}}

\newcommand{\EPchim}{\chi_\mathrm{m}}
\newcommand{\EPCm}{C_\mathrm{m}}
\newcommand{\EPIion}{{\mathcal{I}_{\mathrm{ion}}}}
\newcommand{\EPIapp}{{\mathcal{I}_{\mathrm{app}}}}

\newcommand{\EPdiffTens}{\boldsymbol{D}_{\mathrm{M}}}

\newcommand{\EPrhsGating}{\boldsymbol{H}}

\newcommand{\CondF}[1]{\sigma_\f^{#1}}
\newcommand{\CondS}[1]{\sigma_\s^{#1}}
\newcommand{\CondN}[1]{\sigma_\n^{#1}}
\newcommand{\CondK}[1]{\sigma_{\mathbf{k}}^{#1}}

\newcommand{\aXB}{a_{\mathrm{XB}}}

\newcommand{\CircRhs}{\boldsymbol{D}}

\newcommand{\Cai}{{[\mathrm{Ca}^{2+}]_{\mathrm{i}}}}

\newcommand{\SL}{\mathrm{SL}} % sarcomere length

\newcommand{\Tens}{T_\mathrm{a}} % Active tension [Pa]
\newcommand{\Stiff}{K_\mathrm{a}} % Active stiffness [Pa]
\newcommand{\TensConst}{\overline{T}_\mathrm{a}} % Active tension [Pa]
\newcommand{\TensTilde}{\widetilde{T}_\mathrm{a}} % Active tension [Pa]
\newcommand{\TensInit}{T_{\mathrm{a},0}} % Active tension [Pa]
\newcommand{\ActStateHF}{\mathbf{z}}

\newcommand{\ActRhs}{\boldsymbol{K}}

\newcommand{\Inv}[1]{{\mathcal{I}_{#1}}}
\newcommand{\IIVf}{\Inv{4f}}

\newcommand{\IIVs}{\Inv{4s}}
\newcommand{\IIVn}{\Inv{4n}}

\newcommand{\BCmecCepiN}{{C_\bot^{\mathrm{epi}}}}
\newcommand{\BCmecCepiNpf}{{C_\bot^{\mathrm{epi,PF}}}}
\newcommand{\BCmecCepiNfat}{{C_\bot^{\mathrm{epi,EAT}}}}
\newcommand{\BCmecKepiN}{{K_\bot^{\mathrm{epi}}}}
\newcommand{\BCmecKepiNpf}{{K_\bot^{\mathrm{epi,PF}}}}
\newcommand{\BCmecKepiNfat}{{K_\bot^{\mathrm{epi,EAT}}}}

\newcommand{\mecF}{\mathbf{F}}
\newcommand{\mecC}{\mathbf{C}}
\newcommand{\tenspiola}{\mathbf{P}}
\newcommand{\identity}{\mathbf{I}}

\newcommand{\mecNref}{{\mathbf{N}}}

\newcommand{\Clrv}{C_\mathrm{lrv}}
\newcommand{\XiHat}{\hat{\xi}}

\newcommand{\VLA}[1][(t)]{V_{\LA}#1}
\newcommand{\VLV}[1][(t)]{V_{\LV}#1}
\newcommand{\VRA}[1][(t)]{V_{\RA}#1}
\newcommand{\VRV}[1][(t)]{V_{\RV}#1}
\newcommand{\VI}[1][(t)]{V_{i}#1}

\newcommand{\ParSYS}{p_{\mathrm{AR}}^{\mathrm{SYS}}(t)}
\newcommand{\ParSYSn}{p_{\mathrm{AR}}^{\mathrm{SYS},n}}
\newcommand{\ParPUL}{p_{\mathrm{AR}}^{\mathrm{PUL}}(t)}
\newcommand{\ParPULn}{p_{\mathrm{AR}}^{\mathrm{PUL},n}}
\newcommand{\PvnSYS}{p_{\mathrm{VEN}}^{\mathrm{SYS}}(t)}
\newcommand{\PvnPUL}{p_{\mathrm{VEN}}^{\mathrm{PUL}}(t)}
\newcommand{\QarSYS}{Q_{\mathrm{AR}}^{\mathrm{SYS}}(t)}
\newcommand{\QarPUL}{Q_{\mathrm{AR}}^{\mathrm{PUL}}(t)}
\newcommand{\QvnSYS}{Q_{\mathrm{VEN}}^{\mathrm{SYS}}(t)}
\newcommand{\QvnPUL}{Q_{\mathrm{VEN}}^{\mathrm{PUL}}(t)}

\newcommand{\PLA}[1][(t)]{p_{\LA}#1}
\newcommand{\PLV}[1][(t)]{p_{\LV}#1}
\newcommand{\PRA}[1][(t)]{p_{\RA}#1}
\newcommand{\PRV}[1][(t)]{p_{\RV}#1}
\newcommand{\PAO}[1][(t)]{p_{\AO}#1}
\newcommand{\PPT}[1][(t)]{p_{\PT}#1}

\newcommand{\PLAconst}{\overline{p}_{\LA}}
\newcommand{\PLVconst}{\overline{p}_{\LV}}
\newcommand{\PRAconst}{\overline{p}_{\RA}}
\newcommand{\PRVconst}{\overline{p}_{\RV}}
\newcommand{\PAOconst}{\overline{p}_{\AO}}
\newcommand{\PPTconst}{\overline{p}_{\PT}}
\newcommand{\PIconst}{\overline{p}_{i}}

\newcommand{\PLAtilde}{\widetilde{p}_{\LA}}
\newcommand{\PLVtilde}{\widetilde{p}_{\LV}}
\newcommand{\PRAtilde}{\widetilde{p}_{\RA}}
\newcommand{\PRVtilde}{\widetilde{p}_{\RV}}
\newcommand{\PAOtilde}{\widetilde{p}_{\AO}}
\newcommand{\PPTtilde}{\widetilde{p}_{\PT}}
\newcommand{\PItilde}{\widetilde{p}_{i}}

\newcommand{\PIinit}{p_{i,0}}

\newcommand{\QAV}{Q_{\AV}}
\newcommand{\QMV}{Q_{\MV}}
\newcommand{\QTV}{Q_{\TV}}
\newcommand{\QPV}{Q_{\PV}}

\newcommand{\VLVzerodim}{V_{\LV}^{\ZeroDim}}
\newcommand{\VRVzerodim}{V_{\RV}^{\ZeroDim}}
\newcommand{\VLAzerodim}{V_{\LA}^{\ZeroDim}}
\newcommand{\VRAzerodim}{V_{\RA}^{\ZeroDim}}
\newcommand{\VIzerodim}{V_{i}^{\ZeroDim}}
\newcommand{\VKzerodim}{V_{k}^{\ZeroDim}}

\newcommand{\VIthreedim}{V_{i}^{\ThreeDim}}

\newcommand{\CarSYS}{C_{\mathrm{AR}}^{\mathrm{SYS}}}
\newcommand{\CarPUL}{C_{\mathrm{AR}}^{\mathrm{PUL}}}
\newcommand{\CvnSYS}{C_{\mathrm{VEN}}^{\mathrm{SYS}}}
\newcommand{\CvnPUL}{C_{\mathrm{VEN}}^{\mathrm{PUL}}}

\newcommand{\RarSYS}{R_{\mathrm{AR}}^{\mathrm{SYS}}}
\newcommand{\RarPUL}{R_{\mathrm{AR}}^{\mathrm{PUL}}}
\newcommand{\RvnSYS}{R_{\mathrm{VEN}}^{\mathrm{SYS}}}
\newcommand{\RvnPUL}{R_{\mathrm{VEN}}^{\mathrm{PUL}}}

\newcommand{\LarSYS}{L_{\mathrm{AR}}^{\mathrm{SYS}}}
\newcommand{\LarPUL}{L_{\mathrm{AR}}^{\mathrm{PUL}}}
\newcommand{\LvnSYS}{L_{\mathrm{VEN}}^{\mathrm{SYS}}}
\newcommand{\LvnPUL}{L_{\mathrm{VEN}}^{\mathrm{PUL}}}

\newcommand{\Rmin}{R_{\mathrm{min}}}
\newcommand{\Rmax}{R_{\mathrm{max}}}
\begin{document}

\begin{frontmatter}

    \title{A comprehensive and biophysically detailed computational model of the whole human heart electromechanics}
    
    \author[inst1]{Marco Fedele} \corref{cor1}\ead{marco.fedele@polimi.it}
    \author[inst1]{Roberto Piersanti}
    \author[inst1]{Francesco Regazzoni}
    \author[inst1]{Matteo Salvador}
    \author[inst1]{Pasquale Claudio Africa}
    \author[inst1]{Michele Bucelli}
    \author[inst1]{Alberto Zingaro}
    \author[inst1]{Luca Dede'}
    \author[inst1,inst2]{Alfio Quarteroni}
    
    \cortext[cor1]{Corresponding author.}
    
    \affiliation[inst1]{organization={MOX - Department of Mathematics, Politecnico di Milano},%Department and Organization
                addressline={Piazza Leonardo da Vinci, 32}, 
                city={Milano},
                postcode={20133},
                %state={Italy},
                country={Italy}}
    
    \affiliation[inst2]{organization={Mathematics Institute (Professor Emeritus), École Polytechnique Fédérale de Lausanne}, %Department and Organization
                addressline={Av. Piccard},
                city={Lausanne},
                postcode={CH-1015},
                %state={State Two},
                country={Switzerland}}
    
\begin{abstract}
While ventricular electromechanics is extensively studied in both physiological and pathological conditions, four-chamber heart models have only been addressed recently;
most of these works however neglect atrial contraction.
Indeed, as atria are characterized by a complex anatomy and a physiology that is strongly influenced by the ventricular function, developing computational models able to capture the physiological atrial function and atrioventricular interaction is very challenging.
In this paper, we propose a biophysically detailed electromechanical model of the whole human heart that considers both atrial and ventricular contraction.
Our model includes:
i) an anatomically accurate whole-heart geometry;
ii) a comprehensive myocardial fiber architecture;
iii) a biophysically detailed microscale model for the active force generation;
iv) a 0D closed-loop model of the circulatory system, fully-coupled with the mechanical model of the heart;
v) the fundamental interactions among the different \textit{core models}, such as the mechano-electric feedback or the fibers-stretch and fibers-stretch-rate feedbacks;
vi) specific constitutive laws and model parameters for each cardiac region.
Concerning the numerical discretization, we propose an efficient segregated-intergrid-staggered scheme and we employ recently developed stabilization techniques -- regarding the circulation and the fibers-stretch-rate feedback -- that are crucial to obtain a stable formulation in a four-chamber scenario.
We are able to reproduce the healthy cardiac function for all the heart chambers, in terms of pressure-volume loops, time evolution of pressures, volumes and fluxes, and three-dimensional cardiac deformation, with unprecedented matching (to the best of our knowledge) with the expected physiology.
We also show the importance of considering atrial contraction, fibers-stretch-rate feedback and suitable stabilization techniques, by comparing the results obtained with and without these features in the model.
The proposed model represents the state-of-the-art electromechanical model of the iHEART ERC project -- an Integrated Heart Model for the Simulation of the Cardiac Function -- and is a fundamental step toward the building of physics-based digital twins of the human heart.
\end{abstract}
    
    %     \includegraphics{grabs}
    
    % https://www.elsevier.com/authors/tools-and-resources/highlights
    \begin{highlights}
    \item We propose a novel whole-heart electromechanical model including atrial contraction
    \item Unprecedented match with healthy cardiac physiology
    \item Physiological atrial eight-shaped pressure-volume loops
    \item Fibers-stretch-rate feedback essential to avoid unphysiologically large fluxes
    \item Crucial interplay among accurate mathematical models and stable numerical methods
    \end{highlights}
    
    \begin{keyword}
    Multiphysics and multiscale modeling \sep Whole-heart modeling \sep Cardiac Electromechanics \sep Computational Cardiology \sep High Performance Computing \sep Cardiac Digital Twin
    % \sep Atrial function modeling
    \end{keyword}
    
\end{frontmatter}

    %\linenumbers  % switch line numbers on for the whole article

%%%

\newcommand{\AcroLine}[1]{\acs{#1}&\acl{#1}}
\newcommand{\AcroLineDouble}[2]{\acs{#1}&\acl{#1}&&\acs{#2}&\acl{#2}}
\newcommand{\AcroLineTriple}[3]{\acs{#1}&\acl{#1}&&\acs{#2}&\acl{#2}&&\acs{#3}&\acl{#3}}

% to define without listing ...
\newcommand{\MyAcro}[2]{\acrodef{#1}{#2}}

% to define with listing ...
%\newcommand{\MyAcro}[2]{\acro{#1}{#2}}
% \section{List of abbreviations}
% \label{sec:appendix_acro}

%\begin{acronym}
\MyAcro{LV}{Left Ventricle}
\MyAcro{RV}{Right Ventricle}
\MyAcro{LA}{Left Atrium}
\MyAcro{RA}{Right Atrium}
\MyAcro{MV}{Mitral Valve}
\MyAcro{TV}{Tricuspid Valve}
\MyAcro{AV}{Aortic Valve}
\MyAcro{PV}{Pulmonary Valve}
\MyAcro{AO}{Aorta}
\MyAcro{PT}{Pulmonary Trunk}
\MyAcro{PVs}{Pulmonary Veins}
\MyAcro{SupVC}{Superior Vena Cava}
\MyAcro{InfVC}{Inferior Vena Cava}
\MyAcro{PaMs}{Papillary Muscles}
\MyAcro{ChT}{Chordae Tendineae}
\MyAcro{PF}{Pericardial Fluid}
\MyAcro{LAA}{Left Atrial Appendage}
\MyAcro{RAA}{Right Atrial Appendage}
\MyAcro{BB}{Bachmann's Bundle}
\MyAcro{EAT}{Epicardial Adipose Tissue}
\MyAcro{SAN}{SinoAtrial Node}
\MyAcro{AVN}{AtrioVentricular Node}
\MyAcro{HB}{His Bundle}
\MyAcro{LBB}{Left Bundle Branch}
\MyAcro{RBB}{Right Bundle Branch}
\MyAcro{PFs}{Purkinje Fibers}
\MyAcro{CrT}{Crista Terminalis}
\MyAcro{PeMs}{Pectinate Muscles}

\MyAcro{IVR}{IsoVolumetric Relaxation}
\MyAcro{VPF}{Ventricular Passive Filling}
\MyAcro{AC}{Atrial Contraction}
\MyAcro{IVC}{IsoVolumetric Contraction}
\MyAcro{VE}{Ventricular Ejection}
\MyAcro{EDV}{End Diastolic Volume}
\MyAcro{ESV}{End Systolic Volume}
\MyAcro{SV}{Stroke Volume}

\MyAcro{RBM}{Rule-Based-Method}
\MyAcro{LDRBM}{Laplace-Dirichlet \acs{RBM}}

\MyAcro{CRN}{\citet{CRN}}
\MyAcro{TTP06}{\citet{TTP06}}
\MyAcro{RDQ20}{\citet{regazzoni2020biophysically}}

\MyAcro{MEF}{Mechano-Electric Feedback}

\MyAcro{BDF}{Backward Differentiation Formula}
\MyAcro{BDF1}{\acs{BDF} of order 1}
\MyAcro{BDF2}{\acs{BDF} of order 2}
\MyAcro{FE}{Finite Element}
\MyAcro{HPC}{High Performance Computing}
\MyAcro{IMEX}{Implicit-Explicit}
\MyAcro{ICI}{Ionic Current Interpolation}
\MyAcro{DOFs}{Degrees Of Freedom}
%\end{acronym}

\begin{table}[ht]
\centering
  \begin{tabular}{@{}lllll@{}}
  \toprule
  %\cline{1-2}
  %\cline{4-5}
    \multicolumn{2}{l}{\textbf{Cardiac anatomy}} && \multicolumn{2}{l}{\textbf{Cardiac cycle}} \\
  %\cline{1-2}
  %\cline{4-5}
    \AcroLineDouble{AO}{AC}\\
    \AcroLineDouble{AV}{EDV}\\
    \AcroLineDouble{AVN}{ESV}\\
    \AcroLineDouble{BB}{IVC}\\
    \AcroLineDouble{ChT}{IVR}\\
    \AcroLineDouble{CrT}{SV}\\
    \AcroLineDouble{EAT}{VE}\\
    \AcroLineDouble{HB}{VPF}\\
    \AcroLine{InfVC} \\
  %\cline{4-5}
    \AcroLine{LA} && \multicolumn{2}{l}{\textbf{Modeling}} \\
  %\cline{4-5}
    \AcroLineDouble{LAA}{BDF}\\
    \AcroLineDouble{LBB}{BDF1}\\
    \AcroLineDouble{LV}{BDF2}\\
    \AcroLineDouble{MV}{CRN}\\
    \AcroLineDouble{PF}{DOFs}\\
    \AcroLineDouble{PFs}{FE}\\
    \AcroLineDouble{PT}{HPC}\\
    \AcroLineDouble{PV}{ICI}\\
    \AcroLineDouble{PVs}{IMEX}\\
    \AcroLineDouble{PaMs}{LDRBM}\\
    \AcroLineDouble{PeMs}{MEF}\\
    \AcroLineDouble{RA}{RBM}\\
    \AcroLineDouble{RAA}{RDQ20}\\
    \AcroLineDouble{RBB}{TTP06}\\
    \AcroLine{RV}\\
    \AcroLine{SAN}\\
    \AcroLine{SupVC}\\
    \AcroLine{TV}\\
  \bottomrule
  \end{tabular}
  \caption{List of abbreviations.}
  \label{tab:acro}
\end{table}
\acresetall

%%%%%%%%
\section{Introduction}\label{sec:intro}

We propose a biophysically detailed, numerically stable and accurate computational model of the electromechanics of the whole human heart, considering an active contraction model for both atria and ventricles.
Our model can accurately reproduce the healthy cardiac function, representing a fundamental step toward the building of physics-based digital twins of the human heart.

Computational models of the cardiac function are progressively increasing their role in cardiology, revealing diagnostic information, contributing to the development of new therapies and promising patient-specific treatments based on individual pathophysiology \citep{Trayanova2012,Gray2018,Niederer2019}.
Successful examples can be found in the context of cardiac electrophysiology \citep{Trayanova,Gillette2021,Arevalo2016,Prakosa2018,Frontera2022}, electromechanics \citep{Marx2020,Jung2022,Salvador2021,Peirlinck2021,Peirlinck2022} and fluid-dynamics \citep{Karabelas2022,Mittal2016,Santiago2018,Verzicco2022}.

The growing demand for computational models in clinical applications requires the development of increasingly detailed mathematical models and efficient numerical methods \citep{Augustin2016,Quarteroni2017,Gerbi2019,Viola2020,Strocchi2020,Regazzoni2022,Piersanti2022biv,Stella2022,Zingaro2022,Cicci2023}.
In the context of cardiac electromechanics, a biophysically detailed model of the human heart encompasses all the multiscale and multiphysics processes underlying the cardiac function, ranging from the cellular (microscale) to the organ (macroscale) level, such as the propagation of the electrical signal, the active and passive mechanics, and the interaction with the circulatory system~\citep{Quarteroni2017}.
Moreover, the biophysics of the heart tissue is substantially different among atria, ventricles and non-conductive regions (e.g. valves, arteries).
Further modeling difficulties are given by the complex anatomy made up of many components with non-trivial shapes, each of which plays an important role in the cardiac function~\citep{Sanchez2015,Katz2010}.

All these complex aspects make accurate simulation of the cardiac cycle -- characterized by highly coordinated electrical, mechanical and valvular events -- a very challenging subject still not fully addressed.
In particular, the literature lacks electromechanical models of the entire human heart that take into account both atrial and ventricular contraction in detailed whole-heart geometries.
While ventricular electromechanics in image-based geometries is extensively studied in both physiological and pathological conditions~\cite{Usyk2002,Smith,Goktepe,Nordsletten,Trayanova,Genet2014,Quarteroni2017,Quarteroni2019,Salvador2021,Regazzoni2022,Piersanti2022biv},
whole-heart models emerged only in recent years~\citep{Sugiura2012a,Fritz2014,Baillargeon2014,Augustin2016,Land2018,Santiago2018,Pfaller2019,Strocchi2020,Strocchi2020b,Piersanti2021fibers,Gerach2021,DelCorso2022ep4ch}.
Some studies focus only on electrophysiology \citep{Piersanti2021fibers,DelCorso2022ep4ch} or, if they consider electromechanics, include the atrial muscle only as passive tissue \citep{Sugiura2012a,Fritz2014,Augustin2016,Santiago2018,Pfaller2019,Strocchi2020b}.
More specifically, \citet{Sugiura2012a} review the essential methodologies for a multiscale and multiphysics heart model using the University of Tokyo whole-heart simulator.
However, the electromechanical results are limited to the ventricles, as well as those of related papers using this simulator \citep{Yoneda2021,Sugiura2022}.
\citet{Fritz2014} propose a whole-heart image-based model of the ventricular contraction that considers the interaction with passive atria, pericardium and surrounding organs, demonstrating their impact on the modeling of a physiological heart deformation.
% especially in terms of atrioventricular plane displacement.
\citet{Augustin2016} focus their study on the importance of considering anatomically accurate image-based geometries of the entire heart. % a bidirectional coupling between electrophysiology and mechanics.
They also develop novel numerical techniques that allow solving these complex problems in high-resolution computational meshes.
\citet{Santiago2018} present a fluid-electro-mechanical model of the heart focusing on the ventricles and the arterial flow.
They perform simulations in the anatomically accurate Zygote Solid 3D Heart Model \citep{Zygote2014}, considering simplified passive atria filled with a soft material in their cavity.
Despite this simplification, they show the impact of including atria to achieve physiological ventricular motion.
\citet{Pfaller2019} analyze the importance of proper epicardial boundary conditions in the mechanical model to correctly surrogate the effect of the pericardium and surrounding organs.
\citet{Strocchi2020b} propose simular boundary conditions but considering spatially varying coefficients, to take into account the different stiffness of the surrounding organs.
Both of these studies are based on whole-heart geometries reconstructed from medical images, that are also used to validate the results.
Finally, \citet{Strocchi2020} release a publicly available cohort of four-chamber heart meshes reconstructed from CT-images to facilitate the study of the whole-heart electromechanics.
They also perform simulations of the ventricular electrical activation and contraction on this cohort.

All the aforementioned electromechanical models neglect atrial contraction, that, to the best of our knowledge, is instead considered only in a few works, namely \citep{Baillargeon2014,Land2018,Gerach2021}.
\citet{Baillargeon2014} present the Living Heart project, a simulator of the human cardiac function that includes a phenomenological representation of both ventricular and atrial active contraction.
This simulator has been extensively used in recent years, but mainly to study ventricular pathologies \citep{Baillargeon2015,Genet2016,Peirlinck2021,Peirlinck2022}, while more details on atrial contraction (such as pressure-volume loops) have never been shown.
The work of \citet{Land2018} is the first one focusing on the influence of atrial contraction on the cardiac function, investigating also an atrial fibrillation scenario.
Active contraction is taken into account using the lumped-parameter model previously proposed for the ventricles \citep{Land2017}, by adapting some parameters to the atrial case.
This work shows, as a result, atrial pressure-volume loops that qualitatively tend to the distinguishing physiological eight-shape.
\citet{Gerach2021} use the same active contraction model while also including a three-dimensional representation of the pericardium, the adipose tissue and the beginning of the major vessels.
They show atrial pressure-volume loops that qualitatively match the characteristic eight-shape, representing the most realistic result concerning atrial function available in the literature.
However, blood fluxes across the semilunar valves thereby shown substantially exceed the physiological values.
As we show in our paper, these anomalies can be explained by the lack of the fibers-stretch-rate feedback (between passive mechanics and active force generation model).
Indeed, this feedback is commonly neglected since it can generate, at a numerical level, strong non-physical oscillations~\citep{regazzoni2021oscillation,Gerach2021}.

Compared to the ventricles, the atria exhibit a more complex anatomy and physiology, characterized by a thinner and weaker muscle strongly influenced by ventricular contraction and relaxation.
Consequently, computational models of the atrial function are very challenging and must consider properly calibrated biophysically detailed models of the four chambers in order to obtain physiologically meaningful results.

In this paper, we propose a novel mathematical model of whole-heart electromechanics endowed with biophysically detailed \textit{core models} for electrophysiology, passive mechanics, and ventricular and atrial active contraction.
Specifically, our mathematical model -- that extends the left-ventricular model we have recently proposed in~\citep{Regazzoni2022} -- features several innovative contributions:
\begin{itemize}
  \item an anatomically accurate whole-heart model consisting of detailed geometries for the four chambers, simplified valves acting as electrically insulating regions, and the initial tracts of the arteries;
  \item an accurate myocardial fiber architecture using a novel whole-heart \ac{RBM} that takes into account also the characteristic atrial fiber bundles \citep{Piersanti2021fibers,Piersanti2021phd};
  \item chamber-specific and accurate ionic models for atria and ventricles~\citep{CRN,TTP06};
  \item a biophysically detailed microscale model for the active force generation~\citep{Regazzoni2020};
  \item a 0D closed-loop model of the circulatory system, fully-coupled with the mechanical model~\citep{Regazzoni2022}; %to determine pressure and volume load of the heart chambers
  \item specific spring-damper Robin boundary conditions to model the pericardium and the surrounding organs~\citep{Pfaller2019}.
\end{itemize}
The core models are specifically calibrated for each cardiac compartment and coupled with each other taking into account the most important feedbacks, such as the \ac{MEF} or the fibers-stretch and fibers-stretch-rate feedbacks.
Concerning the numerical discretization, we use the efficient segregated-intergrid-staggered scheme proposed in \citep{Regazzoni2022,Piersanti2022biv} and we employ recently developed stabilization terms -- related to the circulation~\citep{regazzoni2022stabilization} and the fibers-stretch-rate feedback~\citep{regazzoni2021oscillation} -- that are crucial to obtain a stable formulation in a four-chamber scenario.
The numerical models proposed in this work are characterized by high dimensionality and huge computational complexity, thus calling for efficient and accurate computational tools.
To this aim, the solver that we developed for the numerical simulation of the whole-heart electromechanics has been built upon \lifex\footnote[1]{\url{https://lifex.gitlab.io/}}, an in-house \ac{FE} library focused on large-scale cardiac applications in a \ac{HPC} framework.

% To conclude, TAB ... shows an overview of the main features and results of our model compared to the whole-heart electromechanical models existing in the literature.
% \todoinline{
% Aggiungere una tabella di sintesi (politically-correct, ma che dimostri come il nostro sia il miglio modello 4ch attuale) completando \href{https://docs.google.com/spreadsheets/d/1_7yG252FxdAfufBvpa3bW9EP945YbVKsjM79QjZQBSc/edit?usp=sharing}{QUESTO EXCEL}?\\
% @MF per il momento eliminerei quest'ultima parte in arancione.
% }

This paper is structured as follows:
in \cref{sec:heart} we shortly review the anatomy and physiology of the heart;
in \cref{sec:meth_math} we describe the full electromechanical model;
\cref{sec:meth_num} is devoted to the numerical discretization;
in \cref{sec:res} we discuss the numerical results;
finally, in \cref{sec:concl} we draw our conclusions.

%%%%%%%%
%%%%%%%%
\section{Cardiac anatomy and physiology}
\label{sec:heart}
In this section we briefly review the anatomy of the human heart aiming at introducing all the cardiac components that -- with different level of details -- we consider in our electromechanical model.
We also describe the phases of the cardiac cycle, focusing on the differences of the atrial and ventricular function.
For a more in-depth overview of the cardiac anatomy and physiology, we refer to \citep{Opie2004,Anderson2004,Iaizzo2010,Katz2010,Klabunde2011,Askari2019,Verzicco2022}.

\begin{figure}[t]
    \centering
    \includegraphics[width=1.0\textwidth]{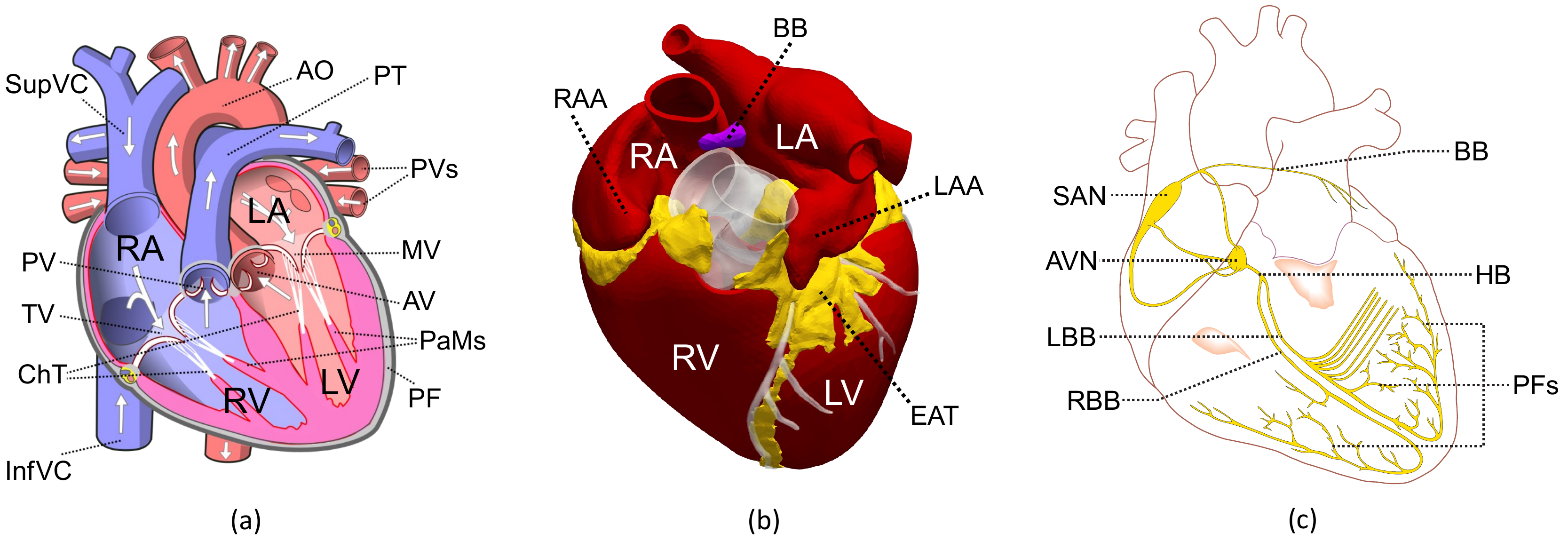}
    \caption{
    The anatomy of the heart:
    (a) a sketch of the internal view of the four chambers and their components (picture elaborated  from \url{https://commons.wikimedia.org/w/index.php?curid=830253});
    (b) an  external view of the cardiac anatomy;
    (c) a sketch of the electrical conduction system (picture elaborated  from \url{https://commons.wikimedia.org/w/index.php?curid=10197958}).
    All the abbreviations are defined in \cref{tab:acro}.}
    \label{fig:heart_anatomy}
\end{figure}

As shown in \cref{fig:heart_anatomy}, (a)-(b), the human heart
% -- a muscular organ that pumps blood in the circulatory system --
is characterized by a very complex anatomy and is made up of several components, each of which plays a crucial role in the cardiac function.
The heart is made up of four muscle chambers:
%whose contraction is regulated by the cardiac electrical conduction system:
\ac{RA} and \ac{LA} on the upper part, \ac{RV} and \ac{LV} on the lower part.
Commonly, \ac{RA} and \ac{RV} are collectively referred to as \textit{right heart} and their left counterparts (\ac{LA} and \ac{LV}) as \textit{left heart}.
The right heart pumps the oxygen-depleted blood -- coming from the systemic venous return and flowing through the \ac{SupVC} and the \ac{InfVC} -- toward the \ac{PT} into the lungs, where oxygenation takes place;
the left heart pumps oxygenated blood -- coming from the lungs through the \ac{PVs} -- toward the \ac{AO} into the systemic circulation, closing the loop of the circulatory system.

The blood flow is regulated by four cardiac valves made of strong fibrous tissue:
the \ac{TV} and \ac{MV} lie in the atrioventricular plane and divide the \ac{RA} and \ac{LA} from the \ac{RV} and \ac{LV}, respectively, also acting as electrical insulators between atria and ventricles;
the \ac{PV} and \ac{AV} connect the \ac{RV} and \ac{LV} to the \ac{PT} and \ac{AO}, respectively.
Valves passively open and close depending on the pressure exerted on their leaflets;
\ac{TV} and \ac{MV} are also supported by \ac{ChT} and \ac{PaMs} to avoid valve prolapse while closed.

The tissue of the cardiac chambers is made up of three layers:
the \textit{endocardium} is the thin innermost layer in direct contact with the blood;
%and where the branches of the electrical conduction system are immersed;
the \textit{myocardium} is the thick muscle layer made of \textit{cardiomyocytes}, the cells responsible for generating contractile force in the heart;
the \textit{epicardium} forms the thin outermost layer.
This latter layer, mainly characterized by a smooth surface, features a complex rough anatomy in some regions (see \cref{fig:heart_anatomy}, (b)).
In particular, the atrioventricular regions, the initial part of the arteries (\ac{PT} and \ac{AO}) and the presence of the \ac{LAA} and the \ac{RAA} contribute to create some empty regions among the different cardiac components.
These regions are filled of the \ac{EAT}, a visceral fat deposit that creates a sort of soft pillow among the nearby cardiac components and contributes to make the external surface of the heart a smooth surface.

The whole heart -- including \ac{EAT} -- is surrounded by the \textit{pericardium}, a sac that holds the heart in place.
This sac is filled with the \ac{PF} which allows the free sliding of the heart external surface, thus also allowing the volume of the four chambers to increase or decrease during the different phases of the cardiac cycle.

The cardiac cycle is a highly coordinated, temporally related series of electrical, mechanical, and valvular events \citep{Pagel2018}.
The contraction of the four cardiac chambers is regulated by the electrical conduction system of the heart whose main components are sketched in \cref{fig:heart_anatomy}, (c).
The pacemaking electrical signal arises in the \ac{SAN}, located in the \ac{RA} near the junction of the \ac{SupVC}.
From \ac{SAN} the signal propagates into the \ac{RA} myocardium and reaches the \ac{LA} through specific interatrial bundle connections, of which the most important is the \ac{BB}.
On the other side, the conduction network continues toward the \ac{AVN} where the signal is delayed until the end of the atrial contraction.
Then, the signal travels through the \ac{HB} and, in the interventricular septum, splits between the \ac{LBB} and the \ac{RBB} to end in the respective \ac{PFs} network located in the subendocardial layer.
Through the connection of the \ac{PFs} with the cardiomyocytes the signal transmurally propagates from the endocardium to the epicardium stimulating the ventricular contraction.
%\todoinline{Aggiungere qualcosa su polarizzazione/depolarizzazione e dinamica del calcio?}

%
\begin{figure}[t]
    \centering
    \includegraphics[width=0.9\textwidth]{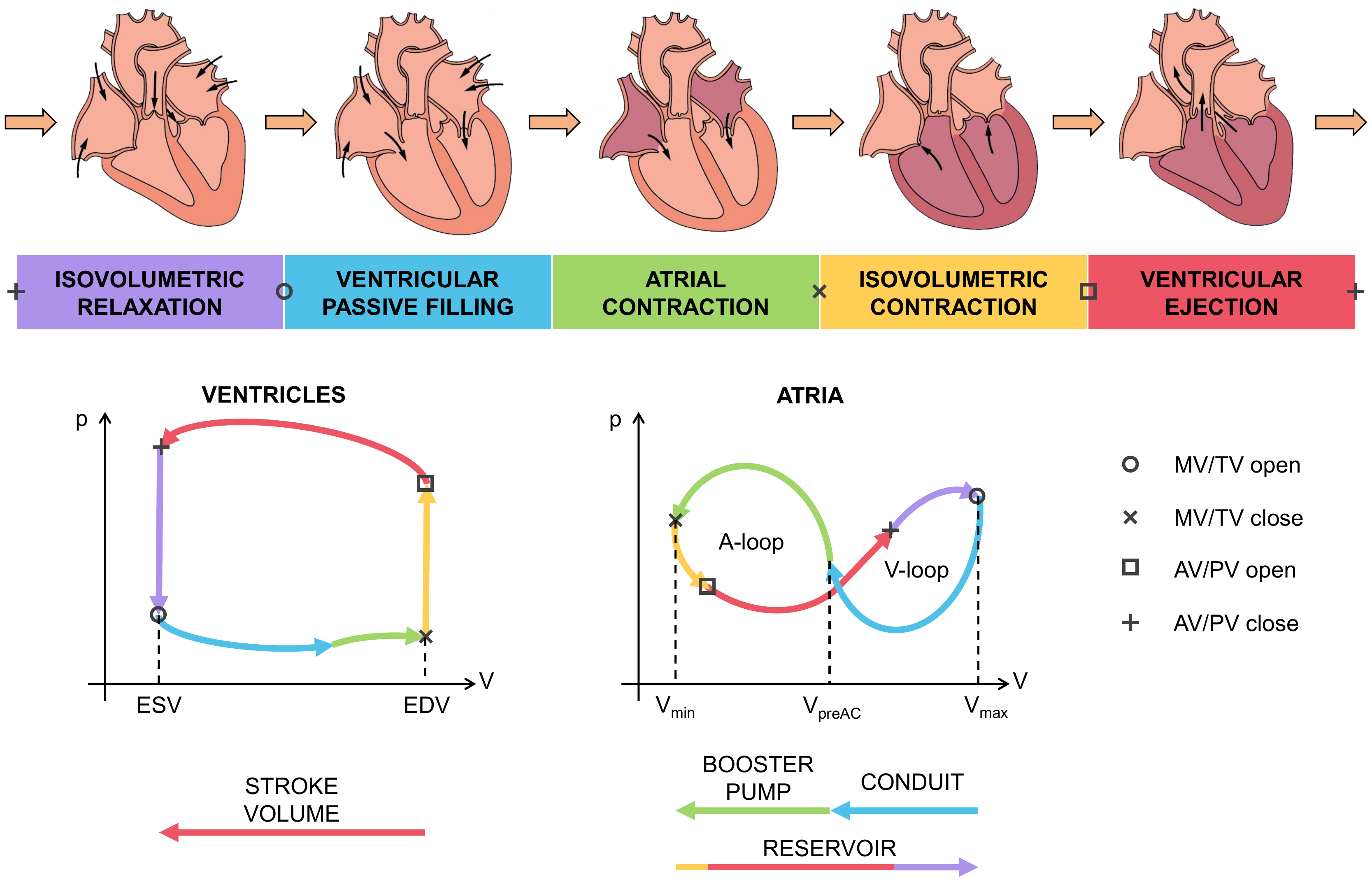}
    \caption{The five phases of the cardiac cycle:
    on the top, a sketch of the direction of the blood flow, the status of the valves and the contraction of the chambers (darker color) during the different phases
    (pictures elaborated from \url{https://commons.wikimedia.org/w/index.php?curid=30148227});
    on the bottom, schematic ventricular and atrial pressure-volume loops with the opening and closing of the valves and colored with the five phases.}
    \label{fig:cardiac_cycle}
\end{figure}

The main mechanical events of the cardiac cycle are sketched in \cref{fig:cardiac_cycle}, top, where the heart chambers are darker in color when contracting -- i.e. during \textit{systole} -- and lighter when relaxing -- i.e. during \textit{diastole}.
Atrial and ventricular systole and diastole occur in different phases of the cycle.
During the \ac{IVR} phase, the ventricular muscle is relaxing after the end of the ejection phase of the previous heartbeat.
Since all the valves are closed, the ventricular volume remains constant while its pressure quickly drops down until it reaches the atrial pressure, causing the opening of the atrioventricular valves (\ac{TV} and \ac{MV}).
At this moment the \ac{VPF} phase begins, blood flows from the atria to the ventricles and the volumes of both ventricles increase, driven by the muscle relaxation that thins the myocardium and moves the ventricular base upwards.
At the same time, the atrial volumes decrease due to the passive deformation of the atrial myocardium squeezed by the movement of the atrioventricular plane.
We remark that none of the four chambers is contracting during these first two phases.
When the passive filling slows down, the \ac{AC} begins, stimulated by the pacemaking of the \ac{SAN}.
The active deformation of the atrial muscle pushes additional blood toward the ventricles giving an additional ventricular preload;
for this peculiarity, this phase is also called \textit{atrial kick}.
The \ac{AVN} delays the electrical signal, allowing the ventricular contraction to begin only when the atrial contraction has ended.
When the ventricular muscle begins to contract, the ventricular pressure suddenly rises, exceeding the atrial pressure and determining the closure of the \ac{TV} and \ac{MV}.
This event starts the \ac{IVC} phase, in which all the cardiac valves are closed again.
This short phase ends when the pressures of the ventricles (\ac{RV} and \ac{LV}) reach the pressures of the respective arteries (\ac{PT} and \ac{AO}), triggering the opening of the semilunar valves (\ac{PV} and \ac{AV}) and the beginning of the \ac{VE} phase.
During this last phase, the ventricular volumes drop down -- driven by the myocardial thickening and the downward movement of the base -- and the blood flows toward the pulmonary and systemic circulation.
At the same time the atria fill, passively dilating due to the downward movement of the atrioventricular plane.

A plot of the pressure against the volume has long been used to measure the work done by a system and is also widely applied to assess the efficiency of the cardiac pump.
In \cref{fig:cardiac_cycle}, bottom, the different phases of the cardiac cycle are shown in this kind of diagram -- usually called \textit{pressure-volume loop} -- for both ventricles and atria.
The squared shape of the ventricular diagram is a direct consequence of the state of the valves and of the two isovolumetric phases.
To this counterclockwise loop is associated the (positive) ejection work exerted by the tissue on the blood \citep{Verzicco2022}, that increases if the maximum difference in pressures (or in volumes) arises.
The difference among the ventricular \ac{EDV} and \ac{ESV} is called \ac{SV} and represents the volume of blood pumped out from each ventricle during a heartbeat.
Instead, the interpretation of the eight-shaped atrial loop needs more explanations and can be divided into two parts.
%We can split the diagram into the following two loops.
The \textit{V-loop} is dominated by the effect of the ventricles on the atria.
Indeed, for a substantial part of the cardiac cycle, the atria fill or empty only passively, dragged by the contraction or relaxation of the ventricles.
Thus, the clockwise V-loop represents the (negative) work exerted by the ventricles on the atria, since the atria are not contracting during this period (i.e. their muscle is not consuming energy).
Conversely, the counterclockwise \textit{A-loop} is dominated by the atrial contraction and relaxation and is associated with the (positive) work exerted by the atrial muscle on the blood.
This complex pressure-volume loop is related to the threefold atrial function of \textit{reservoir}, \textit{conduit}, and \textit{booster pump} \cite{Spencer2001,Cui2008,Abhayaratna2008,Blume2011a,Rocsca2011,Hoit2014,Lang2015} (see \cref{fig:cardiac_cycle}, bottom-right):
while the atrioventricular valves (\ac{TV} and \ac{MV}) are closed, the atria store blood for later delivery to the ventricles (reservoir)
when the valves open, atria release blood to the ventricles, passively driven by the ventricular relaxation (conduit) \citep{Marino2021};
finally, at the end of ventricular diastole, the atrial muscle contraction actively supplies additional blood to the ventricles (booster pump), increasing the efficiency of the heart pump as is also evident from the effect on the ventricular pressure-volume loop.

This partial introduction to the heart anatomy and to the cardiac cycle aims at highlighting the complexity of concurring events that contribute to a physiological heart function and regulation.
In order to capture these events in a computational framework, an electromechanical model must accurately grasp the interaction among the heart chambers, the complex biophysics and all the multiscale and multiphysics aspects underlying the cardiac function.

%%%%%%%%
%%%%%%%%
\section{Mathematical models}
\label{sec:meth_math}

The mathematical model that we propose is based on the model by \citet{Regazzoni2022}, a cardiac electromechanical model fully coupled with a lumped-parameter model of blood circulation.
In that work, the only heart chamber considered as a 3D domain was the \ac{LV}.
This work was then extended by \citet{Piersanti2022biv} to a 3D domain for both ventricles.
However, in both cases, the 3D computational domain consists only of the ventricular muscle, which can be considered as a unique tissue with homogeneous electrical and mechanical properties.
Here, we propose the extension of these ventricular models to the whole heart, taking into account the heterogeneity of the cardiac tissue in the different cardiac components (\textit{e.g.} atria, ventricles, valves and vessels) for what concerns both electrical signal conduction and active/passive mechanics.

In the following sections we detail our mathematical model, focusing on the novelties introduced with respect to \citep{Regazzoni2022,Piersanti2022biv} to account for the extension to the entire human heart.
In \cref{subsec:domain} we define the whole-heart computational domain and we describe the choices made in terms of domain partitioning and boundaries;
in \cref{subsec:fibers} we discuss the modeling of the cardiac fibers;
in \cref{subsec:full_model} we present the full mathematical model by highlighting each core model;
in \cref{subsec:ref_conf_init_displ} we describe the strategy employed to recover the unloaded (i.e. stress-free) configuration and to subsequently compute the initial displacement.

%%%%%%%%%%%
\subsection{Computational domain}
\label{subsec:domain}
%%%
\begin{figure}[t]
    \centering
    \includegraphics[width=0.9\textwidth]{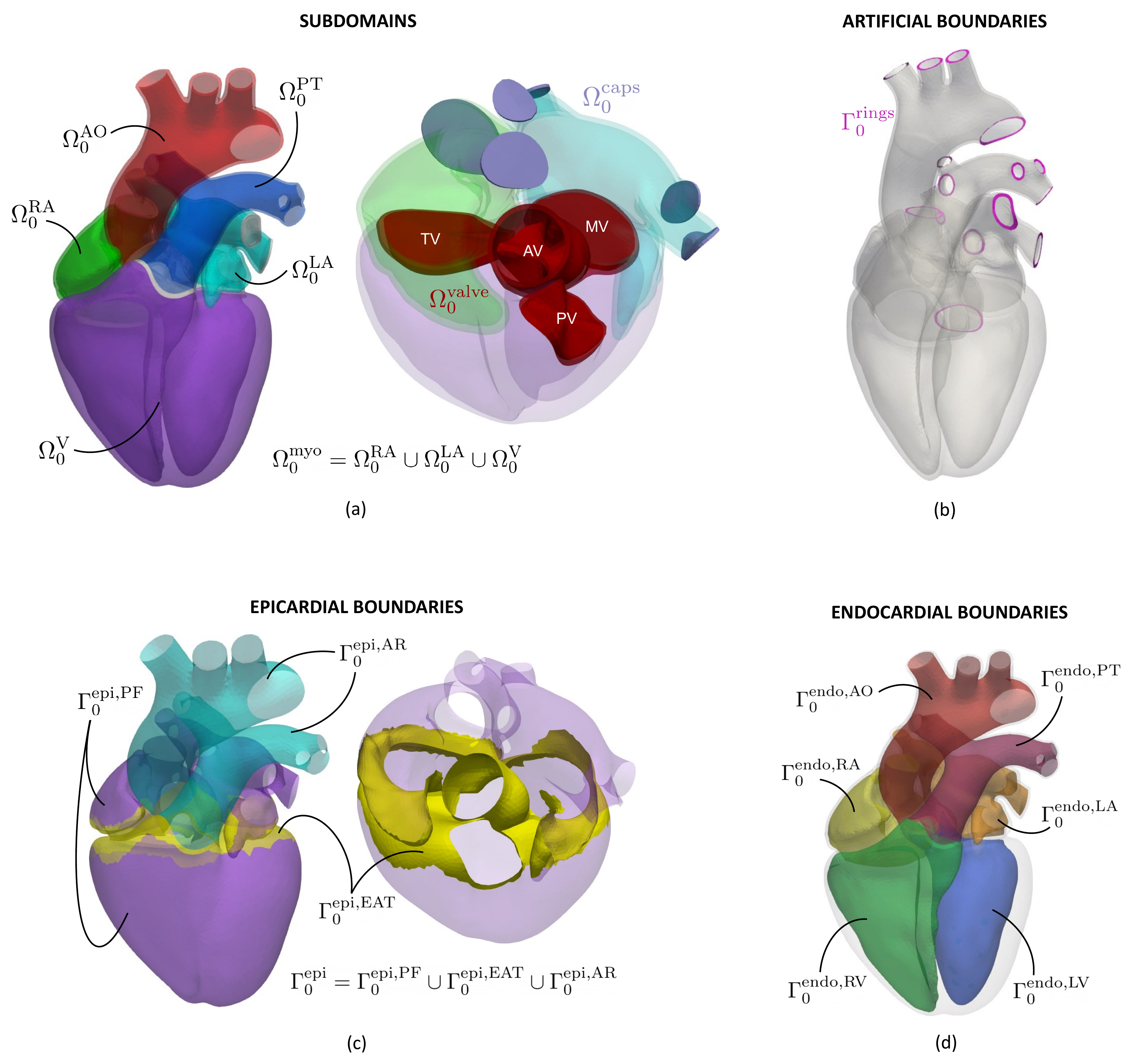}
    \caption{The computational domain $\Omega_0$:
    (a) the division in subdomains;
    (b) the artificial boundaries;
    (c) the epicardial boundaries;
    (d) the endocardial boundaries.}
    \label{fig:domain}
\end{figure}

% ANATOMY and SUBDOMAINs
In \cref{fig:domain}, (a), we show the computational domain $\Omega_0 \subset \mathbb{R}^3$ of the entire human heart, subdivided in the following subdomains:
\begin{itemize}
    \item the myocardium of the four cardiac chambers -- named $\OmegaMyo$ -- in turn divided into:
    (i) the \ac{RA} and \ac{LA} -- named $\OmegaRA$ and $\OmegaLA$, respectively -- characterized by a detailed anatomy that includes the two appendages (\ac{RAA}, \ac{LAA}) and physically connected to each other through the interatrial septum and the \ac{BB};
    (ii) a unique subdomain for the two ventricles -- named  $\OmegaVentr$ -- with a smoothed endocardium layer deprived of \ac{PaMs};
    \item the two arteries (\ac{PT}, \ac{AO}) -- named $\OmegaPT$ and $\OmegaAO$, respectively -- modeled up to their main bifurcations so that they can be fixed sufficiently far from their connection with the heart, where their movement can be considered negligible;
    \item the cardiac valves (\ac{TV}, \ac{MV}, \ac{PV}, \ac{AV}) -- named $\OmegaValve$ -- as flat simplified geometries filling the valvular orifices and connecting the atria to the ventricles (\ac{TV}, \ac{MV}) and the ventricles to the arteries (\ac{PV}, \ac{AV}).
    Although very simplified anatomycally, this representation allows to model some crucial aspects of valvular functioning such as the role of their \textit{annuli} as stiff and insulating fibrous tissue that connects different cardiac compartments and the high pressure difference across their closed leaflets occurring during some phases of the cardiac cycle;
    \item some artificial caps -- named $\OmegaCaps$ -- placed in all the entry veins (\ac{InfVC}, \ac{SupVC}, \ac{PVs}) and included in the domain in order to close the atrial blood pools, facilitating the calculation of their volumes
    (see \cref{subsubsec:circ_coupling}).
\end{itemize}
% BOUNDARIES
In order to apply proper boundary conditions to the mechanical model, the boundaries of the domain are divided as follows: 
\begin{itemize}
    \item some artificial boundaries (see \cref{fig:domain}, (b)) -- named $\GammaRings_0$ -- placed where the veins (\ac{InfVC}, \ac{SupVC}, \ac{PVs}) and the arteries (\ac{PT}, \ac{AO}) are cut;
    \item the external cardiac surface $\GammaEpi$ (see \cref{fig:domain}, (c)) in turn divided into:
    (i) the regions of the epicardium in contact with the \ac{PF}, named $\GammaEpiPF_0$;
    (ii) the regions of the epicardium in contact with the \ac{EAT}, named $\GammaEpiFat_0$;
    (iii) the epithelium of the two arteries, named $\GammaEpiArter_0$;
    % The dist 
    %
    \item the internal cardiac surface (see \cref{fig:domain}, (d)) made up of the endocardium of the four cardiac chambers (\ac{RA}, \ac{LA}, \ac{RV}, \ac{LV}) -- named $\GammaEndoRA_0$, $\GammaEndoLA_0$, $\GammaEndoRV_0$ and $\GammaEndoLV_0$, respectively -- and the endothelium of the two arteries (\ac{PT}, \ac{AO}) -- named $\GammaEndoPT_0$ and $\GammaEndoAo_0$, respectively.
\end{itemize}

% confronto con altri modelli 4CH, giustificare l'assenza del grasso e del pericardio "3D" citando le coorti di Strocchi \citep{Strocchi2020} per sottolineare come grasso e pericardio siano difficilmente ricostruibili da imaging;\\
% \citet{Gerach2021} pericardio 3D;\\
% \citet{Pfaller2019} grasso 3D;\\

% FIBERS
\subsection{Modeling the cardiac fibers}
\label{subsec:fibers}
\begin{figure}[t]
    \centering
    \includegraphics[width=0.9\textwidth]{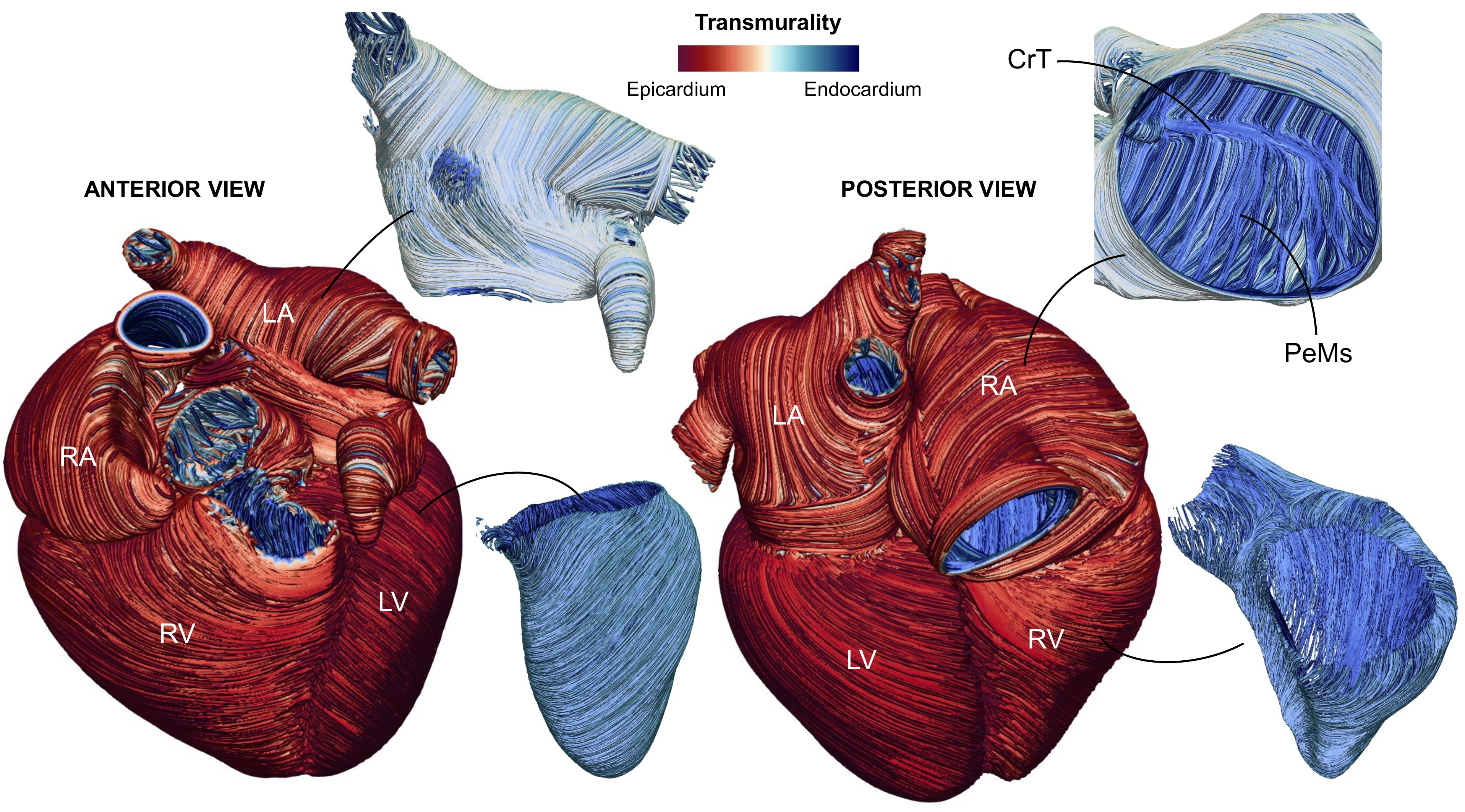}
    \caption{
        The myofibers architecture of the four cardiac chambers obtained using the whole-heart \acl{LDRBM} proposed by \citet{Piersanti2021fibers,Piersanti2021phd}.
        The transmural variation is pointed out through anterior and posterior views of the entire epicardium and chamber-specific views of the endocardium. On the endocardium of \ac{RA} also the \acf{CrT} and the \acf{PeMs} are clearly visible.}
    \label{fig:fibers}
\end{figure}
To prescribe the muscular fiber architecture in the myocardium $\OmegaMyo$, we rely on a particular class of \acfp{RBM}, known as \acp{LDRBM} \citep{bayer2012novel,Quarteroni2017,doste2019rule} recently reviewed in a communal mathematical description and also extended to account for atrial geometries in \citep{Piersanti2021fibers}.
Specifically, we use the whole-heart \ac{LDRBM} proposed by \citet{Piersanti2021fibers} in its improved version detailed in \citep[Chapter 4]{Piersanti2021phd}.

To properly reproduce the characteristic features of the cardiac fiber bundles in all the four chambers, the whole-heart \ac{LDRBM} first defines a transmural distance $\phi$ (from epicardium to endocardium) and several internal distances $\psi_i$
These are obtained by solving Laplace boundary-value problems of the type
\begin{equation}
    \label{eqn: laplace}
    \begin{cases}
        -\Delta \xi=0 &\qquad{\text{in }}\OmegaMyo,
        \\
        \xi = \xiA &\qquad{\text{on }}\GammaA,
        \\
        \xi = \xiB &\qquad{\text{on }}\GammaB,
        \\
        \nabla \xi \cdot \mecNref=0 &\qquad{\text{on }}\GammaN,
    \end{cases}
\end{equation}
where $\xiA,\,\xiB \in \mathbb{R}$ are suitable
Dirichlet data set on generic partitions of the heart boundary $\GammaA,\,\GammaB,\,\GammaN$, with $\GammaClosedA \cup \GammaClosedB \cup \GammaClosedN=\partial\OmegaMyo$.
In particular, the internal distances are used both to discriminate the left from the right heart and the atria from the ventricles, and also to represent different atrial and ventricular distances, characteristic of the four-chambers.
Then, for each point of the cardiac computational domain, the whole-heart \ac{LDRBM} suitably combines the gradients of the heart distances with the aim of defining an orthonormal local coordinate axial system $[\eL, \eN, \eT]$ owing to 
$\eT= \tfrac{\nabla \phi}{\left\lVert \nabla \phi \right\rVert}$, $\eN=\tfrac{\nabla \psi_i - (\nabla \psi_i \cdot \eT )\eT}{\left\lVert \nabla \psi_i - (\nabla \psi_i \cdot \eT )\eT \right\rVert}$ and $\eL=\eN \times \eT$, defined as the unit transmural, normal, and longitudinal directions, respectively.
Finally, the reference frame $[\eL, \eN, \eT]$ is properly rotated to define the myofiber orientations $[\eL, \eN, \eT] \xrightarrow{\alpha_j, \beta_j} [\f_0, \n_0, \s_0]$, where $\f_0$ is the fiber direction, $\n_0$ is the sheet-normal direction, $\s_0$ is the sheet direction, and $\alpha_j$ and $\beta_j$ are suitable helical and sheetlet angles following linear relationships $\theta_j (d_j) = \thetaEpii(1-d_j)+\thetaEndoi d_j$, (with $\theta_j=\alpha_j, \beta_j$) in which $d_j \in [0,1]$ is the transmural normalized distance and $\thetaEndoi$, $\thetaEpii$ are suitable prescribed rotation angles on the endocardium and epicardium of the $j$-th heart fibers bundle.

\cref{fig:fibers} shows that the whole-heart \ac{LDRBM} is able to accurately reproduce the myocardial fiber architecture, capturing the helical structure of \ac{LV}, the characteristic fibers of \ac{RV}, the outflow tracts regions and the fiber bundles of \ac{LA} and \ac{RA}, including the inter-atrial connections, the \acf{CrT} and the \acf{PeMs}. 
For further details about this whole-heart \ac{LDRBM} we refer to \citep{Piersanti2021phd}.

%%%%%%%%%%%
\subsection{The full electromechanical model}
\label{subsec:full_model}

%%%
\begin{figure}[t]
    \centering
    \includegraphics[width=1.0\textwidth]{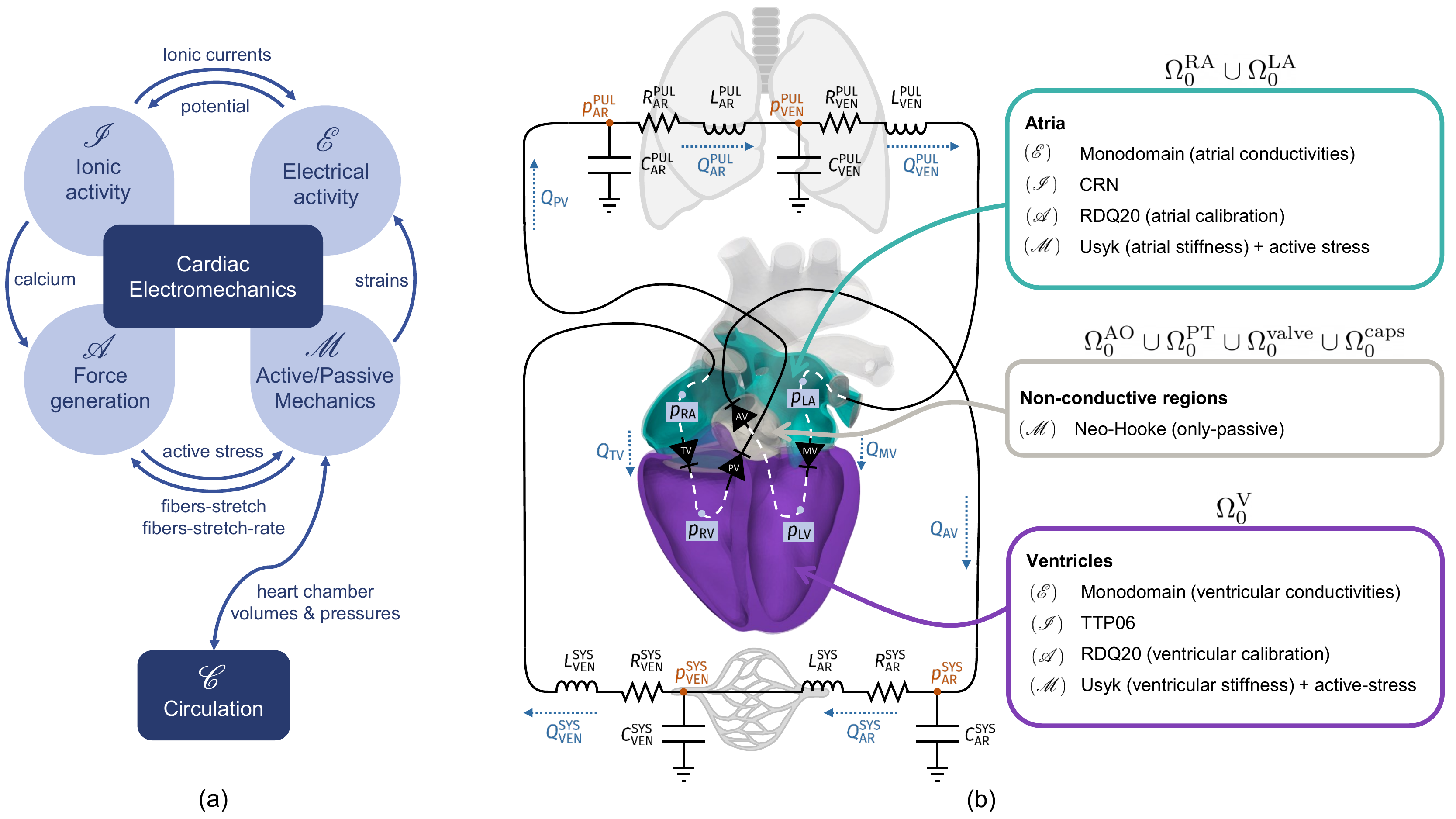}
    \caption{A sketch of the electromechanical model.
    (a) The underlying \textit{core models} and the fundamental quantities for their coupling.
    (b) The 0D model of the circulatory system made of resistance-inductance-capacitance (RLC) circuits for the systemic and pulmonary circulation and diodes for cardiac valves, coupled with the electromechanical model of the four 3D cardiac chambers.
    On the right, we highlight the models employed in the three main regions of the computational domain, i.e. atria, ventricles and non-conductive regions (\textit{CRN}, \citet{CRN}; \textit{TTP06}, \citet{TTP06}; \textit{RDQ20}, \citet{regazzoni2020biophysically}; \textit{Usyk}, \citet{Usyk2002}.}
    \label{fig:model}
\end{figure}
%%%
A multiphysics and multiscale whole-heart electromechanical model consists of several \textit{core models}, each of them describing biophysical processes that occur at different spatial and temporal scales during the cardiac cycle:
electrophysiology, in turn consisting of ionic activity ($\PhyIon$) at the microscale \citep{LuoRudy1,LuoRudy2,AlievPanfilov,CRN,TTP04,TTP06,BuenoOrovio,ToR-ORd} and electrical activity ($\PhyEP$) in terms of propagation of the transmembrane potential at the macroscale \citep{Henriquez,Pullan,Potse,ColliFranzone1,ColliFranzone2};
active force generation of cardiomyocites ($\PhyAct$) \citep{Rossi1,Rossi2,RuizBaier,regazzoni2020biophysically,Regazzoni2020};
active and passive mechanics of the cardiac tissue ($\PhyMec$) \citep{Ogden1997,Guccione1991a,Guccione1993a,Guccione1993b,HolzapfelOgden};
blood circulatory system ($\PhyCirc$) \citep{Blanco2010,Hirschvogel,Regazzoni2022}.
These core models are coupled to each other through some fundamental variables or feedbacks that represent biophysical processes.
In \cref{fig:model}, (a), we sketch the variables, interactions and feedbacks that we consider in our whole-heart electromechanical model:
the electrical and the ionic activities are coupled by the transmembrane potential and the ionic currents, respectively;
the ionic activity determines the calcium dynamics which is of fundamental importance for the active force generation model;
the cardiac mechanics is strongly influenced by the active stress provided by the force generation model which in turn is affected by the fibers-stretch and fibers-stretch-rate provided by the mechanics;
the loop is closed by the influence of the mechanical strains on the electrical activity;
finally, a volume conservation condition on the four cardiac chambers handles the two-way coupling between the 0D model of the circulatory system (see \cref{fig:model}, (b)) and the 3D cardiac mechanics.
More in detail, the proposed model features the following unknowns:
\begin{equation}
\label{eq:unknonws}
\begin{aligned}
&u \colon \OmegaMyo \times [0,T] \to \mathbb{R},
    && \wAtria \colon  \{\OmegaRA \cup \OmegaLA\} \times [0,T] \to \mathbb{R}^{n_{\wAtria}},\\
& \wVentr \colon  \OmegaVentr \times [0,T] \to \mathbb{R}^{n_{\wVentr}},
    &&\ActStateHF \colon \OmegaMyo \times [0,T] \to \mathbb{R}^{n_\ActStateHF},\\
& \Disp \colon \Omega_0 \times [0,T] \to \mathbb{R}^3,
    && \Circ \colon [0,T] \to \mathbb{R}^{n_\Circ}, \\
&p_i \colon [0,T] \to \mathbb{R},
    && \IinCH,
\end{aligned}
\end{equation}
where $u$ denotes the transmembrane potential, $\wAtria$ and $\wVentr$ the ionic variables on atria and ventricles, respectively, $\ActStateHF$ the state variables of the force generation model, $\Disp$ the mechanical displacement of the tissue, $\Circ$ the state vector of the circulation model (including pressures, volumes and fluxes in the different compartments of the vascular network), and $p_\RA$, $p_\LA$, $p_\RV$, and $p_\LV$ the blood pressures inside the four cardiac chambers.
% (\ac{RA}, \ac{LA}, \ac{RV}, \ac{LV}) and of the main arteries (\ac{PT}, \ac{AO}).
The full model reads as follows:
%%%%%%%%%%%%%%%%%%%%%%%
%%% FULL MODEL
%%%%%%%%%%%%%%%%%%%%%%%
\newlength{\EqWidth}
\newcommand{\MonodomainEqDiffTerm}{\nabla \cdot ( J \mecF^{-1} \EPdiffTens \mecF^{-T} \nabla u)}
\newcommand{\MonodomainEq}{
J \EPchim \left[ \EPCm \dfrac{\partial u}{\partial t} + \EPIion(u, \w_1, \w_2) \right] - \MonodomainEqDiffTerm = J \EPchim \EPIapp(t)
}
\newcommand{\CircEq}{\dfrac{d \Circ(t)}{d t} = \CircRhs(t,\,\Circ(t),\,p_\RA(t),\,p_\LA(t),\,p_\RV(t),\,p_\LV(t))\quad}
\settowidth{\EqWidth}{$\CircEq \quad$}
\newcommand{\EqBox}[1]{\makebox[\EqWidth][l]{$#1$}}

\newlength{\DomainWidth}
\settowidth{\DomainWidth}{$\text{in } \{\OmegaRA \cup \OmegaLA\} \times (0,T],$}
\newcommand{\DomainBox}[1]{\makebox[\DomainWidth][r]{$#1$}}

\newlength{\PhysWidth}
\settowidth{\PhysWidth}{$(\PhyAct)\quad$}
\newcommand{\PhysBox}[1]{\makebox[\PhysWidth][r]{$#1$}}

\newcommand{\IonicEq}[2]{\dfrac{ \partial #1 } {\partial t} - \EPrhsGating_{#2}(u, #1 ) = \boldsymbol{0}}
\newcommand{\MecEndoBC}[2]{\tenspiola(\Disp, #1) \, \mecNref = -{#2} \, J \mecF^{-T} \mecNref}
\newcommand{\CouplingEq}[1]{ V_{#1}^{\ThreeDim}(\Disp(t)) = V_{#1}^{\ZeroDim}(\Circ(t))}

%%%%%%%%%%%%%%%%%%%%%%
%%% EP - monodomain
%%%%%%%%%%%%%%%%%%%%%%
\begin{subequations}
\label{eq:ep}
    \begin{empheq}[ left=\PhysBox{ (\PhyEP)\empheqlbrace\, } ]{alignat=3}
        & \EqBox{ J \EPchim \left[ \EPCm \dfrac{\partial u}{\partial t} + \EPIion(u, \w_1, \w_2) \right] + }  &\nonumber \\
        & \EqBox{ \quad - \MonodomainEqDiffTerm = J \EPchim \EPIapp(t) }
            & { \text{in } \OmegaMyo \times (0,T], }
            \label{eq:ep_monodomain}
            \\
        & \left( J \mecF^{-1} \EPdiffTens \mecF^{-T} \nabla u \right) \cdot \mecNref = 0
            & \text{on } \partial\OmegaMyo \times (0,T], 
            \label{eq:ep_monodomain_bc}
    \end{empheq}
\end{subequations}
with $u = u_0$ in $\OmegaMyo$, at time $t=0$;
%%%%%%%%%%%%%%%%%%%%%%
%%% EP - ionics
%%%%%%%%%%%%%%%%%%%%%%
\begin{subequations}
\label{eq:ion}
    \begin{empheq}[ left=\PhysBox{ (\PhyIon)\empheqlbrace\, } ]{alignat=3}
        & { \IonicEq{\wAtria}{1} \qquad}
            & { \text{in } \{\OmegaRA \cup \OmegaLA\} \times (0,T], }
            \label{eq:ep_ionic_atria}
            \\
        & \IonicEq{\wVentr}{2}
            & \text{in } \OmegaVentr \times (0,T],
            \label{eq:ep_ionic_ventr}
    \end{empheq}
\end{subequations}
with $\w_1 = \w_{1,0}$ in $\{\OmegaRA \cup \OmegaLA\}$ and $\w_2 = \w_{2,0}$ in $\OmegaVentr$, at time $t=0$;
%%%%%%%%%%%%%%%%%%%%%%
%%% Activation
%%%%%%%%%%%%%%%%%%%%%%
\begin{empheq}[ left=\PhysBox{ (\PhyAct)\empheqlbrace\, } ]{alignat=3}
\label{eq:act}
    & \EqBox{ \dfrac{\partial \ActStateHF}{\partial t} = \ActRhs\left(\ActStateHF,\,\wCa,\,\SL, \dfrac{\partial\SL}{\partial t}\right) }
        & {\text{in } \OmegaMyo \times (0,T],}\quad
\end{empheq}
with $\ActStateHF = \ActStateHF_0$ in $\OmegaMyo$ at time $t=0$;
%%%%%%%%%%%%%%%%%%%%%%
%%% Mechanics
%%%%%%%%%%%%%%%%%%%%%%
\begin{subequations}
\label{eq:mec}
    \begin{empheq}[ left=\PhysBox{ (\PhyMec)\empheqlbrace\, } ]{alignat=3}
        & { \rho_{\mathrm{s}} \dfrac{\partial^2 \Disp}{\partial t^2} - \nabla \cdot \tenspiola(\Disp, \Tens(\ActStateHF,\SL)) = \boldsymbol{0} }
            & { \text{in } \Omega_0 \times (0,T], }
            \label{eq:mec_pde}
            \\
        & \tenspiola(\Disp, \Tens(\ActStateHF,\SL)) \, \mecNref + & \nonumber\\
        & \quad +(\mecNref \otimes \mecNref) \left( \BCmecKepiN \Disp + \BCmecCepiN \dfrac{\partial \Disp}{\partial t} \right) = \mathbf{0}\qquad
            & \text{on } \GammaEpi_0 \times (0,T],
            \label{eq:mec_bc_epi}
            \\
        & \MecEndoBC{\Tens(\ActStateHF,\SL)}{\PRA}
            & \text{on } \GammaEndoRA_0 \times (0,T],
            \label{eq:mec_bc_endo_RA}
            \\
        & \MecEndoBC{\Tens(\ActStateHF,\SL)}{\PLA}
            & \text{on } \GammaEndoLA_0 \times (0,T],
            \label{eq:mec_bc_endo_LA}
            \\
        & \MecEndoBC{\Tens(\ActStateHF,\SL)}{\PRV}
            & \text{on } \GammaEndoRV_0 \times (0,T],
            \label{eq:mec_bc_endo_RV}
            \\
        & \MecEndoBC{\Tens(\ActStateHF,\SL)}{\PLV}
            & \text{on } \GammaEndoLV_0 \times (0,T],
            \label{eq:mec_bc_endo_LV}
            \\
        & \MecEndoBC{\Tens(\ActStateHF,\SL)}{\PPT}
            & \text{on } \GammaEndoPT_0 \times (0,T],
            \label{eq:mec_bc_endo_PT}
            \\
        & \MecEndoBC{\Tens(\ActStateHF,\SL)}{\PAO}
            & \text{on } \GammaEndoAo_0 \times (0,T],
            \label{eq:mec_bc_endo_Ao}
            \\
        & \Disp = \boldsymbol{0}
            & \text{on } \GammaRings_0 \times (0,T],
            \label{eq:mec_bc_rings}
    \end{empheq}
\end{subequations}
with $\Disp = \Disp_0$ and $\dfrac{\partial \Disp}{\partial t} = \Dot{\Disp}_0$ in $\Omega_0$ at time $t=0$;
%%%%%%%%%%%%%%%%%%%%%%
%%% Circulation
%%%%%%%%%%%%%%%%%%%%%%
\begin{empheq}[ left=\PhysBox{ (\PhyCirc)\empheqlbrace\, } ]{alignat=3}
\label{eq:circ}
    & \EqBox{ \dfrac{d \Circ(t)}{d t} = \CircRhs(t,\,\Circ(t),\,p_\RA(t),\,p_\LA(t),\,p_\RV(t),\,p_\LV(t)) }
        & { \text{for } t \in (0,T], }\quad
\end{empheq}
with $\Circ(0) = \Circ_0$ at time $t=0$;
%%%%%%%%%%%%%%%%%%%%%%
%%% 3D-0D coupling
%%%%%%%%%%%%%%%%%%%%%%
\begin{subequations}
\label{eq:coupl}
    \begin{empheq}[ left=\PhysBox{ (\PhyCoupl)\empheqlbrace\, } ]{alignat=3}
        & \EqBox{ \CouplingEq{\RA} }
            & { \text{for } t \in (0,T], } 
            \label{eq:coupling_RA}
            \\
        & \CouplingEq{\LA}
            & \text{for } t \in (0,T],
            \label{eq:coupling_LA}
            \\
        & \CouplingEq{\RV}
            & \text{for } t \in (0,T],
            \label{eq:coupling_RV}
            \\
        & \CouplingEq{\LV}
            & \text{for } t \in (0,T].
            \label{eq:coupling_LV}
    \end{empheq}
\end{subequations}
We remark that $(\PhyEP)$ and $(\PhyAct)$ are both defined in the whole domain $\OmegaMyo$, but with specific parameters for atria and ventricles. 
Thus, since the myocardial domain is composed of the two disconnected parts $\{\OmegaRA \cup \OmegaLA\}$ and $\OmegaVentr$, they behave independently in the atria and ventricles.
Instead, $(\PhyIon)$ is composed of two distinct ionic models, each one characterized by different variables and equations for the atria and ventricles.
The variability of the parameters and models employed in the different regions of the heart is sketched in \cref{fig:model}, (b).

In \crefrange{subsubsec:ep}{subsubsec:circ_coupling} we describe each core model (\crefrange{eq:ep}{eq:coupl}), detailing how they are coupled to each other and how they vary along the heart domain.

%%%
\subsubsection[Electrophysiology]{Electrophysiology ($\PhyEP$)-($\PhyIon$)}
\label{subsubsec:ep}
% Domain
\cref{eq:ep,eq:ion} represent the electrophysiological model and are solved only in the cardiac chambers $\OmegaMyo$, since the rest of the domain is made up of non-conductive regions.
% Model description
The model consists of the monodomain equation \eqref{eq:ep} -- describing the propagation of the transmembrane potential $u$ \citep{Pullan,ColliFranzone2} -- coupled with suitable ionic models, one for the atria in $\{\OmegaRA \cup \OmegaLA\}$ (\cref{eq:ep_ionic_atria}) and one for the ventricles in $\OmegaVentr$ (\cref{eq:ep_ionic_ventr}).
The vectors $\wAtria = \{w_{1,i}\}_{i=1}^{n_{\wAtria}}$ and $\wVentr = \{w_{2,j}\}_{j=1}^{n_{\wVentr}}$ constitute the gating variables and the concentrations of ionic species.
Among them, the intracellular calcium ions concentration $\Cai$ plays a crucial role for active force generation.
We denote this quantity with $\wCaAtria$ and $\wCaVentr$ for the atria and ventricles, respectively, and we define in the whole myocardium the corresponding function $\wCa$ (used in \cref{eq:act}) as:
\begin{align}
\label{eq:w_ca}
\wCa &=
    \begin{cases}
        \wCaAtria,
            & \qquad \text{in } \{\OmegaRA \cup \OmegaLA\} \times [0,T],\\
        \wCaVentr,
            & \qquad \text{in } \OmegaVentr \times [0,T].
    \end{cases}
\end{align}
% Ionic models and monodomain-I_ion
We use the \ac{CRN} model for the atria and the \ac{TTP06} model for the ventricles, respectively.
These two models are used to define the nonlinear reaction term $\EPIion(u,\wAtria,\wVentr)$ of \cref{eq:ep_monodomain} that models the ionic currents taking into account the multiscale effects from the cellular to the tissue level:
\begin{align}
\label{eq:Iion}
\EPIion(u, \w_1, \w_2) &=
    \begin{cases}
        \EPIion(u, \w_1),
            & \qquad \text{in } \{\OmegaRA \cup \OmegaLA\} \times [0,T],\\
        \EPIion(u, \w_2),
            & \qquad \text{in } \OmegaVentr \times [0,T].
    \end{cases}
\end{align}
% BC and A-V insulation
The monodomain model is finally closed by the no-flux Neumann boundary condition of \cref{eq:ep_monodomain_bc} that represents an electrically insulated domain.
Moreover, since the domain $\OmegaMyo$ is composed of the two disjoint parts ($\OmegaRA \cup \OmegaLA$ and $\OmegaVentr$) separated by the insulating fibrous tissue of the atrioventricular valves (the \ac{TV} and \ac{MV} parts of $\OmegaValve$), also the atrial and ventricular muscles are electrically insulated from each other.

% Monodomain: diffusion term, MEF, diffusion tensor, sigma and fast endocardial layer.
The anisotropic transmission of the transmembrane potential $u$ is regulated by the diffusion term $\MonodomainEqDiffTerm$ of the monodomain model (\cref{eq:ep_monodomain}) \citep{Regazzoni2022}, where $\EPdiffTens$ represents the diffusion tensor in the deformed configuration and $\mecF = \identity + \nabla \Disp$ and $J = \det(\mecF)$ are the deformation gradient tensor and the deformation Jacobian, respectively.
Being $\mecF$ dependent on the unknown displacement $\Disp$ of the mechanical model ($\PhyMec$) (\cref{eq:mec}), this diffusion term takes into account the variation of the electrical properties due to the tissue deformation, modeling the so-called \acf{MEF} \citep{Kohl2003}.
The diffusion tensor $\EPdiffTens$ regulates the anisotropic conduction of the electrical signal using the local orthonormal coordinates system $(\f_0,\s_0,\n_0)$ (see \cref{subsec:fibers}) by prescribing three different conductivities $\CondF{*}$, $\CondS{*}$, and $\CondN{*}$ along the fiber, sheet normal and crossfiber directions, respectively~\citep{Regazzoni2022,Piersanti2022biv}:
\begin{equation}
\label{eq:ep_diff_tens}
\EPdiffTens = \CondF{*} \frac{\mecF\f_0 \otimes \mecF\f_0}{\|\mecF \f_0 \|^2} + \CondS{*} \frac{\mecF\s_0 \otimes \mecF\s_0}{\|\mecF \s_0 \|^2} + \CondN{*} \frac{\mecF\n_0 \otimes \mecF\n_0}{\|\mecF \n_0 \|^2}.
\end{equation}
In addition to varying along the local direction, the conductivities vary in space depending on the cardiac compartment:
\begin{align}
\label{eq:ep_cond}
    \CondK{*} &=
        \begin{cases}
            \CondK{\Atria} 
                &\text{ in } \{\OmegaRA \cup \OmegaLA\}, \\[0.5em]
            \CondK{\Ventr}(\phi) =
            \begin{cases}
            \CondK{\Ventr,\Myo} \qquad \text{ if } \phi > \epsilon , \\[0.5em]
            \CondK{\Ventr,\Endo} \qquad \text{ if } \phi \leq \epsilon ,
            \end{cases}
                & \text{ in } \OmegaVentr,
        \end{cases}
    & \text{ for } \mathbf{k} = \f, \s, \n.
\end{align}
Following \citet{Piersanti2022biv}, the conductivities in the ventricles $\OmegaVentr$ also depend on a scalar function $\phi$ that smoothly connects the endocardium to the epicardium, allowing the definition of an endocardial layer where the electric signal propagates faster.
This surrogates the \ac{PFs} network \citep{Lee2019,DelCorso2022ep4ch} and represents a valid alternative (at least in sinus rhythm) to the generation of the \ac{PFs} as a 1D network \citep{Vergara2014,Vergara2016,Costabal2016,Landajuela2018}.
Instead, in the atria $\{\OmegaRA \cup \OmegaLA\}$, different conduction velocities of the various bundles characterize the atrial fibers morphology (see \cref{subsec:fibers}), varying from fast to slow conduction regions \citep{ferrer2015detailed}.
This feature is of paramount importance in the modeling of atrial electrical disorders and related pathologies \citep{dossel2012computational,lemery2007normal}.
However, as this work is focused on a healthy scenario, we do not vary the conductivities $\CondF{\Atria}$, $\CondS{\Atria}$, and $\CondN{\Atria}$ in space, considering only the variation along the local fibers orientation, as done in \citep{Piersanti2021fibers}.

% Monodomain: Iapp
Finally, the forcing term $\EPIapp(t)$ of \cref{eq:ep_monodomain} represents an applied current that triggers the action potential of the myocardium at specific locations and times.
This term is used to model a series of electrical impulses that mimic the behavior of the electrical conduction system (see \cref{sec:heart,fig:heart_anatomy}, (c)), starting from the \ac{SAN} and ending into a series of points on the ventricular endocardium which, combined with the fast endocardial layer, surrogate the effect of the \ac{PFs}.

%%%%%%%%%%%
\subsubsection[Active force generation]{Active force generation ($\PhyAct$)}
\label{subsubsec:active_force}

We model the subcellular processes by which cardiomyocytes generate an active force in response to changes in calcium concentration $\wCa$ using the model proposed by \ac{RDQ20}.
\ac{RDQ20} is based on a biophysically accurate description of the subcellular mechanisms of force generation and regulation.
Despite its computational lightness (its state $\ActStateHF$ has only 20 variables), this model explicitly describes the end-to-end interactions of tropomyosin, which are responsible for the cooperative tissue response to calcium ion concentration, manifested in a markedly enhanced sensitivity to calcium around the half-maximal effective concentration (so-called $\text{EC}_{50}$).
Moreover, the \ac{RDQ20} model takes into account the effect of sarcomere length on the total force generated and, thanks to its explicit representation of the attachment-detachment mechanism of crossbridges, it is able to reproduce the force-velocity relationship, according to which the generated force decreases while the muscle fibers are shortening.
These subcellular mechanisms are responsible for two organ-level feedbacks, namely the fibers-stretch and the fibers-stretch-rate feedback, which regulate the force generated in each region of the myocardium depending on how much and how quickly it deforms \cite{regazzoni2019reviewXB}. 
The former is related to the dependence of the model \eqref{eq:act} on $\SL$, while the latter is related to the dependence of the model on ${\partial\SL}/{\partial t}$.
The variable $\SL$ represents the local sarcomere length, obtained as $\SL = \SL_0 \|\mecF \f_0 \|$, where $\SL_0$ is the sarcomere length at rest.
The regulatory and feedback mechanisms mentioned above play a key role in the cardiac function. 
Nevertheless, some of them are sometimes neglected in multiscale models, % @FR CITE?
due to the difficulty of capturing them in mathematical models of low computational cost and because of the difficulties involved in their numerical approximation.

The \ac{RDQ20} model describes subcellular mechanisms inherent to both atrial and ventricular cells.
The model can be adapted to reproduce experimental measurements of different cell types by calibrating the parameters, which reflect the different calcium-sensitivity and kinetics of protein interactions. 
See \cite{regazzoni2022machine} for an adaptation to ventricular cells and \cite{Mazhar2021} for a calibration to atrial cells.
Therefore, we use the same model throughout the computational domain, but with different parameter calibration to reflect the specificities of the cells belonging to the different chambers. 

The tissue level active tension $\Tens$ of the \ac{RDQ20} can be defined as a nonlinear function of the state $\ActStateHF$ and of the sarcomere length $\SL$ \citep{regazzoni2020biophysically}.
% TODO: AQ: dovremmo dire cosa rappresenta ActStateHF, nelle due eq. precedenti s veniva usato per gli sheet.
This quantity determines the coupling with the mechanical model \eqref{eq:mec} and contributes to the active stress part of the Piola-Kirchhoff stress tensor (see \cref{subsubsec:mechanics}).
More specifically, $\Tens$ can be written as:
\begin{equation}
  \Tens(\ActStateHF,\SL)=\aXB^i G(\ActStateHF,\SL), \qquad \text{for }\IinCH,
\label{eq:active_tension}
\end{equation}
where the microscale crossbridge stiffness $\aXB^i$ links the microscopic force with the macroscopic active tension and $G(\ActStateHF,\SL)$ is a nonlinear function (see \citep{regazzoni2020biophysically}).
Thus, the organ-level contractility of each chamber is calibrated using the $\aXB^i$ parameter.
Moreover, in order to set a specific contractility also in the \ac{RV} and \ac{LV} (that belong to the same subdomain $\OmegaVentr$), we use the same strategy proposed by \citet{Piersanti2022biv} (for the previous version of the active force generation model \citep{regazzoni2018redODE}) defining the ventricular microscale crossbridge stiffness $\aXB^\mathrm{V}\colon \OmegaVentr \to \mathbb{R}$ as a function of space:
\begin{equation}
\aXB^\mathrm{V}(\x) = \aXB^\mathrm{LV} \,\left( \XiHat(\x) + \Clrv ( 1 - \XiHat(\x) ) \right),
\label{eq:a_XB_ventr}
\end{equation}
where $\XiHat \colon \OmegaVentr \to [0,1]$ is the normalized interventricular distance \citep{Piersanti2021fibers,Piersanti2022biv} -- that smoothly goes from $0$ to $1$ in the interventricular septum -- and $\Clrv \in \mathbb{R}$ is a coefficient that represents the left-right ventricle contractility ratio.
In practical terms, this is equivalent to setting two constant values in the two ventricles (smoothly connected in the septum):
$\aXB^\mathrm{LV}$ in the \ac{LV} and $\aXB^\mathrm{RV}=\Clrv\,\aXB^\mathrm{LV}$ in the \ac{RV}.

%%%%%%%%%%%
\subsubsection[Active and passive mechanics]{Active and passive mechanics ($\PhyMec$)}
\label{subsubsec:mechanics}

The mechanics of the cardiac tissue is modeled by the problem $(\PhyMec)$ of \cref{eq:mec},
describing the dynamics of the tissue displacement $\Disp$ by the momentum conservation (\cref{eq:mec_pde}) under the hyperelasticity assumption \citep{Ogden1997} and employing an active stress approach \citep{Guccione1991b}.
The active and passive mechanical properties are embedded in the Piola-Kirchhoff stress tensor $\tenspiola(\Disp, \Tens(\ActStateHF,\SL))$:
\begin{subequations}
\label{eq:piola}
    \begin{empheq}[]{alignat=3}
        &\tenspiola(\Disp, \Tens(\ActStateHF,\SL)) = \dfrac{\partial \mathcal{W}(\mecF)}{\partial \mecF} + \nonumber\\
        &\quad + \Tens(\ActStateHF,\SL) \left[ n_\f \frac{\mecF \f_0 \otimes \f_0}{ \sqrt{\IIVf}} + n_\s \frac{\mecF \s_0 \otimes \s_0}{ \sqrt{\IIVs}} + n_\n \frac{\mecF \n_0 \otimes \n_0}{ \sqrt{\IIVn}}\right]\quad
            & \text{in } \OmegaMyo ,
            \label{eq:piola_myo}
            \\
        &\tenspiola(\Disp, \Tens(\ActStateHF,\SL)) = \dfrac{\partial \mathcal{W}(\mecF)}{\partial \mecF}
            & \text{in } \left\{ \Omega_0 \setminus \OmegaMyo \right\}.
            \label{eq:piola_nocond}
    \end{empheq}
\end{subequations}

The passive part of the tensor is modeled by the term $\partial \mathcal{W}(\mecF) / \partial \mecF$ where $\mathcal{W}$ is the hyperelastic strain energy density function.
In the myocardium $\OmegaMyo$ we employ the exponential constitutive law of \citet{Usyk2002}, with a volumetric term enforcing quasi-incompressibility \citep{Cheng2005,Doll2000,Yin1996,Regazzoni2022}.
In the non-conductive regions $\{\Omega_0 \setminus \OmegaMyo\}$, instead, we use a Neo-Hookean model \citep{Ogden1997}.
The resulting strain energy density function reads:
\begin{subequations}
\label{eq:w_all}
\begin{empheq}[left={\mathcal{W}(\mecF) =\,\empheqlbrace\,}]{align}
  & \dfrac{C^i}{2} \left( e^Q  - 1 \right) + \dfrac{B}{2} \left( J - 1 \right) \log(J),
    & \text{in } \OmegaMyo\label{eq:w_guccione},\\
  &\frac{\mu^j}{2}\left(J^{-\frac{2}{3}}\mecF:\mecF-3\right) + \dfrac{\kappa^j}{4}\left[ \left( J - 1 \right)^2 + \log^2(J)\right], 
    & \text{in } \left\{\Omega_0 \setminus \OmegaMyo\right\}, \label{eq:w_neo_hooke}
\end{empheq}
\end{subequations}
where, in the \citet{Usyk2002} model \eqref{eq:w_guccione},
$B \in \mathbb{R}^+$ represents the bulk modulus contributing to the term that realizes a weakly incompressible constraint \citep{Regazzoni2022},
$C^i$, for $\IinCHV$, is the stiffness scaling parameter that assumes a specific value in each subdomain of the myocardium $\OmegaRA$, $\OmegaLA$, and $\OmegaVentr$.
Instead, in the Neo-Hookean model \eqref{eq:w_neo_hooke},
$\mu^j$ and $\kappa^j$, for $j \in \{\mathrm{valve},\mathrm{caps},\AO,\PT\}$, are the shear modulus and the bulk modulus, respectively, and assume specific values in each non-conductive region $\OmegaValve$, $\OmegaCaps$, $\OmegaAO$, and $\OmegaPT$.
Finally, the term $Q$ of the \citet{Usyk2002} model \eqref{eq:w_guccione} reads:
\begin{equation*}
\begin{aligned}
&Q =
  b_{\mathrm{ff}} E_{\mathrm{ff}}^2  + b_{\mathrm{ss}} E_{\mathrm{ss}}^2 + b_{\mathrm{nn}} E_{\mathrm{nn}}^2+ b_{\mathrm{fs}} \left( E_{\mathrm{fs}}^2 + E_{\mathrm{sf}}^2 \right) + b_{\mathrm{fn}} \left( E_{\mathrm{fn}}^2 + E_{\mathrm{nf}}^2 \right) + b_{\mathrm{sn}} \left( E_{\mathrm{sn}}^2 + E_{\mathrm{ns}}^2 \right),\\
&E_\text{ab} =
  \textbf{E} \boldsymbol{a}_\text{0} \cdot \boldsymbol{b}_\text{0}, \qquad \text{for }a, b \in \{ f, s, n \},
\end{aligned}
\end{equation*}
where $\textbf{E} = \tfrac{1}{2} \left( \mecC - \identity \right)$ is the Green-Lagrange strain energy tensor, being $\mecC = \mecF^{T} \mecF$ the right Cauchy-Green deformation tensor.
%We refer to \ref{sec:appendix_params} for the specific calibration of the model parameters in each subdomain.

The active part of the Piola-Kirchhoff stress tensor acts only in the conductive subdomains $\OmegaMyo$.
This tensor depends on the active tension $\Tens(\ActStateHF,\SL)$, provided by the active force generation model \eqref{eq:act}, and on the fiber orientation in the deformed configuration.
%More specifically, following \citet{Piersanti2022biv}, 
We consider the orthotropic active stress tensor \eqref{eq:piola_myo}, where the coefficients $\IIVf$, $\IIVs$, and $\IIVn$ (equal to $\mecF \mathbf{k} \cdot \mecF \mathbf{k}$, for $\mathbf{k}=\f_0,\s_0,\n_0$) represent the tissue stretches along the fiber, sheet, and sheet-normal directions, respectively,
while $n_\f$, $n_\s$, and $n_\n$ model the proportion of active tension along these directions~\citep{Piersanti2022biv}.
In this way the active stress tensor can mainly act in the fiber direction $\f$ while also being applied on the cross-fiber directions $\s$ and $\n$, to surrogate the contraction caused by the dispersed myofibers \citep{Guan2020,Guan2021}.

The mechanical model is closed by the boundary conditions of \crefrange{eq:mec_bc_epi}{eq:mec_bc_rings}.
On the epicardium $\GammaEpi_0$ we apply the Robin-like condition \eqref{eq:mec_bc_epi} originally proposed in the whole-heart context by \citet{Pfaller2019}.
This condition surrogates the pressure exerted by the pericardium and surrounding organs on the external cardiac surface by penalizing only the normal displacement \citep{Pfaller2019,Strocchi2020b}.
No constraints are added on the other directions as the pericardial fluid allows free sliding within the pericardial sac \citep{Pfaller2019,Strocchi2020b}.
Instead, an additional constraint on the tangential direction can be necessary to avoid rigid rotation when the computational domain consists of the sole ventricles \citep{Regazzoni2022}.
The calibration of the pericardial stiffness $\BCmecKepiN$ of \cref{eq:mec_bc_epi} plays a fundamental role in the realistic movement of the heart \citep{Pfaller2019,Strocchi2020b}.
\citet{Pfaller2019} have tested different constant values on the whole external cardiac surface, but they model the \ac{EAT} as a 3D subdomain.
More recently, \citet{Strocchi2020b} have proposed a spatially varying coefficient to surrogate the different stiffness of the organs in contact with the pericardial sac, without including the \ac{EAT} as a 3D subdomain.
% Both these works exploit medical image data to validate the used coefficients.
Inspired by both of these works, we vary the $\BCmecKepiN$ only between two regions:
we prescribe a stiffer value on $\GammaEpiPF_0$ -- where the external organs are in contact with the pericardium -- and a much lower value $\BCmecKepiNfat$ on $\GammaEpiFat_0$ -- where the presence of the \ac{EAT} leaves the ventricular base and the lower part of the \ac{LAA} and \ac{RAA} more free to move.

On the endocardium and endothelium surfaces, we apply the normal stress boundary conditions of \crefrange{eq:mec_bc_endo_RA}{eq:mec_bc_endo_Ao} that model the pressure exerted by the blood.
The blood pressure of the various chambers and arteries depends on the circulation model \eqref{eq:circ}, as detailed in \cref{subsubsec:circ_coupling}.
Finally, we apply the homogeneous Dirichlet boundary condition \eqref{eq:mec_bc_rings} on all the artificial boundaries $\GammaRings_0$, since the arteries and veins can be considered almost fixed where we cut the computational domain (see \cref{fig:domain}).

%%%%%%%%%%%
\subsubsection[Blood circulation and 3D-0D coupling]{Blood circulation ($\PhyCirc$) and 3D-0D coupling ($\PhyCoupl$)}
\label{subsubsec:circ_coupling}
We model the blood circulatory system using the 0D lumped-parameter closed-loop model proposed by \citet{Regazzoni2022} and inspired by \citet{Blanco2010,Hirschvogel}.
In this model, as sketched in \cref{fig:model}, (b), resistance-inductance-capacitance (RLC) circuits represent the systemic ($\mathrm{SYS}$) and pulmonary ($\mathrm{PUL}$) circulations in both their arterial ($\mathrm{AR}$) and venous ($\mathrm{VEN}$) compartments, while non-ideal diodes model the four cardiac valves.
The state vector $\Circ$ comprises the volumes of the cardiac chambers and the systemic/pulmonary arterial/venous pressures and flow rates:
\begin{equation*}
    \begin{aligned}
    \Circ(t) = \big(&\VRA, \VLA, \VRV, \VLV,\\
     &\ParSYS, \PvnSYS, \ParPUL, \PvnPUL,\\
     &\QarSYS, \QvnSYS, \QarPUL, \QvnPUL \big).
    \end{aligned}
\end{equation*}
The corresponding ODE system ($\PhyCirc$), summarized by \cref{eq:circ}, reads:
\begin{equation}
\label{eq:circ_extended}
\left\{
\begin{aligned}
    &\CvnSYS \dfrac{d \PvnSYS}{d t} = \QarSYS - \QvnSYS, \\ % \label{eq:circ_PvnSYS}, \\
    &\CvnPUL \dfrac{d \PvnPUL}{d t} = \QarPUL - \QvnPUL, \\ % \label{eq:circ_PvnPUL}, \\
    &\dfrac{\LvnSYS}{\RvnSYS} \dfrac{d \QvnSYS}{d t}  = - \QvnSYS - \dfrac{\PRA    - \PvnSYS}{\RvnSYS}, \\ % \label{eq:circ_QvnSYS}\\
    &\dfrac{\LvnPUL}{\RvnPUL} \dfrac{d \QvnPUL}{d t}  = - \QvnPUL - \dfrac{\PLA    - \PvnPUL}{\RvnPUL}, \\ % \label{eq:circ_QvnPUL}
    %%%%%%%%%%%%%%%%%%%%%%%%%%%%%%%%%%%%%%%%%%%%%%%%%%
    &\dfrac{d \VRA}{d t}  = \QvnSYS - \QTV(\PRA,\PRV),  \\ % \label{eq:circ_VRA},
    &\dfrac{d \VLA}{d t}  = \QvnPUL - \QMV(\PLA,\PLV),  \\ % \label{eq:circ_VLA},
    &\dfrac{d \VRV}{d t}  = \QTV(\PRA,\PRV) - \QPV(\PRV, \ParPUL),  \\ % \label{eq:circ_VRV},
    &\dfrac{d \VLV}{d t}  = \QMV(\PLA,\PLV) - \QAV(\PLV, \ParSYS),  \\ % \label{eq:circ_VLV},
    %%%%%%%%%%%%%%%%%%%%%%%%%%%%%%%%%%%%%%%%%%%%%%%%%%
    &\CarPUL \dfrac{d \ParPUL}{d t} = \QPV(\PRV, \ParPUL)    - \QarPUL, \\ % \label{eq:circ_ParPUL}, \quad
    &\CarSYS \dfrac{d \ParSYS}{d t} = \QAV(\PLV, \ParSYS)    - \QarSYS, \\ % \label{eq:circ_ParSYS}, \quad
    &\dfrac{\LarPUL}{\RarPUL} \dfrac{d \QarPUL}{d t}  = - \QarPUL - \dfrac{\PvnPUL - \ParPUL}{\RarPUL}, \\ % \label{eq:circ_QarPUL}\\
    &\dfrac{\LarSYS}{\RarSYS} \dfrac{d \QarSYS}{d t}  = - \QarSYS - \dfrac{\PvnSYS - \ParSYS}{\RarSYS}, \\ % \label{eq:circ_QarSYS}\\
  \end{aligned}
\right.
\end{equation}
with $t \in [0,T]$ and where the flow rates of the valves read:
\begin{equation}
    Q_{\text{i}}(p_1, p_2) =
    \left\{ 
    \begin{split}
      \dfrac{p_1 - p_2}{\Rmin}, & \qquad \text{if } p_1 < p_2 \\
      \dfrac{p_1 - p_2}{\Rmax}, & \qquad \text{if } p_1 \geq p_2 \\
    \end{split}
    \right.
    \quad \text{for } \IinValve,
\end{equation}
% \begin{equation}
% \begin{aligned}
%     \QMV & = \dfrac{\PLA -\PLV    }{\RMV(\PLA, \PLV)}, \qquad&%\label{eq:circ_QMV} \\
%     \QAV & = \dfrac{\PLV -\ParSYS }{\RAV(\PLV, \ParSYS)}, \\%\label{eq:circ_QAV} \\
%     \QTV & = \dfrac{\PRA -\PRV    }{\RTV(\PRA, \PRV)},  \qquad&%\label{eq:circ_QTV}\\
%     \QPV & = \dfrac{\PRV -\ParPUL }{\RPV(\PRV, \ParPUL)}, %\label{eq:circ_QPV}
% \end{aligned}
% \label{eq:circ_valve_flow}
% \end{equation}
% where the functions $\RI$ define the behavior of the valves as diodes:
% \begin{equation}
%     R_{\text{i}}(p_1, p_2) =
%     \begin{cases}
%         \Rmin, & p_1 < p_2 \\
%         \Rmax, & p_1 \geq p_2 \\
%     \end{cases}
%     \quad \text{for } \IinValve,
% \end{equation}
where $p_1$ and $p_2$ denote the proximal and distal pressures of the valve, whereas $\Rmin$ and $\Rmax$ are its minimum and maximum resistance \citep{Regazzoni2022}.

While in the fully 0D model the four cardiac chambers consists of time-varying elastance elements \citep{Regazzoni2022}, in the 3D-0D whole-heart model the pressure-volume relationships of each chamber is provided by the 3D electromechanical model and must satisfy the volume-consistency conditions ($\PhyCoupl$) of \cref{eq:coupl},
where $\VIzerodim(\Circ(t))=\VI$, for $\IinCH$, represent the volumes of the four cardiac chambers in the 0D circulation model, while the 3D volumes are computed using the divergence (Gauss) theorem on the closed endocardial surfaces of the four cardiac chambers:
\begin{equation}
  \VIthreedim(\Disp(\x,t)) =
    \frac{1}{3} \int_{\GammaEndoI} J(\x,t) \left(\x + \Disp(\x,t) \right) \cdot \mecF^{-T}(\x,t) \, \mecNref(\x) \, d \mathbf{x}, \qquad \IinCH.
\end{equation}
We remark that these volumes can be exactly computed since the endocardial surfaces $\GammaEndoI$ are closed surfaces thanks to the presence of the valves $\OmegaValve$ and of the artificial caps $\OmegaCaps$ (see \cref{fig:domain}). 
The resulting model ($\PhyCirc$)-($\PhyCoupl$) of \cref{eq:circ,eq:coupl} consists of $n_\Circ + 4$ equations and unknowns, where the four additional unknowns are the chamber pressures $(\PRA,\PLA,\PRV,\PLV)$ that act as Lagrange multipliers enforcing the volume-consistency constraints.
These four pressures take into account the coupling with the ($\PhyMec$) model through the normal stress boundary conditions of \crefrange{eq:mec_bc_endo_RA}{eq:mec_bc_endo_LV} applied on the endocardium of the four chambers.
Instead, on the endothelium of the \ac{PT} and \ac{AO} we apply the pulmonary and systemic arterial pressures by setting $\PPT=\ParPUL$ and $\PAO=\ParSYS$ in \cref{eq:mec_bc_endo_PT,eq:mec_bc_endo_Ao}, respectively.

%%%%%%%%%%%
\subsection{Reference configuration and initial displacement}
\label{subsec:ref_conf_init_displ}

The most interesting applications of computational cardiac electromechanics occur when the human heart domain is directly reconstructed from medical images, with the aim of performing patient-specific simulations.
However, these reconstructed geometries correspond to a configuration $\OmegaTilde$ loaded by the internal blood pressure while, on the contrary, the stress-strain relationship at the basis of the mechanical model ($\PhyMec$) is formulated in an unloaded (stress-free) configuration $\Omega_0$ (see \cref{eq:piola}).
In order to recover this reference configuration $\Omega_0$ from the imaging configuration $\OmegaTilde$ we extend the procedure proposed by \citet{Regazzoni2022} for the \ac{LV} to the whole-heart case:
starting from $\OmegaTilde$, we recover the configuration $\Omega_0$ by virtually deflating the whole-heart domain previously subject to the internal pressures $\PItilde$, for $\IinCHAR$;
then, by applying on the endocardium and endothelium the pressures $\PIinit$, we inflate the domain again in order to compute the displacement $\Disp_0$ for the initial condition of the mechanical problem ($\PhyMec$) of \cref{eq:mec}.
Both these two steps are performed by assuming a quasi-static approximation of the mechanical problem \eqref{eq:mec} \citep{Regazzoni2022}.
This hypothesis is reasonable only in a few moments of the cardiac cycle, such as at the end of the \ac{VPF} phase just before the beginning of the \ac{AC} phase (see \cref{sec:heart,fig:cardiac_cycle}).
Indeed, at this time of the diastole, the ventricular filling slows down and the movement of the four chambers becomes negligible.
This moment is usually captured in standard cardiac medical images, because on the one hand it is easy to identify using the ECG signal, on the other hand the image quality is better when the heart moves slowly.
Furthermore, the small blood pressures that load the heart chambers during this phase make the associated numerical problem less challenging to solve.

More in detail, the procedure is based on the following quasi-static approximation, obtained by neglecting the time derivative term of \cref{eq:mec_pde} in the mechanical problem ($\PhyMec$):
\begin{subequations}
\label{eq:mec_static}
    \begin{empheq}[ left={ \left(\PhyMec^{\mathrm{static}}\right)\empheqlbrace\, } ]{align}
        & \nabla \cdot \tenspiola(\Disp, \TensConst) = \boldsymbol{0}
            & \text{in } \Omega_0 \times (0,T],
            \label{eq:mec_pde_static}
            \\
        & \tenspiola(\Disp, \TensConst) \, \mecNref + (\mecNref \otimes \mecNref) \left( \BCmecKepiN \Disp + \BCmecCepiN \dfrac{\partial \Disp}{\partial t} \right) = \mathbf{0}
            & \text{on } \GammaEpi_0 \times (0,T],
            \label{eq:mec_bc_epi_static}
            \\
        & \MecEndoBC{\TensConst}{\PRAconst}
            & \text{on } \GammaEndoRA_0 \times (0,T],
            \label{eq:mec_bc_endo_RA_static}
            \\
        & \MecEndoBC{\TensConst}{\PLAconst}
            & \text{on } \GammaEndoLA_0 \times (0,T],
            \label{eq:mec_bc_endo_LA_static}
            \\
        & \MecEndoBC{\TensConst}{\PRVconst}
            & \text{on } \GammaEndoRV_0 \times (0,T],
            \label{eq:mec_bc_endo_RV_static}
            \\
        & \MecEndoBC{\TensConst}{\PLVconst}
            & \text{on } \GammaEndoLV_0 \times (0,T],
            \label{eq:mec_bc_endo_LV_static}
            \\
        & \MecEndoBC{\TensConst}{\PPTconst}
            & \text{on } \GammaEndoPT_0 \times (0,T],
            \label{eq:mec_bc_endo_PT_static}
            \\
        & \MecEndoBC{\TensConst}{\PAOconst}
            & \text{on } \GammaEndoAo_0 \times (0,T],
            \label{eq:mec_bc_endo_Ao_static}
            \\
        & \Disp = \boldsymbol{0}
            & \text{on } \GammaRings_0 \times (0,T],
            \label{eq:mec_bc_rings_static}
    \end{empheq}
\end{subequations}
where $\TensConst>0$ represents the residual active tension and $\PIconst$, for $\IinCHAR$, are the constant pressures loading the endocardium and the endothelium.
% We denote with $\x_0$ and $\xtilde$ the coordinates corresponding to the unloaded configuration $\Omega_0$ and the imaging configuration $\OmegaTilde$, respectively.
Being $\x_0$ the coordinates associated to $\Omega_0$,
the solution $\overline{\Disp}=\Disp(\x_0, \PIconst, \TensConst)$ of \cref{eq:mec_static} can be used to move the coordinate $\x_0$ into a coordinate $\overline{\x}= \x_0 + \overline{\Disp}$ corresponding to a loaded configuration $\overline{\Omega}$.
Thus, in order to recover %the coordinate $\x_0$ corresponding to
the unloaded configuration $\Omega_0$ starting from the imaging configuration $\OmegaTilde$, we need to solve the following inverse problem:
find the domain $\Omega_0$ such that, if we displace $\x_0$ by the solution $\DispTilde=\Disp(\x_0, \PItilde, \TensTilde)$ of \cref{eq:mec_static}, we get the coordinate $\xtilde$ of the domain $\OmegaTilde$, i.e. $\xtilde = \x_0 + \DispTilde$.
To solve this problem we employ the algorithm proposed in \cite{regazzoni2020thesis,Regazzoni2022}, that is based on a fixed point method augmented with an adaptive step continuation method to ensure stability and boost convergence speed.

Finally, once the reference configuration $\Omega_0$ has been recovered, we can set proper values of $\PIconst = \PIinit$ and $\TensConst = \TensInit$ corresponding to the phase of the cardiac cycle at the initial time $t=0$ of the unsteady electromechanical model
and solve again \cref{eq:mec_static}.
In this way, we obtain the initial condition $\Disp_0=\Disp(\x_0,\PIinit,\TensInit)$ for the unsteady mechanical problem ($\PhyMec$) of \cref{eq:mec}.
Note that, in principle, the phase of the cardiac cycle corresponding to the initial time $t=0$ and the time when the imaging configuration $\OmegaTilde$ is acquired can be different, justifying possible different values of $\TensConst$ and $\PIconst$ during the reference configuration recovery and the initial displacement computation.

%%%%%%%%
%%%%%%%%
\section{Numerical approximation}
\label{sec:meth_num}

\begin{figure}[t]
    \centering
    \includegraphics[width=0.75\textwidth]{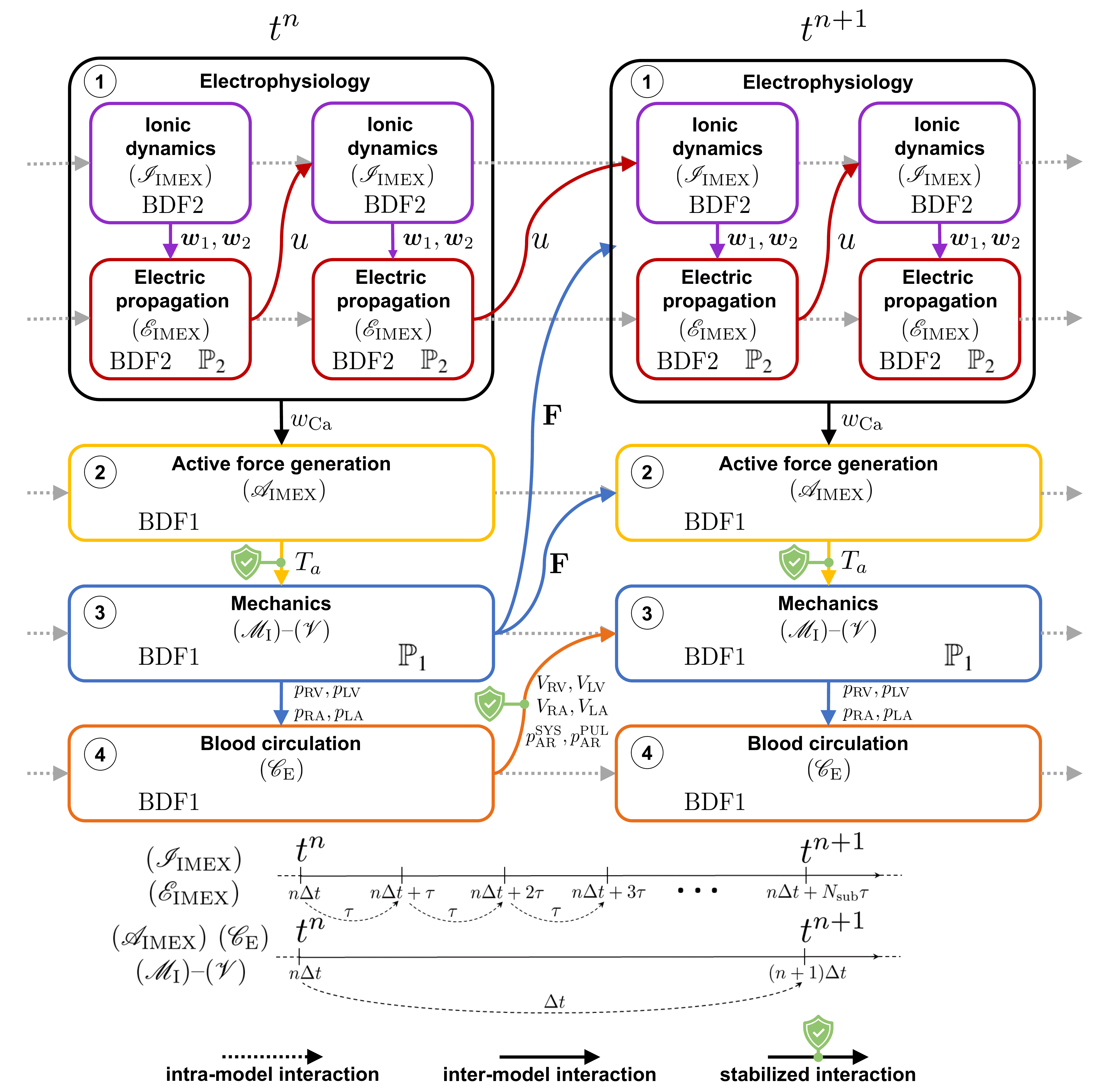}
    \caption{A sketch of the segregated-intergrid-staggered numerical scheme.
    Each block represents a core model and we show their order of resolution, which variables interconnect them and which interactions need stabilization.
    We also highlight the time and space discretization employed.
    Note that the electrophysiological block, being solved using a smaller timestep, features several repeated solutions of the $(\PhyEPIMEX)$ and $(\PhyIonIMEX)$ blocks for each time $t^n$ of the $(\PhyActIMEX)$--$(\PhyMecI)$--$(\PhyCoupl)$--$(\PhyCircE)$ blocks (in the figure, for illustrative purposes, only two sub-steps are displayed).
    }
    \label{fig:num_scheme}
\end{figure}

For the numerical approximation of the whole-heart electromechanical model (\crefrange{eq:ep}{eq:coupl}) we employ the segregated-intergrid-staggered numerical approach introduced for the ventricular cases in~\citep{Regazzoni2022,Piersanti2022biv}.
In this numerical scheme the core models are sequentially solved in a segregated manner, using different resolutions in space and time to properly take into account the heterogeneous space and time scales that characterize the different core models~\citep{nordsletten2011coupling,Quarteroni2017,Quarteroni2019}.
In \cref{fig:num_scheme} we show how the different core models are separately solved and which variables interconnect them, also highlighting which interactions need numerical stabilization.

\subsection{Numerical approximation of the core models}
\label{subsec:meth_num_core_models}

For the time discretization, we employ \ac{BDF} schemes~\citep{quarteroni2009numerical}.
The $(\PhyEP)$ and ($\PhyIon$) models are solved using a \ac{BDF2}, using an \ac{IMEX} scheme, denoted by $(\PhyEPIMEX)$ and $(\PhyIonIMEX)$, respectively, where the diffusion term is treated implicitly, the reaction term is treated explicitly and the ionic variables are advanced through the \ac{IMEX} scheme of~\cite{Regazzoni2022,Piersanti2022biv}. 
Moreover, the discretization of the ionic current term $\EPIion$ is performed following the \ac{ICI} approach \cite{krishnamoorthi2013numerical}.  
Both $(\PhyMec)$ and $(\PhyAct)$ models are advanced in time with a \ac{BDF1} scheme, with an \ac{IMEX} scheme for the activation $(\PhyActIMEX)$ \cite{regazzoni2020biophysically} and a fully implicit scheme for the mechanical problem $(\PhyMecI)$~\cite{Regazzoni2022,Piersanti2022biv}. 
Finally, we use an explicit \ac{BDF1} scheme for the circulation $(\PhyCircE)$~\cite{Piersanti2021phd}.

Concerning the space discretization, we use the \ac{FE} Method with continuous \acp{FE} and tetrahedral meshes~\cite{quarteroni2009numerical}.
We consider a unique mesh $\mathcal{T}_{\h}$ ($\h$ represents the mesh size) for the entire computational domain~$\Omega_0$ (see \cref{fig:mesh_stimuli_at}, (a)).
We employ a scalable and efficient intergrid transfer operator on the unique mesh $\mathcal{T}_{\h}$ that enables the use of arbitrary \acp{FE} among the different core models.
In particular, we consider \ac{FE} of order 2 ($\mathbb{P}_2$) for $(\PhyEPIMEX)$ to properly capture the dynamics of traveling waves, and \ac{FE} of order 1 ($\mathbb{P}_1$) for both $(\PhyActIMEX)$ and $(\PhyMecI)$~\cite{Augustin2016,ColliFranzone3,Regazzoni2022,Piersanti2022biv}.

Regarding the $(\PhyEPIMEX)$-$(\PhyIonIMEX)$-$(\PhyActIMEX)$ models, which are defined only on the subdomain $\OmegaMyo$, we assemble and solve the \ac{FE} system on the cells and \ac{DOFs} of the mesh $\mathcal{T}_{\h}$ corresponding to $\OmegaMyo$, neglecting the cells and \ac{DOFs} belonging only to non-conductive regions $\{\Omega_0 \setminus \OmegaMyo\}$.
In $(\PhyEPIMEX)$, this approach models the atrioventricular valves as electrical insulators between atria and ventricles, representing the discrete counterpart of the homogeneous Neumann condition \eqref{eq:ep_monodomain_bc} on the internal interfaces between the conductive and non-conductive regions.
At the same time, it allows to use a unique mesh for all the core models, making the intergrid transfer operator more efficient and easier to define.

\subsection{Numerical coupling of the core models}

We adopt a segregated approach to couple the different core models, solving them in a sequential manner.
Moreover, we make use of two different time steps, a larger one (denoted by $\Delta t$) for $(\PhyActIMEX)$--$(\PhyMecI)$--$(\PhyCoupl)$--$(\PhyCircE)$ and a finer one (that is $\tau = \Delta t/\Nsub$) for $(\PhyEPIMEX)$--$(\PhyIonIMEX)$, with $\Nsub \in \mathbb{N}$, see \cref{fig:num_scheme}(b).
As shown in \cref{fig:num_scheme}, we update the variables in the following order: first, we update $(\PhyIonIMEX)$ and $(\PhyEPIMEX)$, by performing $\Nsub$ sub-steps; then, we update $(\PhyActIMEX)$; successively, we update $(\PhyMecI)$ together with the constraint $(\PhyCoupl)$ (more details are provided below); finally, we update $(\PhyCircE)$.

This ordering of the core models is defined to reflect the main direction of the interactions among the core models. 
The interactions that occur in the opposite direction, the so-called feedbacks, are instead evaluated using the solution available from the previous time-step (see, e.g., feedback from mechanics to electrophysiology).
To evaluate the feedbacks between $(\PhyIonIMEX)$--$(\PhyEPIMEX)$ and $(\PhyActIMEX)$--$(\PhyMecI)$, 
we employ the intergrid transfer operator described in \cref{subsec:meth_num_core_models}.
We refer to~\cite{salvador2020intergrid,Regazzoni2022,Piersanti2022biv} for further details.

\subsection{Stabilizing the coupling of the core models}
\label{sec:numerics:stabilization}

The use of segregated schemes can lead to numerical instabilities, especially when feedbacks play a non-negligible role. 
In the case of cardiac electromechanics, numerical instabilities can arise, on the one hand, due to feedbacks between mechanics and activation \cite{whiteley2007soft,niederer2008improved,pathmanathan2009numerical,pathmanathan2010cardiac,regazzoni2020thesis,regazzoni2021oscillation} and, on the other hand, due to feedbacks between circulation and active-passive mechanics \cite{Hirschvogel,regazzoni2022stabilization}.
These instabilities, which yield non-physical oscillations, do not affect monolithic methods, which however require higher computational costs than segregated schemes. 
Furthermore, they force the use of a single time step size for all the core models.
With the aim of preserving the advantages of segregated schemes, we use stabilization terms aimed at curing the numerical oscillations.
Specifically, we employ the stabilization schemes that we proposed in \cite{regazzoni2021oscillation} and \cite{regazzoni2022stabilization}.
Since both schemes act on the $(\PhyMecI)$--$(\PhyCoupl)$ substep, in what follows we provide more detail on this block.

We update the mechanical displacement variable under the constraint of assigned chamber volumes.
The chamber pressures ($\PRA$, $\PLA$, $\PRV$ and $\PLV$) are determined simultaneously with the displacement and play in this context the role of Lagrange multipliers enforcing the volume conservation constraints~$(\PhyCoupl)$. 
Introducing the discrete times $t^n=n\Delta t$ (with $n\ge 0$) and denoting by $\generica_{\h}^{n} \simeq \generica_{\h}(t^n)$ the fully discretized \ac{FE} approximation of the generic (scalar $a$, vectorial $\generica$ or tensorial $\mathbf{A}$) variable $\generica(t)$, we consider the following fully discretized version of the coupled $(\PhyMec)$--$(\PhyCoupl)$ models of \cref{eq:mec,eq:coupl}.

\noindent
For each time step $t^{n+1}$, given ${\Tens}_{\h}^{n+1}$ and $\Circ^{n}$, find $\Disp_{\h}^{n+1}$, $p_\LA^{n+1}$, $p_\RA^{n+1}$, $p_\LV^{n+1}$ and $p_\RV^{n+1}$ by solving:
\begin{equation}
    \left\{
    \begin{split}
        &\int_{\Omega_0} \rho_s \dfrac{\Disp_{\h}^{n+1}-2\Disp_{\h}^{n}+\Disp_{\h}^{n-1}}{\Delta t^2} \cdot \boldsymbol{\varphi}_{\h} \, d \Omega_0 +
        \int_{\Omega_0} \tenspiola(\Disp_{\h}^{n+1}, {\Tens}_{\h}^{n+1}):\nabla \boldsymbol{\varphi}_{\h} \, d \Omega_0 \: + 
        \\
        &\quad+\int_{\GammaEpi_0} \BCmecCepiN\dfrac{\Disp_{\h}^{n+1}-\Disp_{\h}^{n}}{\Delta t}(\mecNref_{\h} \otimes \mecNref_{\h}) \cdot \boldsymbol{\varphi}_{\h} \, d \Gamma_0 + 
        \int_{\GammaEpi_0}  \BCmecKepiN(\mecNref_{\h} \otimes \mecNref_{\h}) \, \Disp_{\h}^{n+1} \cdot \boldsymbol{\varphi}_{\h} \, d \Gamma_0 \: + 
        \\       
        &\quad+\sum_{\KinCH} p_k^{n+1} \int_{\GammaEndok_0} J_{\h}^{n+1} ({\mecF_{\h}^{n+1}})^{-T} \mecNref_{\h} \cdot \boldsymbol{\varphi}_{\h} \, d \Gamma_0  +\\
        &\quad+\sum_{k \in \{ \AO, \PT \}} p_k^{n} \int_{\GammaEndok_0} J_{\h}^{n+1} ({\mecF_{\h}^{n+1}})^{-T} \mecNref_{\h} \cdot \boldsymbol{\varphi}_{\h} \, d \Gamma_0  = 0 \qquad  
        \\ 
        &\pushright{\forall \boldsymbol{\varphi}_{\h} \in [\mathcal{X}_{\text{h}}^s]^3},
        \\
        &\VLAthreedimH(\Disp_{\h}^{n+1}) = \VLAzerodim(\Circ^{n}) \\
        &\VLVthreedimH(\Disp_{\h}^{n+1}) = \VLVzerodim(\Circ^{n}) \\
        &\VRAthreedimH(\Disp_{\h}^{n+1}) = \VRAzerodim(\Circ^{n}) \\
        &\VRVthreedimH(\Disp_{\h}^{n+1}) = \VRVzerodim(\Circ^{n}) \\
    \end{split}
    \right.
    \label{eqn: TMSI-heart}
\end{equation}
with $\mecF^{n+1}_{\text{h}}= \identity + \nabla \Disp_{\h}^{n+1}$, $J^{n+1}_{\h} = \det(\mecF^{n+1}_{\h})$ and $\boldsymbol{\varphi}_{\h}$ being a generic test function for the finite dimensional space $[\mathcal{X}_{\text{h}}^s]^3$ with $\mathcal{X}_{\text{h}}^s= \{ v \in C^0(\overline{\Omega}_0) : v|_\text{K} \in \mathbb{P}_\text{s}(K), \, s \geq 1, \; \forall K \in \mathcal{T}_{\h}, v = 0 \text{ on } \GammaRings_0 \}$, where $\mathbb{P}_\text{s}(K)$ stands for the set of polynomials with degree smaller than or equal to $s$ over a mesh element $K$.
We also remark that, unlike chamber pressures, arterial pressures ($p_\PT^n$ and $p_\AO^n$) are evaluated at the time step $t^n$ since they are equal to the pulmonary and systemic arterial pressures ($\ParPULn$ and $\ParSYSn$, respectively) of the circulation state vector $\Circ^{n}$.

As mentioned above, the formulation of \cref{eqn: TMSI-heart} typically exhibits numerical oscillations when coupled with an active force model on the one hand, and a circulation model on the other hand.

One source of instability is represented by the fibers-stretch-rate feedback, i.e. the influence that the rate at which fibers shorten has on the amount of force generated at each point in the domain.
As shown in \cite{regazzoni2021oscillation}, these numerical oscillations originate from an inconsistent description of strain, which is represented in Eulerian coordinates at the microscale, i.e. in activation models, in Lagrangian coordinates instead at the macroscale, i.e. in the tissue mechanics model. 
This can be corrected by introducing an additional term in the formulation, which constitutes a numerically consistent stabilization term.
This numerical scheme is obtained by replacing the Piola tensor expression in \cref{eqn: TMSI-heart} with the following expression:
\begin{equation}\label{eqn:stab:active-stress}
    \tenspiola\left(\Disp_{\h}^{n+1}, {\Tens}_{\h}^{n+1} + {\Stiff}_{\h}^{n+1} \left(
        \sqrt{\mecF^{n+1}_{\h} \f_0 \cdot \mecF^{n+1}_{\h} \f_0}
        -
        \sqrt{\mecF^{n}_{\h} \f_0 \cdot \mecF^{n}_{\h} \f_0}
    \right) \right)
\end{equation}
where $\Stiff$ represents the total stiffness of the attached crossbridges, and is obtained from the activation model (see~\citep{regazzoni2021oscillation} for more details).

A second source of instability is related to the interaction between active-passive mechanics and circulation.
As discussed in \cite{regazzoni2022stabilization}, the staggered scheme of \cref{eqn: TMSI-heart} is not unconditionally stable, but can exhibit non-physical oscillations for given values of the parameters and $\Delta t$.
This occurs, for example, for sufficiently large values of inertia, viscous dissipation and stiffness, or again as a consequence of fibers-stretch-rate feedback, which leads to an increase in apparent stiffness.
In order to cure these oscillations without resorting to a monolithic scheme, we take inspiration from \cite{regazzoni2022stabilization} and we correct the volume constraint in \cref{eqn: TMSI-heart}, namely $\Vthreedimk(\Disp_{h}^{n+1}) = \VKzerodim(\Circ^{n})$ for $\KinCH$.
In particular, the volumes derived from the circulation model at time $t_n$ are replaced by their extrapolation at time $t_{n+1}$, which takes into account the effect that the variation of the pressures in the four chambers will have on the fluxes through the valves.
More precisely, the volume constraints of the stabilized scheme read:
\begin{equation}\label{eqn:volume_constraints_stabilized}
    \left\{
    \begin{split}
        \VLAthreedimH(\Disp_{\h}^{n+1}) &= \VLAzerodim(\Circ^{n}) 
        + \Delta t \left[ 
                Q_{\mathrm{VEN}}^{\mathrm{PUL,n}} -
                \QMV(p_\LA^{n+1},p_\LV^{n+1}) 
          \right]\\
        \VLVthreedimH(\Disp_{\h}^{n+1}) &= \VLVzerodim(\Circ^{n})
        + \Delta t \left[ 
                \QMV(p_\LA^{n+1},p_\LV^{n+1}) - 
                \QAV(p_\LV^{n+1}, p_{\mathrm{AR}}^{\mathrm{SYS, n}}) 
          \right]\\
        \VRAthreedimH(\Disp_{\h}^{n+1}) &= \VRAzerodim(\Circ^{n}) 
        + \Delta t \left[ 
                Q_{\mathrm{VEN}}^{\mathrm{SYS,n}} -
                \QTV(p_\RA^{n+1},p_\RV^{n+1}) 
          \right]\\
        \VRVthreedimH(\Disp_{\h}^{n+1}) &= \VRVzerodim(\Circ^{n}) 
        + \Delta t \left[ 
                \QTV(p_\RA^{n+1},p_\RV^{n+1}) - 
                \QPV(p_\RV^{n+1}, p_{\mathrm{AR}}^{\mathrm{PUL, n}}) 
          \right]\\
    \end{split} 
    \right.
\end{equation}
We remark that in \cref{eqn:volume_constraints_stabilized}, while the pressures in the four chambers are evaluated at time $t_{n+1}$, the state variables of the circulation model are evaluated at time $t_{n}$. 
In other words, the \cref{eqn:volume_constraints_stabilized} does not invalidate the staggered nature of the scheme.
Nevertheless, the additional terms allow for the removal of numerical oscillations.
Indeed, it is shown in \cite{regazzoni2022stabilization} that this scheme is absolutely stable for any choice of parameters and $\Delta t$.
These stabilization terms are also straightforward to implement and have no impact on the computational cost.
Indeed, the fully discretized version of the stabilized version of system~\eqref{eqn: TMSI-heart} can be compactly written as: 
\begin{equation}
    \begin{cases}\label{eqn: saddle-point-heart}
        \mathbf{r}_{\Disp}(\Disp_{\h}^{n+1}, p_\LA^{n+1}, p_\LV^{n+1}, p_\RA^{n+1}, p_\RV^{n+1}) &= \mathbf{0}, \\
        r_{p_\LA}         (\Disp_{\h}^{n+1}, p_\LA^{n+1}, p_\LV^{n+1}) &= 0, \\
        r_{p_\LV}         (\Disp_{\h}^{n+1}, p_\LA^{n+1}, p_\LV^{n+1}) &= 0, \\
        r_{p_\RA}         (\Disp_{\h}^{n+1}, p_\RA^{n+1}, p_\RV^{n+1}) &= 0, \\
        r_{p_\RV}         (\Disp_{\h}^{n+1}, p_\RA^{n+1}, p_\RV^{n+1}) &= 0, \\
    \end{cases}
\end{equation}
where we moved all the terms to the left hand side and $r_{p_\RA}$, $r_{p_\LA}$,$r_{p_\RV}$, $r_{p_\LV}$ and $\mathbf{r}_{\Disp}$ are suitable functions. 
\cref{eqn: saddle-point-heart} is a nonlinear saddle-point problem, that we solve by means of the Newton algorithm using the Schur complement reduction~\cite{benzi2005numerical,Regazzoni2022,Piersanti2022biv}.
As shown in \cite{regazzoni2022stabilization}, this can be done at the cost of 5 solutions of the linear system (that is, the number or chambers plus one) associated with the Jacobian matrix of the standalone mechanical subproblem for each Newton iteration.

%%%%%%%%
%%%%%%%%
\section{Numerical Simulations and Discussion}
\label{sec:res}
In this section we display and discuss the results obtained using our whole-heart electromechanical model.
More specifically, in \cref{subsec:res_settings} we summarize the common settings for all the numerical simulations.
In \cref{subsec:res_baseline} we show the results of a baseline simulation.
Eventually, in \cref{subsec:res_effects_math_mod,subsec:res_effects_num_mod}, we show the impact of some features of our computational model, as the atrial contraction, the fibers-stretch-rate feedback and the numerical stabilization terms.

%%%%%%%%%
\subsection{Simulation setup}
\label{subsec:res_settings}

\begin{figure}[t]
    \centering
    \includegraphics[width=1\textwidth]{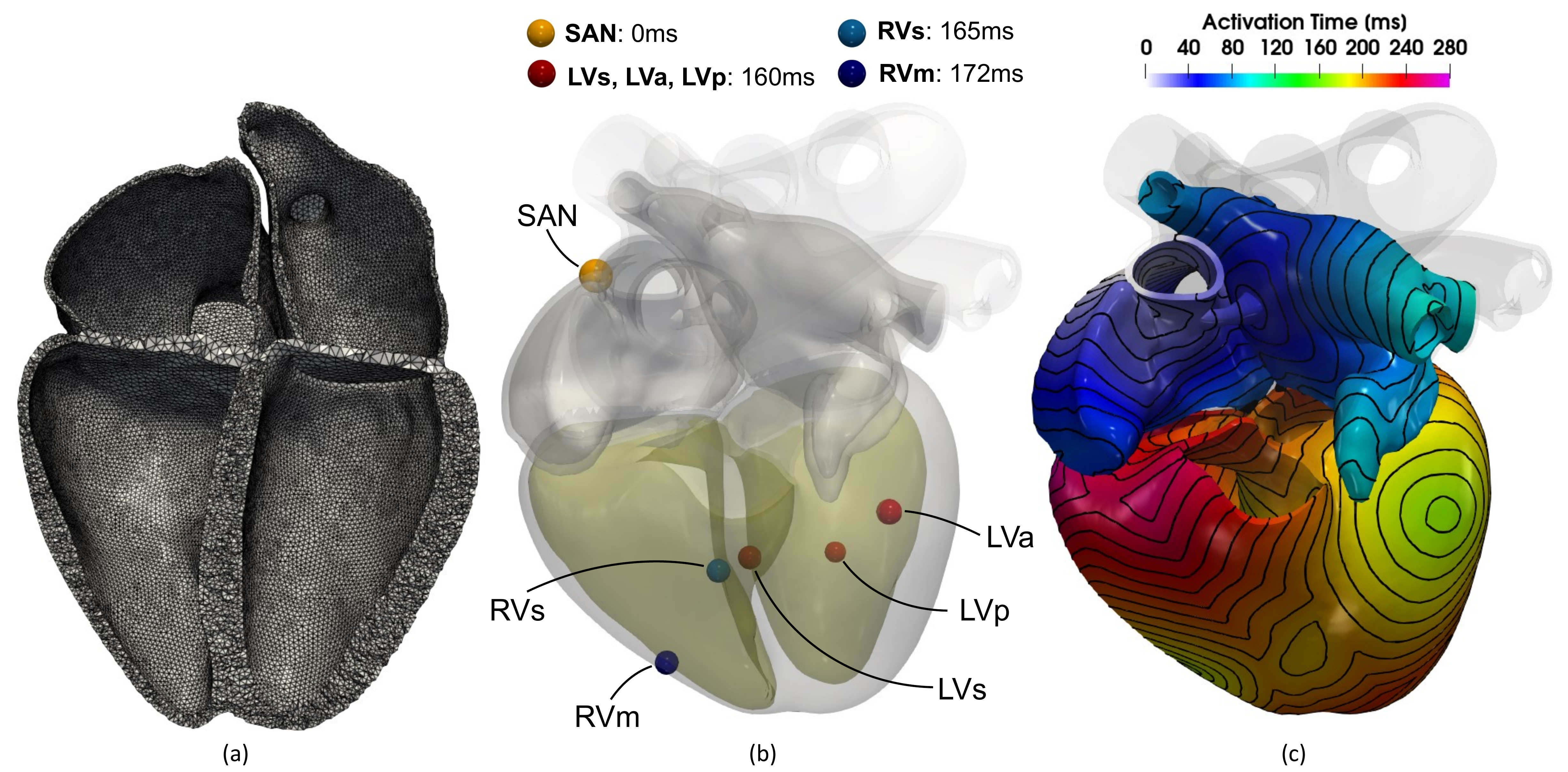}
    \caption{
        (a) A cut view of the computational mesh.
        (b) The stimulation protocol highlighting the location and time of the applied spherical impulses.
        (c) The baseline simulation results in term of activation time.}
    \label{fig:mesh_stimuli_at}
\end{figure}

We generate the computational mesh starting from the Zygote Solid 3D Heart Model \citep{Zygote2014},
an anatomically accurate CAD model of the entire human heart reconstructed from high-resolution CT scans and representing a healthy male subject from the $50^\text{th}$ percentile of the United States population.
The original model -- made of disjoint parts of the various cardiac compartments -- has been processed to fit the domain features described in \cref{subsec:domain}.
With this purpose, we rely on the algorithms recently proposed by \citet{Fedele2021} to facilitate the surface processing and mesh generation of cardiac geometries, implemented in the open source software \texttt{vmtk}\footnote[2]{\url{https://github.com/marco-fedele/vmtk}}~\citep{Antiga2008}.
In particular, we extensively use the \texttt{surface-connection}, \texttt{boolean-connection}, \texttt{surface-tagger}, and \texttt{mesh-connector} algorithms \citep{Fedele2021}.
The final tetrahedral computational mesh is shown in \cref{fig:mesh_stimuli_at}, (a).
This mesh is characterized by a mesh size of about $\SI{1.5}{\milli\meter}$ in the myocardium -- i.e. in the conductive regions where also the electrophysiology and the active force generation model are solved -- and of about $\SI{3}{\milli\meter}$ in the non-conductive regions -- where only the mechanical model is solved with a less demanding isotropic Neo-Hookean constitutive law..
% We remark that the electrophysiology is solved using $\mathbb{P}_2$ elements, thus, in terms of convergence properties, the used mesh-size is equivalent or better than a mesh-size of $0.5\,mm$ CITARE.
Starting from this mesh -- that represents the domain in the imaging configuration $\OmegaTilde$ -- we recover the reference configuration $\Omega_0$ by solving the problem illustrated in \cref{subsec:ref_conf_init_displ}.
The resulting deformed mesh is then remeshed to improve the quality of the elements that can be adversely affected by the deformation procedure, especially in the anatomically complex and thin regions of the atria.

The aforementioned mesh is used only for the baseline simulation (\cref{subsec:res_baseline}), while for the tests described in \cref{subsec:res_effects_math_mod,subsec:res_effects_num_mod}, in order to reduce the computational burden of the numerical simulations, we take advantage of a coarser mesh characterized by a mesh-size of about $\SI{3}{\milli\meter}$ also in the conductive regions.
Indeed, those tests aim at describing the qualitative effects of some changes in the models and the quantities analyzed are not significantly affected by the coarsening of the mesh.
% This choice is motivated by the fact that these tests aim at describing the qualitative effects of some changes in the models.
% Moreover, the quantities analyzed are not significantly affected by the coarsening of the mesh % PIU' DETTAGLI QUANTITATIVI change of pressures of maximum ... and of volumes of maximum ....
The fine and coarse meshes are made up of $1.34\text{M}$ and $270\text{K}$ elements and $229\text{K}$ and $51\text{K}$ vertexes, respectively.
The corresponding number of \ac{DOFs} relative to the electrical $(\PhyEPIMEX)$ and mechanical $(\PhyMecI)$ \ac{FE} problems are $1.71\text{M}$ and $687\text{K}$, respectively, for the fine mesh, and $337\text{K}$ and $154\text{K}$, respectively, for the coarse mesh.

Concerning the time steps, we use $\tau=50\si{\micro \second}$ for the electrophysiology and $\Delta t = 1000\si{\micro \second}$ for the mechanical, activation and circulation problems \cite{Piersanti2022biv,Piersanti2021phd}.
All the other parameters of the baseline simulation (\cref{subsec:res_baseline}) are listed in \ref{sec:appendix_params}.
% variations for other tests are specified from time to time.
We simulate $9$ and $6$ heartbeats for the baseline simulation and the other tests, respectively, showing the results of the last two heartbeats, when the circulation variables reach their limit cycle.

In all the presented simulations the cardiac electrical conduction system (see \cref{fig:heart_anatomy}, (c)) is modeled using the same series of spherical impulses (see \cref{subsubsec:ep}).
We first stimulate the atrial muscles at the \ac{SAN} allowing the propagation of the signal in the \ac{RA} and, through the \ac{BB} and the atrial septum, toward the \ac{LA};
waiting for the natural delay governed by the \ac{AVN}, we then stimulate a series of points on the endocardium of the two ventricles that, together with the fast endocardial layer, surrogate the effect of the \ac{PFs}.
More in detail, the ventricular stimuli are first applied to the \ac{LV} and soon after to the \ac{RV} in order to model the physiological lag of the \ac{RBB} with respect to the \ac{LBB}.
The whole stimulation protocol -- detailed in \cref{fig:mesh_stimuli_at}, (b) -- is periodically repeated every heartbeat, representing a simplified but effective model of the pacemaking activity of the \ac{SAN} and of the entire electrical conduction system.
% \todoinline{@RP Aggiungere referenze sui tempi/ritardi scelti?}

We initialize the ionic models by running a 1000-cycle long single-cell simulation for each model.
Similarly, we run single-cell simulations for the force generation models with a constant calcium input ($\wCa=0.1 \si{\micro \mole}$) and
a reference sarcomere length $\SL = 2.2 \si{\micro \meter}$ \cite{Piersanti2022biv}. 

The numerical framework presented in Section~\ref{sec:meth_num} has been implemented in \lifex, an in-house high-performance \texttt{C++} FE library for cardiac applications, based on the \texttt{deal.II}\footnote[3]{\url{https://www.dealii.org}} FE core~\cite{dealII91}.
A public binary release of \lifex (including the fiber generation package) is freely available online, under an open license\footnote[4]{\url{https://doi.org/10.5281/zenodo.5810269}}~\cite{africa2022lifex}.
All the numerical simulations were performed using either the \texttt{iHeart} cluster (Lenovo SR950 192-Core Intel Xeon Platinum 8160, 2100 MHz and 1.7TB RAM) at MOX, Dipartimento di Matematica, Politecnico di Milano or the \texttt{GALILEO100} supercomputer at Cineca (24~nodes endowed with 48-Core Intel CascadeLake 8260, 2.4GHz, 384 GB RAM).
A simulation of one heartbeat lasts for the fine mesh about 4 hours with 1152 cores on the \texttt{GALILEO100} supercomputer, for the coarse mesh about 4.75 hours with 48 cores of the \texttt{iHeart} cluster.

%%%%%%%%%
\subsection{The baseline simulation}
\label{subsec:res_baseline}

\begin{figure}[t]
    \centering
    \includegraphics[width=0.8\textwidth]{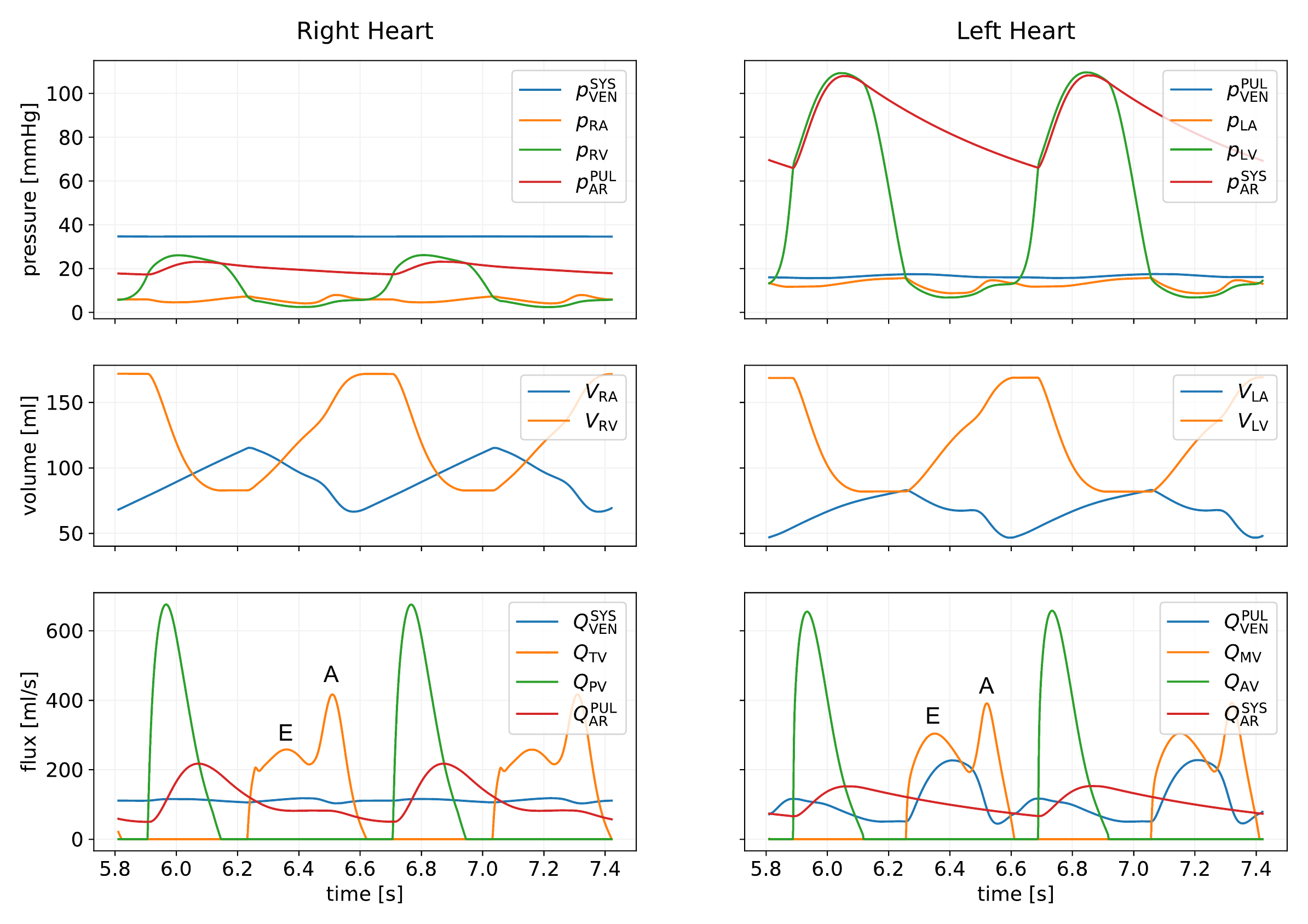}
    \caption{Pressure, volumes, and fluxes evolution over time during the last two heartbeats of the state variables of the coupled mechanics-circulation model for the baseline simulation.
    For the \ac{MV} and \ac{TV} fluxes we also highlight the E-wave and the A-wave.}
    \label{fig:res_baseline_state}
\end{figure}

In \cref{fig:res_baseline_state} we show the temporal variation of some state variables during the last two heartbeats of the baseline simulation.
The obtained curves of pressures, volumes and fluxes qualitatively correspond to those expected for a physiological heart function~\citep{Mitchell2014,Askari2019}.

Concerning the systolic function, we obtain an excellent agreement with reference values for a healthy adult available in the literature.
Indeed, the maximum fluxes obtained through the semilunar valves ($\QPV,\QAV$) during the \ac{VE} phase (about \SI{600}{\milli \liter \per \second}) are in the physiological range usually measured by PC-MRI data (500--600 \si{\milli \liter \per \second}) \citep{Gallo2012,Alastruey2016,Lantz2014}.
This feature is hardly achieved by computational models which tend to largely overestimate these fluxes, even when they reproduce the physiological ventricular output in terms of \ac{SV} (see, e.g., \citep[Fig. 10]{Gerach2021}).
%For instance, the whole heart model of \citet{Gerach2021} obtained the maximum fluxes of about $2000 \, ml/s$ and $1500 \, ml/s$ in the aorta ($\QAV$) and pulmonary artery ($\QPV$), respectively.
As we will show in \cref{subsec:res_effects_math_mod}, a key component of our model to achieve this result is the fibers-stretch-rate feedback accounted for by the RDQ20 model, that homogenizes the fibers shortening velocity and contributes to regulate the blood fluxes.

Concerning the diastolic function, instead, the \textit{atrial kick} is clearly visible during the \ac{AC} phase, for both ventricular and atrial volumes ($\VRA[],\VRV[],\VLA[],\VLV[]$), atrial pressures ($\PRA[],\PLA[]$), and fluxes through the atrioventricular valves ($\QTV,\QMV$).
However, the fluxes during a healthy diastolic function should be characterized by an E-wave -- corresponding to the \ac{VPF} phase -- which is taller than the A-wave -- corresponding to the \ac{AC} phase \citep{Galderisi2005,Nagueh2020}.
In other words, the ventricular filling should be mainly determined by the ventricular relaxation than by the atrial contraction. 
The different behavior that we obtain (see \cref{fig:res_baseline_state}, last row) can be motivated by a
%too strong atrial contraction or a 
too slow ventricular relaxation during the \ac{VPF} phase \citep{Galderisi2005,Nagueh2020}.
We expect that a better agreement with literature data can be obtained by resorting to ionic models with a more realistic decrease transient of calcium concentration \citep{TTP06,ToR-ORd}.

% Add also a comment about pulmonary veins flow and its waves (\citet{Pagel2003})?

\begin{figure}[t]
    \centering
    \includegraphics[width=1.0\textwidth]{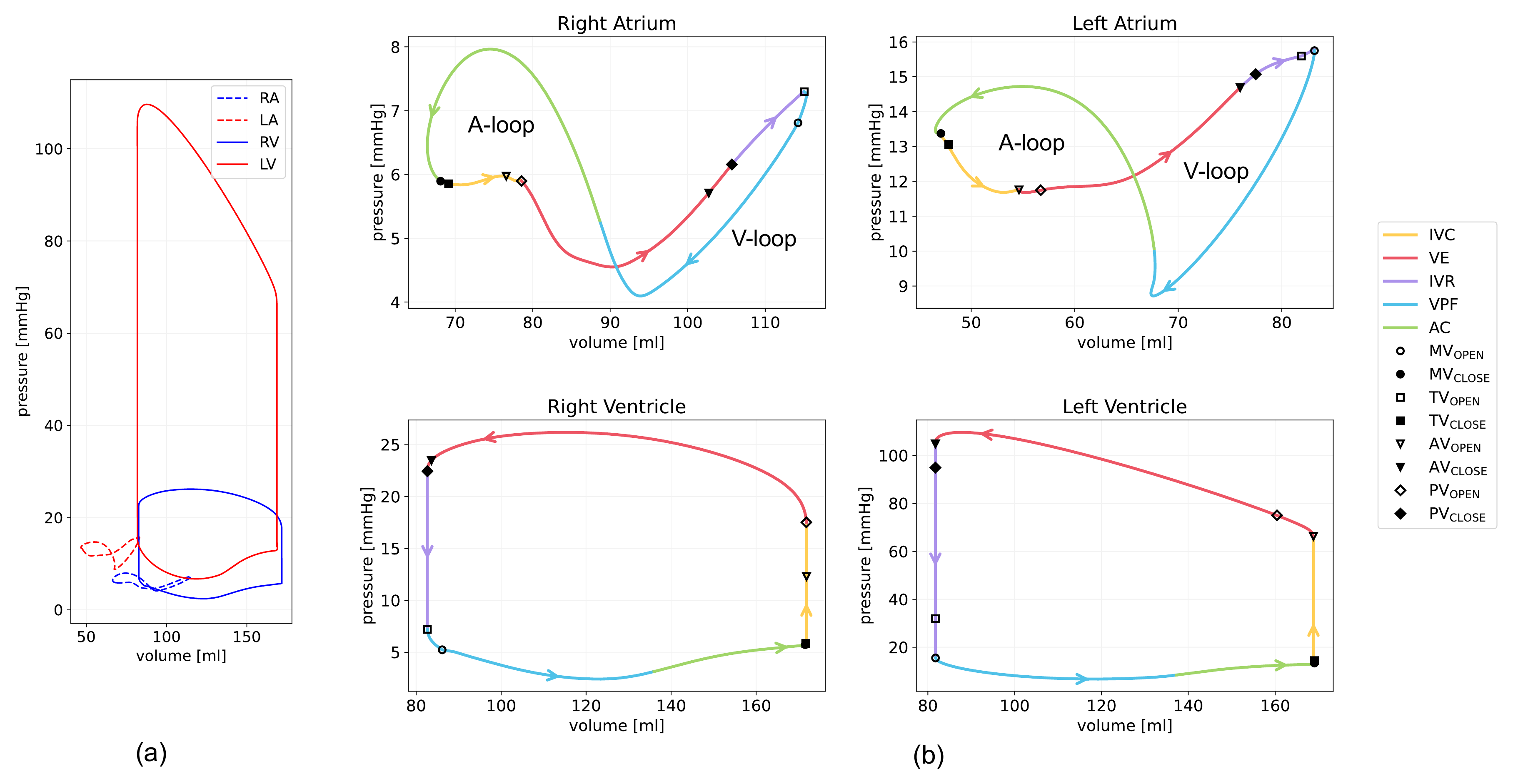}
    \caption{The pressure-volume loops of the last heartbeat for the baseline simulation:
    (a) all the curves in a single figure to highlight the difference in volumes and pressures among the cardiac chambers;
    (b) the curves colored with the phases of the cardiac cycle and highlighting the opening and closing time of each cardiac valve.
    Abbreviations are defined in \cref{tab:acro}.}
    \label{fig:res_baseline_pvloop}
\end{figure}

\begin{figure}[t]
    \centering
    \includegraphics[width=1.0\textwidth]{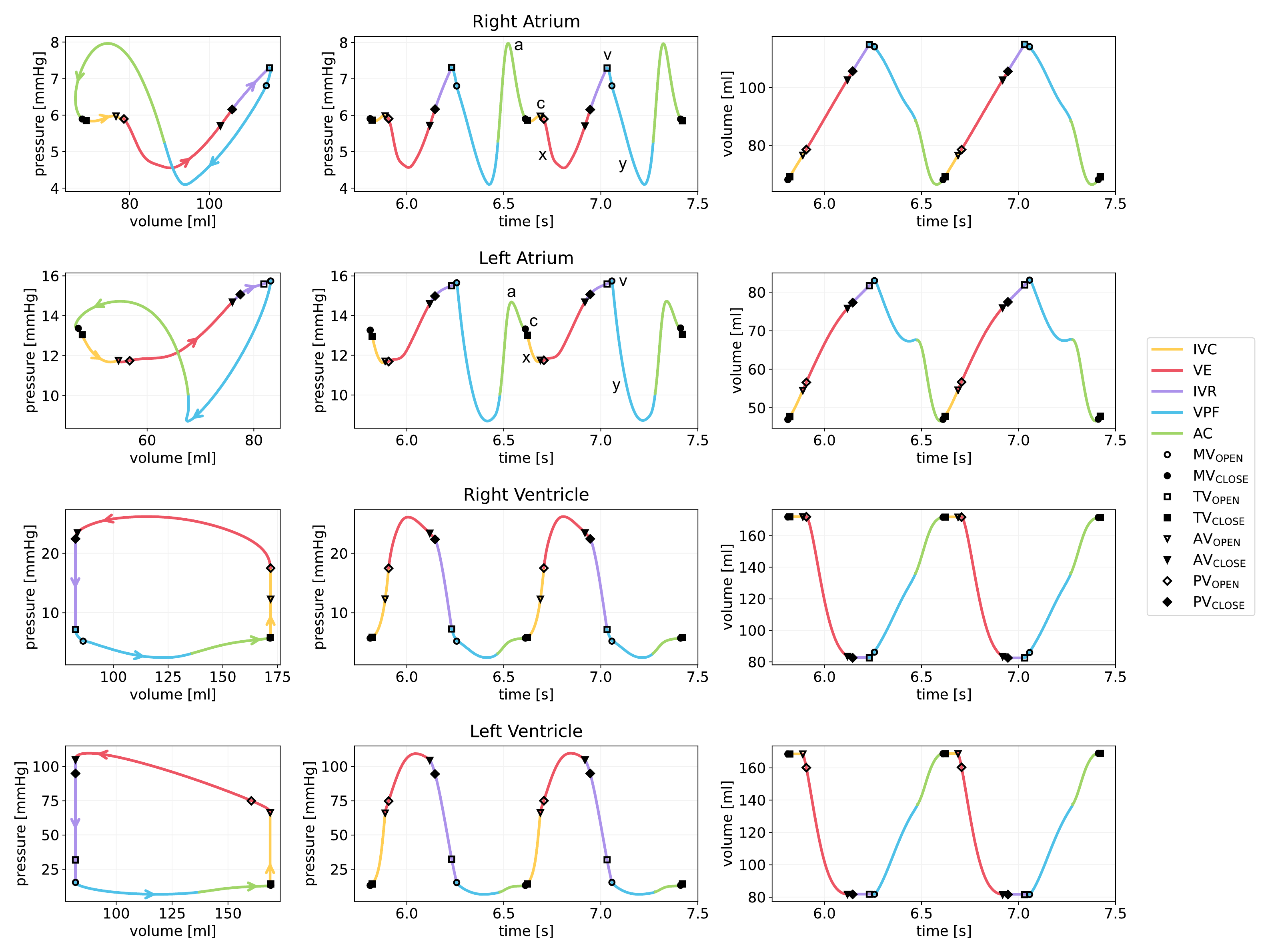}
    \caption{The phases of the cardiac cycle for the four heart chambers both in terms of pressure-volume loop (last heartbeat) and in terms of pressure and volume evolution over time (last two heartbeats).
    The opening and closing time of each cardiac valve is also reported.
    For the atrial pressure evolution we also highlight the a-c-v waves and the x-y descents.
    Abbreviations are defined in \cref{tab:acro}.}
    \label{fig:res_baseline_pvloopext}
\end{figure}

In \cref{fig:res_baseline_pvloop} we show the pressure-volume loops of the four cardiac chambers, while in \cref{fig:res_baseline_pvloopext} we display the same curves together with the evolution over time of the pressures and volumes.
In both figures we color each curve with the five phases of the cardiac cycle described in \cref{sec:heart,fig:cardiac_cycle} and we also represent the opening and closing moments of the four cardiac valves.
As depicted in \cref{fig:res_baseline_pvloop}, (a), the pressure and volume ranges vary significantly among the four cardiac chambers, as reported, e.g., by \citet[Fig. 10]{Verzicco2022}).
The shape of the pressure-volume loops finds a very good agreement with the medical literature \citep{Pagel2003,Blume2011a,Rocsca2011,Klabunde2011,Verzicco2022}.
While this is not the first time that an electromechanical model is able to describe the ventricular physiology \citep{Regazzoni2022,Piersanti2022biv,Gerach2021}, to the best of our knowledge the eight-shaped pressure-volume loops of the atria have never been shown so accurately by a computational model.
Indeed, we obtain, as expected by the literature \citep{Pagel2003,Blume2011a}, A- and V-loops that are similar in size.
On the contrary, A-loops significantly (and abnormally) larger than V-loops are obtained by the few other whole-heart electromechanical models accounting for atrial contraction~\citep{Land2018,Gerach2021}.

The atrial function of reservoir, conduit, and booster pump is also well captured:
the total emptying volume ($V_{\mathrm{max}} - V_{\mathrm{min}}$, reservoir) is divided between the passive emptying volume ($V_{\mathrm{max}} - V_{\mathrm{preAC}}$, conduit) and the active emptying volume ($V_{\mathrm{preAC}} - V_{\mathrm{min}}$, booster pump), with these last two volumes comparable in size \citep{Blume2011a,Peluso2013,Li2017a}.
The contribution of the atrial booster pump function to the ventricular filling falls within the physiological upper limit.
More quantitatively, the \ac{LA} contraction contributes to the $33\%$ of the \ac{LV} filling, while normal healthy values are reported in the range $15\%$--$30\%$ \citep{Pagel2003,Blume2011a,Thomas2020}.
%approximately the $20\%$ according to \citet{Pagel2003}, in the range $15$--$30\%$ according to \citet{Blume2011a}, in the range $15$--$30\%$ according to \citet{Thomas2020}.

The evolution over time of the atrial volume is very well captured.
In particular, the \ac{LA} curve (see \cref{fig:res_baseline_pvloopext}, second row, last plot) matches similar curves reconstructed from medical images (see, e.g., \citet[Fig. 7]{Thomas2020} and \citet[Fig. 3]{Badano2016}):
the volume smoothly increases when the \ac{MV} is closed (\ac{IVC}, \ac{VE} and \ac{IVR} phases);
a sharp decrease followed by a stationary moment occurs during the \ac{VPF} phase;
an additional sharp decrease coincides with the \ac{AC} phase, corresponding to the booster pump function.

The atrial pressure evolution over time is characterized by three waves and two pressure descents \citep{Pagel2003,Askari2019,Chambers2019}:
the a-wave -- corresponding to the increase of pressure due to the atrial contraction;
the c-wave -- caused by the closure of the atrioventricular valves (\ac{TV}, \ac{MV})  that push the blood back toward the atria;
the x-descent -- determined by the initial phase of the ventricular contraction and the consequent downward movement and filling of the atria;
the v-wave -- caused by the continuous venous return while the atrioventricular valves are closed, during the ventricular systole;
the y-descent -- which begins with the opening of the atrioventricular valves and continues during the \ac{VPF} phase.
In \cref{fig:res_baseline_pvloopext}, first two rows, central column, all these complex features are captured.
Additionally, we also obtain an a-wave taller than the v-wave in the \ac{RA} \citep{Chambers2019} and the opposite behavior in the \ac{LA} \citep{Gibson2003,Pagel2003}, as described in the medical literature \citep{Gibson2003,Pagel2003,Chambers2019,Askari2019,Vest2019}.
This behavior is also visible in the atrial pressure-volume loops (\cref{fig:res_baseline_pvloop}, (b), top), where the pressure assumes its maximum value during the A-loop for the \ac{RA} and during the V-loop for the \ac{LA}.
%The C-wave, instead, is clearly visible in the \ac{RA} pressure evolution, while is less evident in the \ac{LA} curve, as sometimes described also in can occur also in in-vivo 
%However, the C-wave is not always detectable in the \ac{LA} even in pressure curves acquired in-vivo , see, e.g., 

Another captured physiological behavior concerns the opening and closing times of the cardiac valves.
%In the ventricular pressure-volume loops we can notice how the \ac{IVR} occurs between \ac{PV} closure and \ac{TV} opening for the \ac{RV} and between \ac{AV} closure and \ac{MV} opening for the \ac{LV}.
Indeed, looking at the pressure and volume evolution over time (\cref{fig:res_baseline_pvloopext}), we observe how the right valves (\ac{TV}, \ac{PV}) close after the left ones (\ac{MV}, \ac{AV}).
More specifically, the closing of the atrioventricular valves (\ac{TV}, \ac{MV}) is almost synchronized, while a longer delay between the closing of the semilunar valves (\ac{PV}, \ac{AV}) occurs.
This behavior corresponds to normal cardiac physiology and can be routinely verified by checking the first and second heart sounds through cardiac auscultation \citep{Chizner2008}.
We obtain these results thanks to our stimulation protocol that, albeit simplified, correctly reproduces the activation delay between \ac{LV} and \ac{RV}.

\begin{figure}[t]
    \centering
    \includegraphics[width=1.0\textwidth]{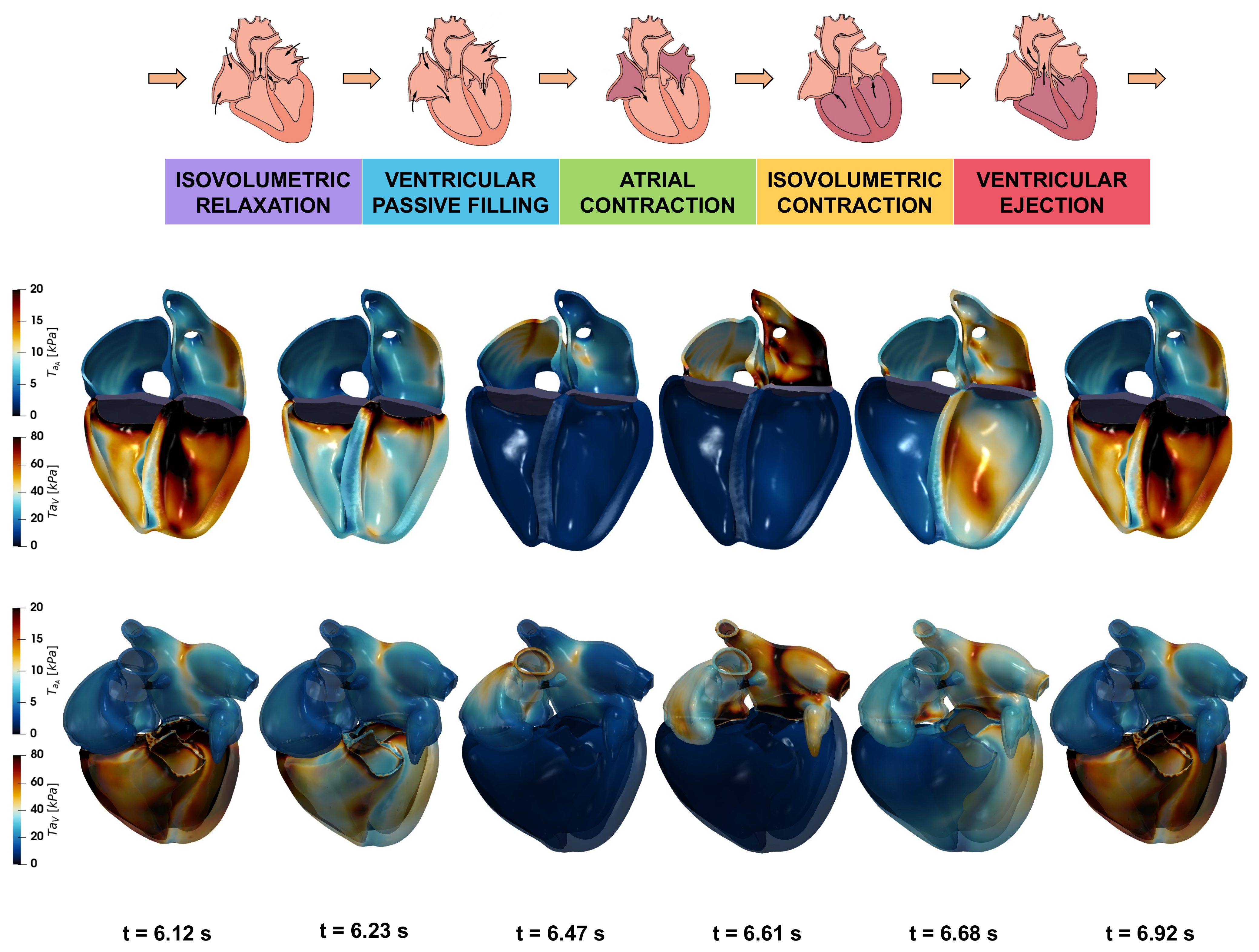}
    \caption{Deformed configuration of the cardiac muscle over time, colored with the active tension $\Tens$ saturated to \SI{20}{\kilo\pascal} and \SI{80}{\kilo \pascal} for the atria and ventricles, respectively:
    the internal view (top) and the external view (bottom) at the initial and final instants of each phase of the cardiac cycle.}
    %The corresponding \href{https://www.dropbox.com/s/ob38iz8z60y1bbp/video_baseline.mp4?dl=0}{Video 1} (online version) shows also the evolution in time of the transmembrane potential $u$ and of the intracellular calcium concentration $\wCa$.}
    \label{fig:res_baseline_active_tension}
\end{figure}

In \cref{fig:res_baseline_active_tension} we show the three-dimensional motion of the heart during a heartbeat.
The cardiac muscle is colored according to the local value of the active tension $\Tens$, to highlight which chambers are contracting and which are relaxing during the five phases of the cardiac cycle.
Specifically, we show the cardiac geometry deformed by the displacement $\Disp$ at the initial and final moments of each phase.
% In the corresponding \href{https://www.dropbox.com/s/ob38iz8z60y1bbp/video_baseline.mp4?dl=0}{Video 1} (online version) we show also the evolution of the transmembrane potential $u$ and of the intracellular calcium concentration $\wCa$.
Our simulation reproduces the expected motion of a healthy heart, as summarized below.
\begin{itemize}
    \item The \ac{IVR} starts with relaxed atria and contracted ventricle.
    During this short phase the ventricular active forces quickly drop down together with the pressure.
    Since both the atrioventricular valves (\ac{TV} and \ac{MV}) and the semilunar valves (\ac{PV} and \ac{AV}) are closed, the ventricular volumes are constant and no significant motion occurs during this phase.
    \item The \ac{VPF} starts when the \ac{TV} and \ac{MV} open and is characterized not only by the active tension that continues to fall, but also by a clear increase of the ventricular volumes and a corresponding decrease of the atrial ones.
    These volumetric changes are mainly caused by the upward movement of the ventricular base, which compresses the atria and dilates the ventricles.
    During this phase the atrial are passively deformed, acting as a conduit.
    \item In late diastole, the \ac{AC} phase starts from the pacemaking stimulation in the \ac{RA} (near the \ac{SupVC}) and propagates toward the \ac{LA}.
    The atrial booster pump function gives an additional preload to the ventricles, visible once again with a clear upward movement of the base.
    We recall that the active tension is influenced by the local fiber stretch.
    Indeed, being the atrial deformation mainly longitudinal, active tension is higher where the fibers are not oriented transmurally.
    Thanks to our anatomically accurate fiber model, this feature is clearly visible in the \ac{PeMs}, where the active tension follows their characteristic orientation.
    %(see \href{https://www.dropbox.com/s/ob38iz8z60y1bbp/video_baseline.mp4?dl=0}{Video 1})
    \item The ventricular contraction starts during the short \ac{IVC} phase, when the atrioventricular valves close again and the ventricular active tension starts to rise from the left to the right part.
    Since also during this phase the ventricular volumes are constant, no clear deformation are visible.
    Meanwhile, the atria begin to fill up, starting their reservoir function fueled by the continuous venous return.
    \item Finally, the \ac{VE} phase is characterized by the opening of the semilunar valves (\ac{PV} and \ac{AV}), the strong ventricular contraction, and the consequent decrease of the ventricular volumes.
    Again, the key factor driving this emptying phase is the downward movement of the atrioventricular plane \citep{Carlsson2007}, which also determines most of the atrial filling during its reservoir function.
\end{itemize}
%\todoinline{@FR / @RP Teniamo quest'ultima frase o non serve?}
The just described physiological motion of the entire heart during the whole cardiac cycle has been obtained thanks to several features of our electromechanical model.
According to our experience, the key factors are the following:
i) the anatomical accuracy of the geometry;
ii) the use of comprehensive and calibrated mathematical models for the atria and the ventricles, in terms of electrophysiology, active-force generation, passive mechanics;
iv) the modeling of the most relevant feedbacks among the different core models, with particular reference to the fibers-stretch and fibers-stretch-rate feedbacks in the force generation model; 
v) the mechanical boundary conditions on the epicardium taking into account both the presence of the \ac{PF} and of the \ac{EAT}.
The latter, in particular, is of fundamental importance for the correct downward and upward movement of the ventricular base.

\subsubsection{A quantitative analysis of volumetric indexes}

\begin{table}[t]
\centering
\footnotesize
\begin{tabular}{@{}lllll@{}}
\toprule
\multirow{2}{*}{Index} & \multirow{2}{*}{Value} & \multicolumn{2}{l}{Reference ranges} & \multirow{2}{*}{Description} \\ \cmidrule(lr){3-4}
                       &                        & mean $\pm$ SD            & [LL, UL]            &                              \\ \midrule
\ac{RA} $V_\mathrm{max}\, $ \mLperMsquare &
  43.68 &
      $ 52 \pm 12$ & $[28, 76]$ & 
        \ac{RA} maximum volume \\
\ac{RA} $V_\mathrm{preAC}\, $ \mLperMsquare &
  32.83 &
      $ 40\pm 10^*$ & $[19 , 61]^*$ & 
        \ac{RA} volume before atrial contraction\\
\ac{RA} $V_\mathrm{min}\, $ \mLperMsquare &
  25.16 &
      $ 27\pm 9$ & $[9, 45]$ & 
        \ac{RA} minimum volume \\
% \ac{RA} $\mathrm{TotEV}\, $ \mLperMsquare &
%   18.52 &
%
%       $ \pm$ & $[ , ]$ &
%         \ac{RA} total emptying volume ($V_\mathrm{max} - V_\mathrm{min}$)\\
% \ac{RA} $\mathrm{PassEV}\, $ \mLperMsquare &
%   10.84 &
%
%       $ \pm$ & $[ , ]$ &
%         \ac{RA} passive emptying volume ($V_\mathrm{max} - V_\mathrm{preAC}$)\\
% \ac{RA} $\mathrm{ActEV}\, $ \mLperMsquare &
%   7.68 &
%
%       $ \pm$ & $[ , ]$ &
%         \ac{RA} active emptying volume ($V_\mathrm{preAC} - V_\mathrm{min}$)\\
\ac{RA} $\mathrm{PassEF}$ [\%] &
  24.82 &
      $ 23 \pm 9^*$ & $[ 4 , 41 ]^*$ & 
        \ac{RA} passive ejection fraction: $(V_\mathrm{max} - V_\mathrm{preAC}) / V_\mathrm{max}$\\
\ac{RA} $\mathrm{ActEF}$ [\%] &
  23.39 &
      $ 33 \pm 10^*$ & $[11 , 55 ]^*$ & 
        \ac{RA} active ejection fraction: $(V_\mathrm{preAC} - V_\mathrm{min}) / V_\mathrm{preAC}$\\
\ac{RA} $\mathrm{TotEF}$ [\%] &
  42.40 &
      $ 49 \pm 10$ & $[29 , 68 ]$ & 
        \ac{RA} total ejection fraction: $(V_\mathrm{max} - V_\mathrm{min}) / V_\mathrm{max}$\\
\midrule
\ac{LA} $V_\mathrm{max}\, $ \mLperMsquare &
  30.65 &
      $41 \pm 8$ & $[24,57]$ & 
        \ac{LA} maximum volume \\
\ac{LA} $V_\mathrm{preAC}\, $ \mLperMsquare &
  24.70 &
      $ 30 \pm 8^*$ & $[15, 46]^*$ & 
        \ac{LA} volume before atrial contraction\\
\ac{LA} $V_\mathrm{min}\, $ \mLperMsquare &
  17.14 &
      $19 \pm 5$ & $[9,28]$ & 
        \ac{LA} minimum volume \\
% \ac{LA} $\mathrm{TotEV}\, $ \mLperMsquare &
%   13.51 &
%
%       $ \pm$ & $[ , ]$ &
%         \ac{LA} total emptying volume ($V_\mathrm{max} - V_\mathrm{min}$)\\
% \ac{LA} $\mathrm{PassEV}\, $ \mLperMsquare &
%   5.96 &
%
%       $ \pm$ & $[ , ]$ &
%         \ac{LA} passive emptying volume ($V_\mathrm{max} - V_\mathrm{preAC}$)\\
% \ac{LA} $\mathrm{ActEV}\, $ \mLperMsquare &
%   7.56 &
%
%       $ \pm$ & $[ , ]$ &
%         \ac{LA} active emptying volume ($V_\mathrm{preAC} - V_\mathrm{min}$)\\
\ac{LA} $\mathrm{PassEF}$ [\%] &
  19.43 &
      $ 26 \pm 9^*$ & $[8 , 44]^*$ & 
        \ac{LA} passive ejection fraction: $(V_\mathrm{max} - V_\mathrm{preAC}) / V_\mathrm{max}$\\
\ac{LA} $\mathrm{ActEF}$ [\%] &
  30.59 &
      $ 37 \pm 10^*$ & $[17, 58]^*$ & 
        \ac{LA} active ejection fraction: $(V_\mathrm{preAC} - V_\mathrm{min}) / V_\mathrm{preAC}$\\
\ac{LA} $\mathrm{TotEF}$ [\%] &
  44.08 &
      $ 54 \pm 8$ & $[37 , 70]$ & 
        \ac{LA} total ejection fraction: $(V_\mathrm{max} - V_\mathrm{min}) / V_\mathrm{max}$\\
\midrule
\ac{RV} \ac{EDV} $$ \mLperMsquare &
  86.78 &
      $88 \pm 17$ & $[53,123]$ & 
        \ac{RV} \acl{EDV}\\
\ac{RV} \ac{ESV} $$ \mLperMsquare &
  41.73 &
      $38 \pm 11$ & $[17,59]$ & 
        \ac{RV} \acl{ESV}\\
\ac{RV} \ac{SV} $$ \mLperMsquare &
  45.06 &
      $52 \pm 12$ & $[28,75]$ & 
        \ac{RV} \acl{SV} ($\mathrm{EDV} - \mathrm{ESV}$)\\
\ac{RV} EF $[\%]$ &
  51.92 &
      $57 \pm 8$ & $[42,72]$ & 
        \ac{RV} ejection fraction ($\mathrm{SV} / \mathrm{EDV}$)\\
\midrule
\ac{LV} \ac{EDV} $$ \mLperMsquare &
  85.28 &
      $77 \pm 15$ & $[47,107]$ & 
      % $83 \pm 14$ & $[55,110]$ & % age 40-50
        \ac{LV} \acl{EDV}\\
\ac{LV} \ac{ESV} $$ \mLperMsquare &
  41.28 &
      $29 \pm 9$ & $[11,47]$ & 
      % $32 \pm 9$ & $[13,50]$ & % age 40-50
        \ac{LV} \acl{ESV}\\
\ac{LV} \ac{SV} $$ \mLperMsquare &
  44.00 &
      $48 \pm 9$ & $[30,66]$ & 
      % $52 \pm 8$ & $[36,68]$ & % age 40-50
        \ac{LV} \acl{SV} ($\mathrm{EDV} - \mathrm{ESV}$)\\
\ac{LV} EF $[\%]$ &
  51.60 &
      $63 \pm 6$ & $[51,76]$ & 
      % $62 \pm 7$ & $[48,76]$ & % age 40-50
        \ac{LV} ejection fraction ($\mathrm{SV} / \mathrm{EDV}$)\\
\bottomrule
\end{tabular}
\normalsize
\caption{Volumetric indexes of the four cardiac chambers:
the values computed from the baseline simulation compared to the reference ranges for cardiac magnetic resonance \citep{Kawel-Boehm2020,Li2017a} (\textit{SD} standard deviation, \textit{LL} lower limit, \textit{UL} upper limit).
All the volumes are indexed by the body surface area.
Reference ranges are taken from the recent meta-analysis by \citet{Kawel-Boehm2020}, with the exception of the values marked with $^*$ taken from \citet{Li2017a}, the only paper included in the meta-analysis in which additional parameters for atrial conduit and booster pump function are analyzed.}
\label{tab:results_index}
\end{table}

Quantitative volume-based indexes of the four cardiac chambers are routinely used in clinics to assess the physiology of the heart.
Reference values for these indexes are available in the medical literature, but their values significantly vary depending on the kind of medical images used or the methods employed to compute the volume.
Echocardiography and cardiac magnetic resonance are the most used techniques, but the former usually underestimates the chambers' volume because of the low spatial resolution.
Indeed, reference values for cardiac magnetic resonance \citep{Kawel-Boehm2020,Li2017a} are consistently larger than the ones for echocardiography \citep{Lang2015,Blume2011a,Peluso2013}.
% However, being the development of the atrial volume-based indexes relatively recent, more atrial parameters are currently available for the echocardiography, a less time-consuming exam easier to routinely acquire CITE.
In both these techniques the volumes can be computed either using surrogate formulas (based on the chambers' area on specific image slices) or employing the more accurate Simpson's method (consisting in the segmentation of stack of contiguous slices that cover the whole cardiac chamber).
Based on these considerations, in \cref{tab:results_index} we compare the volume-based indexes computed from the baseline simulation with the current reference ranges for cardiac magnetic resonance for adult men, computed using the Simpson's method.
To this purpose we use values reported in the recent meta-analysis by \citet{Kawel-Boehm2020} and some additional values focused on the atrial conduit and booster pump function from the study of \citet{Li2017a}, being this paper the unique one of the meta-analysis reporting this kind of indexes -- often computed in echocardiography \citep{Blume2011a,Peluso2013} -- for cardiac magnetic resonance.
Concerning the atria, we compute all the indexes by not considering the appendages and the veins (which contribute about $25\%$ of the total volume), as usually done in the image-based indexes such as the ones taken as reference ranges.

All the indexes calculated, both for the atria and for the ventricles, fall within the reference ranges.
In conclusion, the results presented in this section demonstrate the ability of our whole-heart electromechanical model to capture all fundamental aspects of the healthy physiology of the heart.
To the best of our knowledge, the cardiac function has never been modeled so comprehensively by a computational model of the heart.

%%%%%%%%%
\subsection{The impact of the atrial contraction and of the fiber-stretch-rate feedback}
\label{subsec:res_effects_math_mod}

In this section we aim at showing the critical role that atrial contraction and fiber-stretch-rate feedback play in simulating the physiological cardiac function.
We do not consider the MEF, as it mostly plays a role in pathological conditions involving arrhythmogenic behavior \cite{Salvador2021, Salvador2022}.

\begin{figure}%[t]
    \centering
    \includegraphics[height=0.8\textheight]{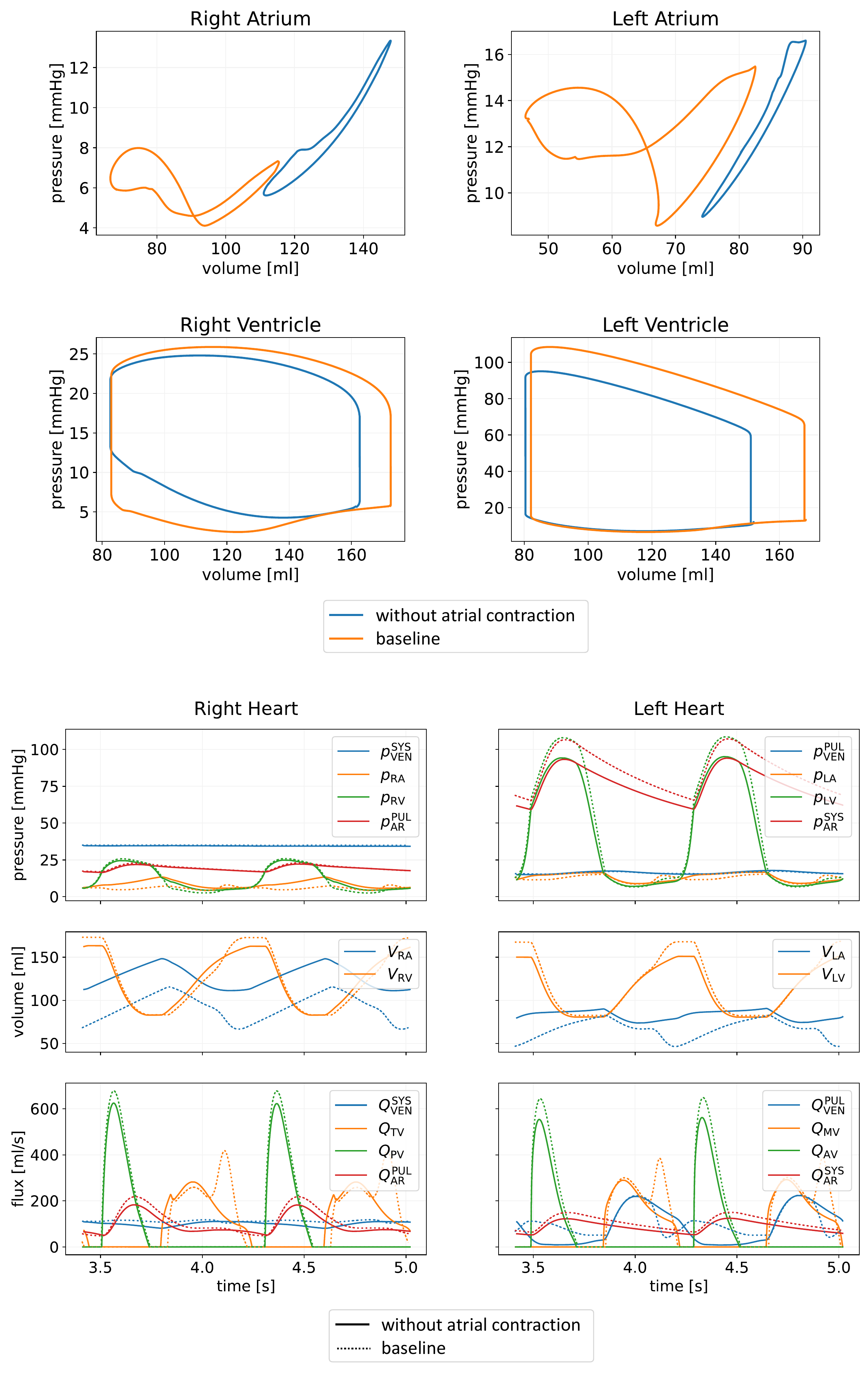}
    \caption{Pressure-volume loops (top) and circulation state variables (bottom) with and without the atrial contraction in the model.}
    \label{fig:res_effect_ac}
\end{figure}

Most of the whole-heart electromechanical models, as discussed in \cref{sec:intro}, neglect the atrial contraction \citep{Sugiura2012a,Fritz2014,Augustin2016,Santiago2018,Pfaller2019,Strocchi2020b}.
To discuss the impact of this choice, in \cref{fig:res_effect_ac} we show the effect of switching off the atrial contraction in the model.
This is simulated by considering the atria as purely passive tissues, by ignoring in the atrial domain $\{\OmegaRA \cup \OmegaLA\}$ the active stress part of the Piola-Kirchhoff stress tensor (see \cref{eq:piola_myo}).
The results show irrefutably the importance of atrial contraction for both atrial and ventricular function:
on the one hand the A-loop disappears from the atrial pressure-volume loops;
on the other hand, the ventricular cycle also changes drastically with a significant decrease in ventricular preload.
This non-physiological behavior is also evident in the evolution of the state variables (\cref{fig:res_effect_ac}, bottom), where the contributions of the atrial contraction in terms of pressures, volumes and fluxes disappear.
In other words, neglecting atrial contraction means neglecting the booster pump function of the atria and its preloading effect on the ventricles, modeling a pathological scenario rather than a healthy one.
For instance, as shown by \citet{Pagel2003}, similar pressure-volume loops are captured during atrial fibrillation, when the chaotic propagation of the electrical signal causes the atrium to lose its booster pump function.

\begin{figure}%[t]
    \centering
    \includegraphics[height=0.8\textheight]{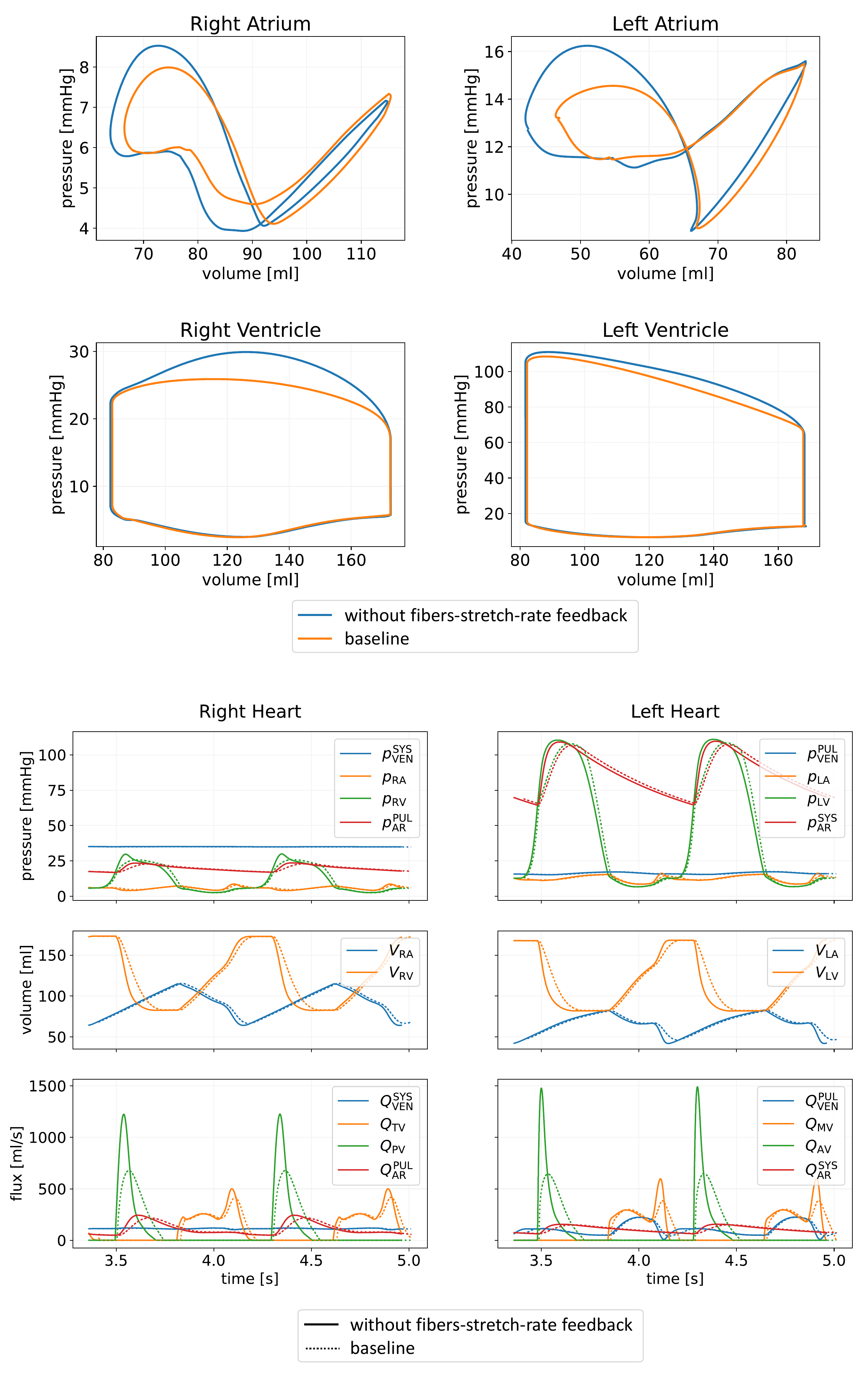}
    \caption{Pressure-volume loops (top) and circulation state variables (bottom) with and without the fibers-stretch-rate feedback in the model.}
    \label{fig:res_effect_fsrf}
\end{figure}

In \cref{fig:res_effect_fsrf} we show the effects on the results of the fibers-stretch-rate feedback off in the model.
In terms of pressure-volume loops (\cref{fig:res_effect_fsrf}, top) no substantial changes are visible, with a small increase of the A-loop size in the atria and a small increase of the ventricular pressures during the \ac{VE} phase.
Conversely, looking at the state variables (\cref{fig:res_effect_fsrf}, bottom), the fluxes through the semilunar valves (\ac{AV} and \ac{PV}) dramatically change.
More specifically, without the fiber-stretch-rate feedback we obtain about \SI{1200}{\milli\liter\per\second} in the \ac{AV} and almost \SI{1500}{\milli\liter\per\second} in the \ac{PV}.
This abnormal values are similar to the ones reported by the whole-heart model of \citet{Gerach2021}, that indeed neglects the fibers-stretch-rate feedback since, without suitable stabilization terms, it yields strong non-physical oscillations in the multi-scale model resulting in an unstable numerical scheme~\citep{Gerach2021}.
A large increase of the fluxes appears also in the atrioventricular valves during the \ac{AC} phase.

We conclude that the fibers-stretch-rate plays a fundamental role in the regulation of the cardiac function.
As a consequence of this feedback, indeed, the active force decreases for the cardiac cells located in regions where fibers are rapidly shortening, thus resulting in slowing down the contraction velocity of fibers.
As this feedback acts locally (i.e. at the cell level), the resulting macroscopic effect is a homogenization of fibers shortening velocity, preventing sharp variations.
From a hemodynamic perspective, this results into a smoothing of the ejected blood flux, as highlighted by our results.
Hence, we postulate that the fibers-stretch-rate feedback, despite originating from the microscale force-velocity relationship of sarcomeres, plays a crucial role in the regulation of blood fluxes.

%%%%%%%%%
\subsection{The need of the numerical stabilization}
\label{subsec:res_effects_num_mod}

To highlight the role of the stabilization terms described in \cref{sec:numerics:stabilization}, we present the results of two numerical simulations obtained by switching off the stabilization terms on active stress (see \cref{eqn:stab:active-stress} and \cite{regazzoni2021oscillation}) and on the 3D-0D mechanics-circulation coupling (see \cref{eqn:volume_constraints_stabilized} and \cite{regazzoni2022stabilization}), respectively.

\begin{figure}%[t]
    \centering
    \includegraphics[height=0.8\textheight]{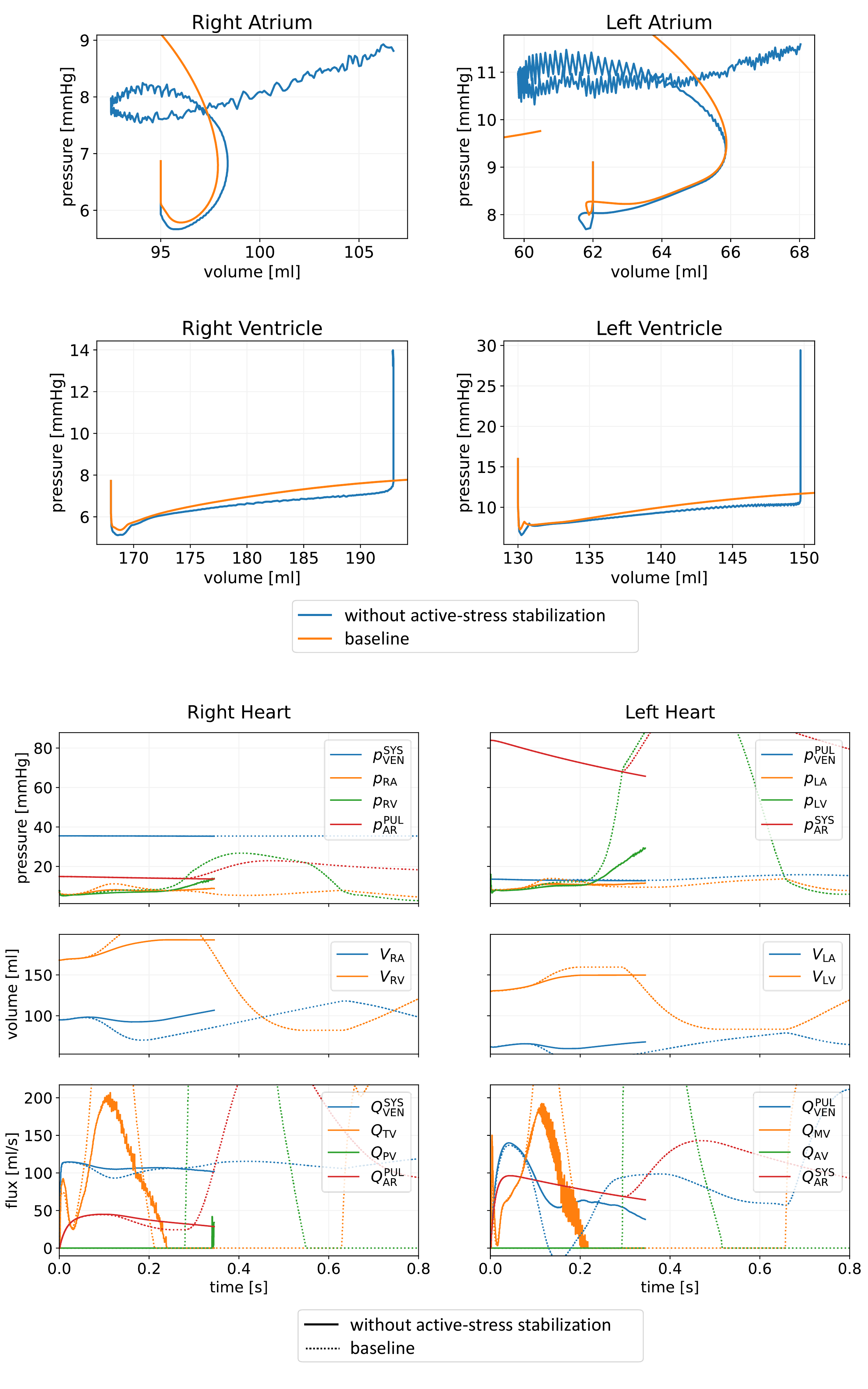}
    \caption{Pressure-volume loops (top) and circulation state variables (bottom) with and without the stabilization terms in the active stress model.}
    \label{fig:res_effect_stab_active_stress}
\end{figure}

In \cref{fig:res_effect_stab_active_stress} we show results obtained without the active stress stabilization.
As soon as active tension is being developed, non-physical oscillations occur, mainly in pressure and flux traces, finally leading to failure of the nonlinear mechanics solver.
In the simulation shown in the figure, failure occurs after nearly \SI{0.35}{\second} of physical time.
While the time of failure depends on the time step size and on the parameters, in our experience the numerical simulation of cardiac active mechanics with realistic parameter values always leads to this kind of numerical oscillations, whenever stretch-rate-feedback is accounted for by the model.
These instabilities cannot be cured by reducing the time step size: on the contrary, as analytically demonstrated \cite{regazzoni2021oscillation}, they are amplified by a rapid exchange of variables between the tissue mechanics and the activation model.

\begin{figure}%[t]
    \centering
    \includegraphics[height=0.8\textheight]{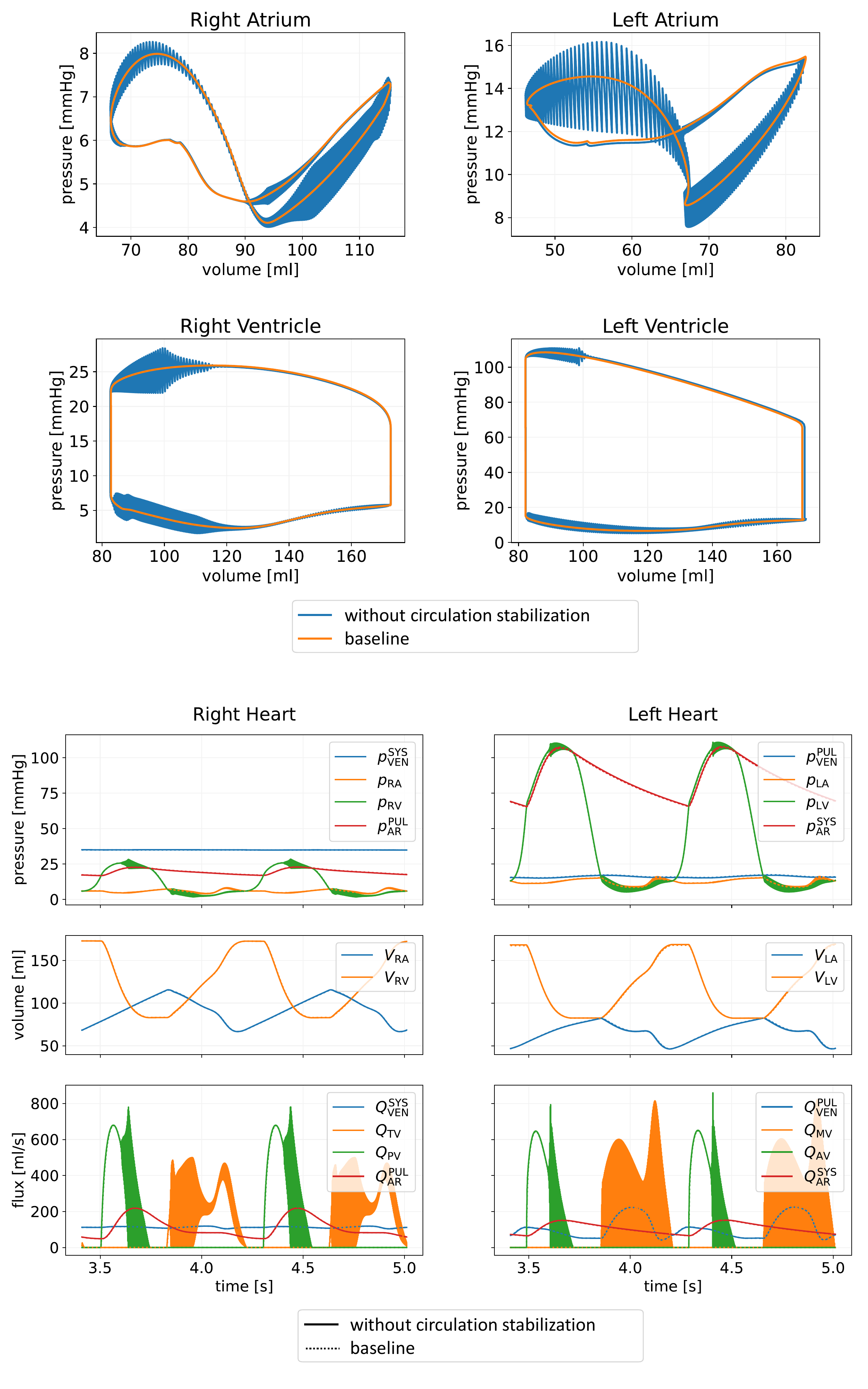}
    \caption{Pressure-volume loops (top) and circulation state variables (bottom) with and without the stabilization terms in the circulation model.}
    \label{fig:res_effect_stab_circ}
\end{figure}

In \cref{fig:res_effect_stab_circ} we report the results obtained by switching off the stabilization term on the 3D-0D mechanics-circulation coupling.
Unlike for the active stress stabilization term, the numerical oscillations obtained in this case do not lead to failure of the simulation.
However, the results are clearly not physically meaningful.
This is particularly evident from the transients of blood fluxes across valves, that exhibit very large oscillations.
As demonstrated in \cite{regazzoni2022stabilization}, also in this case reducing the time step size does not solve this issue, but on the contrary it typically contributes to the onset of oscillations.

Thanks to the stabilization terms of \cref{eqn:stab:active-stress} and \cref{eqn:volume_constraints_stabilized} we are able to remove the non-physical oscillations, for any choice of parameters and of time step size.
These numerical tests demonstrate that the interplay between accurate mathematical models and efficient and stable numerical methods is of fundamental importance to model the cardiac function and to obtain physiological results.

%%%%%%%%
%%%%%%%%
\section{Conclusions}
\label{sec:concl}

In this paper, we proposed a biophysically detailed, numerically stable and accurate computational model of the electromechanics of the whole human heart, by considering an active contraction model for both atria and ventricles.

In developing whole-heart computation models, several aspects are crucial to comprehensively model the cardiac function and to accurately capture the highly coordinated events underlying the cardiac cycle.
In this context, our model embeds different determinant features.
We use an anatomically accurate computational domain including the main cardiac components such as atrial appendages, major arteries, and simplified cardiac valves (\cref{fig:domain}).
In order to characterize the varying biophysical properties of the cardiac tissue, we split the whole domain into several regions, representing cardiac chambers, arteries, and insulating fibrous tissue of the cardiac valves.
To capture the anisotropy of the muscular tissue, we model the myocardial fiber architecture by taking advantage of the anatomically-accurate whole-heart \ac{LDRBM} that we have recently proposed in \citep{Piersanti2021fibers,Piersanti2021phd} (\cref{fig:fibers}).
Our full electromechanical model comprises of several biophysically detailed core models.
We employ chamber-specific and accurate ionic models for atria and ventricles~\citep{CRN,TTP06}, coupled with the monodomain equation to describe the transmembrane potential propagation at the macroscale.
We use the \ac{RDQ20} model~\citep{Regazzoni2020} for the active force generation, a biophysically detailed microscale model that captures the crucial influence of the fiber-stretch and fibers-stretch-rate on the generation of the active forces.
We employ a 0D closed-loop model of the circulatory system, fully-coupled with the mechanical model~\citep{Regazzoni2022}.
We use specific constitutive laws and model parameters for each cardiac region.
The core models are mutually coupled by considering the most important feedbacks that represent the interactions among electric signal propagation, microscopic and macroscopic cardiac tissue contraction and deformation, and blood circulatory system (see \cref{fig:model}).
Among them, in this paper we pay special attention to the fibers-stretch-rate feedback (between passive mechanics and active force generation model).
% investigating its impact on reproducing physiological results and overcoming the releted numerical issues
Concerning the numerical discretization, we use the efficient segregated-intergrid-staggered scheme proposed in \citep{Regazzoni2022,Piersanti2022biv} and we employ recently developed stabilization terms -- related to the circulation~\citep{regazzoni2022stabilization} and the fibers-stretch-rate feedback~\citep{regazzoni2021oscillation} -- that are crucial to obtain a stable formulation in a four-chamber scenario (see \cref{fig:num_scheme}).
To cope with the high computational complexity associated with whole-heart electromechanical simulations, we have developed our solver in \lifex, an efficient in-house \ac{FE} library focused on large-scale cardiac applications in an HPC framework.

We simulate all the phases of the cardiac cycle, showing numerical results that comprehensively capture the atrial and ventricular physiology, the threefold atrial function of reservoir, conduit and booster pump, and the atrioventricular interaction.
To the best of our knowledge, some of the physiological features that we catch have never been shown all together by a computational model of the heart.
Specifically, we mention
the fluxes through the semilunar valves (\cref{fig:res_baseline_state}),
the eight-shaped atrial pressure-volume loops characterized by the correct proportion between A- and V-loops (\cref{fig:res_baseline_pvloop}),
the evolution over time of the atrial volumes (\cref{fig:res_baseline_pvloopext}),
the a-, c-, v-waves of the atrial pressure (\cref{fig:res_baseline_pvloopext}),
and the three-dimensional deformation driven by the upward and downward movement of the atrioventricular plane (\cref{fig:res_baseline_active_tension}).
More quantitatively, we compute volumetric indexes for all the cardiac chambers, finding values that always fall within the reference physiological ranges (\cref{tab:results_index}).

We also analyze the impact of atrial contraction, fibers-stretch-rate feedback and stabilization terms, by comparing the results obtained with and without these features in the model.
Due to the complex anatomy and physiology of the atria, atrial contraction is often neglected in electromechanical models of the whole heart~\citep{Fritz2014,Augustin2016,Santiago2018,Pfaller2019,Strocchi2020b}.
However, we show that neglecting atrial contraction (and the associated atrial booster pump function acting as preload for the ventricles) means modeling a pathological rather than healthy scenario (\cref{fig:res_effect_ac}).
Concerning the fibers-stretch-rate feedback, we show that without this feedback the fluxes across the semilunar valves largely exceed the physiological range (\cref{fig:res_effect_fsrf}).
This feedback originates from the microscale force-velocity relationship of sarcomeres, decreasing the active force in regions where fibers are rapidly shortening.
The macroscopic effect is a homogenization of fibers shortening velocity that, from a hemodynamic perspective, results into a smoothing of the ejected blood flux, as highlighted by our results.
Hence, we postulate that the fibers-stretch-rate feedback, despite originating at the microscale, plays a crucial role in the macroscopic regulation of blood fluxes.
Moreover, if not properly managed at the numerical level, this feedback produces non-physical oscillations that may lead the numerical simulation to fail~\citep{regazzoni2021oscillation,Gerach2021}.
Thus, the interplay between accurate mathematical models and efficient and stable numerical methods is of utmost importance to reproduce the heart physiology.
We show that, thanks to the introduction of the stabilization terms relative to the circulation model
and to the fibers-stretch-rate feedback,
we are able to remove the non-physical oscillations (\cref{fig:res_effect_stab_circ,fig:res_effect_stab_active_stress}).

To conclude, the presented electromechanical model of the whole human heart has shown an unprecedented ability in reproducing the healthy cardiac function of both atria and ventricles and can be considered a fundamental step toward the construction of physics-based digital twins of the human heart.

\section*{Acknowledgements}
\noindent
This project has received funding from the European Research Council (ERC) under the European Union's Horizon 2020 research and innovation program (grant agreement No 740132, iHEART - An Integrated Heart Model for the simulation of the cardiac function, P.I. Prof. A. Quarteroni).
We acknowledge the CINECA award under the class ISCRA B project (CoreMaS - code: HP10BD303V) for the availability of high performance computing resources.
MS and AZ have received funding by the Italian Ministry of University and Research (MIUR) within the PRIN (Research projects of relevant national interest 2017 “Modeling the heart across the scales: from cardiac cells to the whole organ” Grant Registration number 2017AXL54F)
% \begin{center}
% 	\raisebox{-.5\height}{\includegraphics[width=.15\textwidth]{EU_col.jpg}}
% 	\hspace{2cm}
% 	\raisebox{-.5\height}{\includegraphics[width=.15\textwidth]{ERC_col.jpg}}
% \end{center}
\begin{center}
	\includegraphics[width=0.3\textwidth]{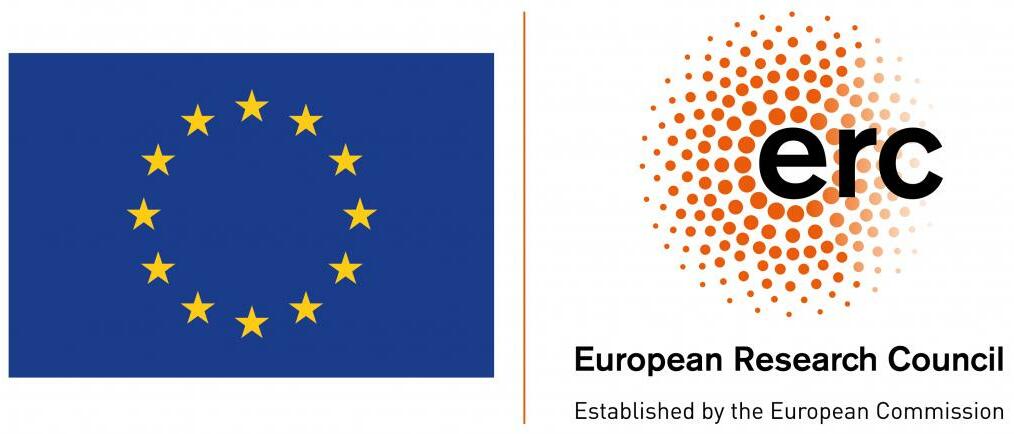}
\end{center}

\appendix

\section{Model and numerical parameters}
\label{sec:appendix_params}

We provide more details about the parameters used for the baseline simulation (\cref{subsec:res_baseline}).
Concerning the fiber generation procedure (\cref{subsec:fibers}), we use the parameters reported in \citet[Chapter~4]{Piersanti2021phd}.
In \cref{tab:params_ep} we report the parameters of the monodomain model ($\PhyEP$) of \cref{eq:ep}.
Concerning the ionic models \eqref{eq:ion}, we use the parameters reported in \acl{CRN} and \acl{TTP06} (endocardium cells) for the atria (\cref{eq:ep_ionic_atria}) and ventricles (\cref{eq:ep_ionic_ventr}), respectively.
Additionally, for the \ac{TTP06} ionic model, we rescale the calcium peak by a factor of $0.48$ to bring the calcium transient into a more physiological range.
%\footnote{Notice that, for the \ac{TTP06} ionic model, we rescaled the calcium peak by a factor of 0.48. This brings the calcium transient to more physiological range  after the 1000-cycle long single-cell simulation initialization of the ionic model.}
% ELECTROPHYSIOLOGHY
\begin{table}%[ht!]
    \centering
    \footnotesize
    \begin{tabular}{lrrr}
            \toprule
            Variable & Value & Unit & Description \\
            \midrule
            %\multicolumn{4}{l}{\textbf{Electrophysiology}} \\
            $T_{\text{hb}}$ & \num{0.8} & \si{\second} & Heartbeat duration \\
            %$\EPchim$ & \num{1} & \si{\micro \farad \per \centi \meter^2} & Surface-to-volume ratio \\
            %$\EPCm$ & \num{1400} & \si{\per\centi\meter}& Transmembrane capacitance \\
            $\epsilon$ & \num{0.05} & $-$ & Threshold of the fast conduction layer \\
            $\left(\sigma^{\Ventr,\Endo}_\text{f}, \sigma^{\Ventr,\Endo}_\text{s}, \sigma^{\Ventr,\Endo}_\text{n}\right)/(\EPchim \EPCm)$ & (\num{8.00}, \num{4.40}, \num{2.20}) $\times$ \num{e-4} & \si{\meter^2 \per \second} & Ventricular fast layer conductivities \\
            $\left(\sigma^{\Ventr,\Myo}_\text{f}, \sigma^{\Ventr,\Myo}_\text{s}, \sigma^{\Ventr,\Myo}_\text{n}\right)/(\EPchim \EPCm)$ & (\num{2.00}, \num{1.10}, \num{0.55}) $\times$ \num{e-4} & \si{\meter^2 \per \second} & Ventricular myocardial conductivities \\
            $\left(\sigma^{\Atria}_\text{f}, \sigma^{\Atria}_\text{s}, \sigma^{\Atria}_\text{n}\right)/(\EPchim \EPCm)$ & (\num{7.00}, \num{1.41}, \num{1.41}) $\times$ \num{e-4} & \si{\meter^2 \per \second} & Atrial conductivities \\
            $\EPIapp/\EPCm$ & \num{25.71} & \si{\volt \per \second} & Applied current value \\
            $\delta t$ & \num{3.0} & \si{\milli \second} & Applied current duration \\
            $t_{\RA}$ & \num{0.0} & \si{\milli \second} & Applied current RA initial time \\
            $t_{\LV}$ & (\num{160}, \num{160}, \num{160}) & \si{\milli \second} & Applied current LV initial times \\
            $t_{\RV}$ & (\num{165}, \num{172}) & \si{\milli \second} & Applied current RV initial times \\
            $r$ & \num{3e-3} & \si{\meter} & Applied current radius 
            \\
            \bottomrule
    \end{tabular}
    \normalsize
    \caption{Parameters of the electrophysiological model ($\PhyEP$).}
    \label{tab:params_ep}
\end{table}
\begin{table}%[ht]
    \centering
    \footnotesize
    \begin{tabular}{lrrr}
        \toprule
        Variable & Value & Unit & Description \\
        \midrule
        \multicolumn{4}{l}{\textbf{Ventricles} ($\OmegaVentr$)} \\
        $\SL_0$          & 1.9           & \si{\micro\meter}                & Reference sarcomere length \\
        %$\Clrv$            & 0.70            & $-$                           & RV contractility ratio \\ % ridondante se si definisce anche aXB^RV
        $(n_\f,n_\s,n_\n)$       & (\num{1}, \num{0}, \num{0.4})  & $-$                     & Share of active tension along the fiber directions \\
        $\aXB^\mathrm{LV}$       & \num{15.0e8}  & \si{\pascal}                     & LV upscaling constant of crossbridge stiffness  \\
        $\aXB^\mathrm{RV}$       & \num{10.5e8}  & \si{\pascal}                     & RV upscaling constant of crossbridge stiffness  \\
        $\overline{k}_\text{d}$       & \num{0.36}  & \si{\micro \mole}                     & Calcium-troponin dissociation constant  \\
        $\alpha_{k_{\text{d}}}$       & \num{-0.2083}  & \si{\micro \mole \per \micro \meter}                     & Sensitivity to sarcomere length of calcium-troponin dissociation constant\\
        $\gamma$       & \num{30}  & $-$                     & End-to-end tropomyosin cooperativity parameter  \\
        $k_{\text{off}}$       & \num{8}  & \si{\per \second}                     & Reaction rate associated with troponin kinetics \\
        $k_{\text{basic}}$       & \num{4}  & \si{\per \second}                     & Reaction rate associated with tropomyosin kinetics \\
        $\mu^{0}_{{f_{P}}}$       & \num{32.225}  & \si{\per \second}                     & Zero order moment of XB attachment rate \\
        $\mu^{1}_{{f_{P}}}$       & \num{0.768}  & \si{\per \second}                     & First order moment of XB attachment rate \\                            
        \midrule
        \multicolumn{4}{l}{\textbf{Atria} ($\OmegaLA \cup \OmegaRA$)} \\
        $\SL_0$          & 1.9           & \si{\micro\meter}                & Reference sarcomere length \\
        $(n_\f,n_\s,n_\n)$       & (\num{1}, \num{0}, \num{0.4})  & $-$                     & Share of active tension along the fiber directions \\
        $\aXB^\mathrm{LA}$       & \num{30.0e7}  & \si{\pascal}                     & LA upscaling constant of crossbridge stiffness  \\
        $\aXB^\mathrm{RA}$       & \num{30.0e7}  & \si{\pascal}                     & RA upscaling constant of crossbridge stiffness \\
        $\overline{k}_\text{d}$       & \num{0.865}  & \si{\micro \mole}                     & Calcium-troponin dissociation constant  \\
        $\alpha_{k_{\text{d}}}$       & \num{-1.25}  & \si{\micro \mole \per \micro \meter}                     & Sensitivity to sarcomere length of calcium-troponin dissociation constant  \\
        $\gamma$       & \num{20}  & $-$                     & End-to-end tropomyosin cooperativity parameter  \\
        $k_{\text{off}}$       & \num{180}  & \si{\per \second}                     & Reaction rate associated with troponin kinetics \\
        $k_{\text{basic}}$       & \num{20}  & \si{\per \second}                     & Reaction rate associated with tropomyosin kinetics \\
        $\mu^{0}_{{f_{P}}}$       & \num{32.225}  & \si{\per \second}                     & Zero order moment of XB attachment rate \\
        $\mu^{1}_{{f_{P}}}$       & \num{0.768}  & \si{\per \second}                     & First order moment of XB attachment rate \\                                    
        \bottomrule
    \end{tabular}
    \normalsize
  \caption{Parameters of the active force generation model ($\PhyAct$) used in the ventricular ($\OmegaVentr$) and atrial ($\OmegaLA \cup \OmegaRA$) domains, if modified from the calibration proposed in \citet{regazzoni2022machine}.}
  \label{tab:params_act}
\end{table}
% @RP/@FR check this paragraph.}
In \cref{tab:params_act} we report the calibration for the \ac{RDQ20} active generation model ($\PhyAct$) both for ventricles and atria.
We only report parameters modified with respect to the original paper of \citet{regazzoni2022machine}.
In particular, the atrial calibration is based on \citet{Mazhar2021}.
We also report the values of the microscale crossbridge stiffness $\aXB^i$, for $\IinCH$, used to define the tissue level active tension of each cardiac chamber (see \cref{eq:active_tension}).
In \cref{tab:params_mec} we report the parameters of the passive mechanical model ($\PhyMec$) of \cref{eq:mec} and the additional parameters specific to its quasi-static approximation ($\PhyMec^{\mathrm{static}}$) of \cref{eq:mec_static} used for the reference configuration recovery (see \cref{subsec:ref_conf_init_displ}). 
\begin{table}%[ht]
  \centering
  \footnotesize
  \begin{tabular}{lrrl}
    \toprule
    Variable & Value & Unit & Description \\
    \midrule
   $B$              & \num{50e3}   & \si{\pascal}                  & Bulk modulus in the myocardium $\OmegaMyo$ \\
   $b_{\text{ff}}$  & 8            & $-$                           & Fiber strain scaling in the myocardium $\OmegaMyo$ \\
   $b_{\text{ss}}$  & 6            & $-$                           & Radial strain scaling in the myocardium $\OmegaMyo$ \\
   $b_{\text{nn}}$  & 3            & $-$                           & Cross-fiber in-plain strain scaling in the myocardium $\OmegaMyo$ \\
   $b_{\text{fs}}$  & 12           & $-$                           & Shear strain in fiber-sheet plane scaling in the myocardium $\OmegaMyo$ \\
   $b_{\text{fn}}$  & 3            & $-$                           & Shear strain in fiber-normal plane scaling in the myocardium $\OmegaMyo$ \\
   $b_{\text{sn}}$  & 3            & $-$                           & Shear strain in sheet-normal plane scaling in the myocardium $\OmegaMyo$ \\
   $C^{\Ventr}$              & \num{0.88e3} & \si{\pascal}                  & Material stiffness in the ventricular domain $\OmegaVentr$ \\
   $C^{\RA}$              & \num{1.47e3} & \si{\pascal}                  & Material stiffness in the right atrial domain $\OmegaRA$\\
   $C^{\LA}$              & \num{1.76e3} & \si{\pascal}                  & Material stiffness in the left atrial domain $\OmegaLA$ \\
   \midrule
   $\mu^{\mathrm{valve,caps}}$  &  \num{10e5}           & \si{\pascal}                           & Shear modulus in the domains $\{\OmegaValve \cup \OmegaCaps\}$ \\
   $\kappa^{\mathrm{valve,caps}}$  & \num{50e5}             & \si{\pascal}                           & Bulk modulus in the domains $\{\OmegaValve \cup \OmegaCaps\}$ \\
   $\mu^{\mathrm{\AO,\PT}}$  & \num{5.25e5}              & \si{\pascal}                           & Shear modulus in the arterial domains $\{\OmegaAO \cup \OmegaPT\}$ \\
   $\kappa^{\mathrm{\AO,\PT}}$  & \num{10e5}              & \si{\pascal}                           & Bulk modulus in the arterial domains $\{\OmegaAO \cup \OmegaPT\}$ \\
   \midrule
   $\rho_\text{s}$  & $10^3$       & \si{\kilogram \per \cubic\meter} & Tissue density in the whole domain $\Omega_0$\\
   \midrule
   $\BCmecKepiNpf$    & \num{2e5}    & \si{\pascal\per\meter}        & Normal stiffness on $\GammaEpiPF$ \\
   $\BCmecCepiNpf$    & \num{2e3}    & \si{\pascal\second\per\meter} & Normal viscosity on $\GammaEpiPF$  \\
   $\BCmecKepiNfat$    & \num{2e2}    & \si{\pascal\per\meter}        & Normal stiffness on $\GammaEpiFat$ \\
   $\BCmecCepiNfat$    & \num{2e0}    & \si{\pascal\second\per\meter} & Normal viscosity on $\GammaEpiFat$  \\
   \midrule
   % \multicolumn{4}{l}{\textbf{Reference Configuration}} \\
   $\PRAtilde$          & 900            & \si{\pascal}                & Residual \ac{RA} pressure for the reference configuration recovery \\
   $\PLAtilde$          & 1200            & \si{\pascal}                & Residual \ac{LA} pressure for the reference configuration recovery \\
   $\PRVtilde$          & 650            & \si{\pascal}                & Residual \ac{RV} pressure for the reference configuration recovery \\
   $\PLVtilde$          & 1150            & \si{\pascal}                & Residual \ac{LV} pressure for the reference configuration recovery \\
   $\PAOtilde$         & 9500            & \si{\pascal}                & Residual \ac{AO} pressure for the reference configuration recovery \\
   $\PPTtilde$        & 1700            & \si{\pascal}                & Residual \ac{PT} pressure for the reference configuration recovery \\
   \bottomrule
  \end{tabular}
  \normalsize
  \caption{Parameters of the mechanical model ($\PhyMec$) and the reference configuration recovery ($\PhyMec^{\mathrm{static}}$).}
  \label{tab:params_mec}
\end{table}
The parameters of the circulation model ($\PhyCirc$) of \cref{eq:circ,eq:circ_extended} are reported in \cref{tab:params_circ}.
\begin{table}%[ht]
  \centering
  \footnotesize
  \begin{tabular}{lrr|lrr}
    \toprule
    Variable & Value & Unit & Variable & Value & Unit \\
    \midrule
    $\RarSYS$     & 0.48             & \si{\mmHg \second \per \milli\liter}   &
      $\LarSYS$     & \num{5e-3}      & \si{\mmHg \second\squared \per \milli\liter}    \\
    $\RarPUL$     & 0.032116          & \si{\mmHg \second \per \milli\liter}   &
      $\LarPUL$     & \num{5e-4}      & \si{\mmHg \second\squared \per \milli\liter}    \\
    $\RvnSYS$     & 0.26            & \si{\mmHg \second \per \milli\liter}   &
      $\LvnSYS$     & \num{5e-4}      & \si{\mmHg \second\squared \per \milli\liter}    \\
    $\RvnPUL$     & 0.035684          & \si{\mmHg \second \per \milli\liter}   &
      $\LvnPUL$     & \num{5e-4}      & \si{\mmHg \second\squared \per \milli\liter}    \\
    $\CarSYS$     & 1.50             & \si{\milli\liter \per \mmHg}   &
      $\Rmin$       & 0.0075          & \si{\mmHg \second \per \milli\liter}            \\
    $\CarPUL$     & 10.0            & \si{\milli\liter \per \mmHg}   &
      $\Rmax$       & 75000           & \si{\mmHg \second \per \milli\liter}            \\
    $\CvnSYS$     & 60.0            & \si{\milli\liter \per \mmHg}   &
                    &                 &                                                 \\
    $\CvnPUL$    & 16.0             & \si{\milli\liter \per \mmHg}   &
                    &                 &                                                 \\
    \bottomrule
  \end{tabular}
  \normalsize
  \caption{Parameters of the circulation model ($\PhyCirc$).}
  \label{tab:params_circ}
\end{table}
Finally, concerning the numerical parameters, we report the setting used for the linear and nonlinear solvers in \cref{tab:lin_sol,tab:nonlin_sol}, respectively.
\begin{table}%[ht]
  \centering
  \footnotesize
  \begin{tabular}{llll}
    \hline\noalign{\smallskip}
    Physics/Fields & Linear solver & Preconditioner & Abs. tol. \\
    \noalign{\smallskip}\hline\noalign{\smallskip}
    Monodomain model & CG & AMG & $10^{-10}$ \\
    Activation & GMRES & AMG & $10^{-10}$ \\
    Mechanics & GMRES & AMG & $10^{-8}$ \\
    \noalign{\smallskip}\hline
  \end{tabular}
  \normalsize
  \caption{Tolerances of the linear solver for the different physics.}
  \label{tab:lin_sol}
\end{table}
\begin{table}%[ht]
  \centering
  \footnotesize
  \begin{tabular}{lllll}
    \hline\noalign{\smallskip}
    Physics/Fields & Nonlinear solver & Rel. tol. & Abs. tol. \\
    \noalign{\smallskip}\hline\noalign{\smallskip}
    Mechanics & Newton & $10^{-8}$ & $10^{-6}$ \\
    Reference configuration & Newton & $10^{-8}$ & $10^{-6}$ \\
    \noalign{\smallskip}\hline
  \end{tabular}
  \normalsize
  \caption{Tolerances of the nonlinear solver for the mechanical problem.}
  \label{tab:nonlin_sol}
\end{table}

    \bibliographystyle{elsarticle-num-names}
    \footnotesize
    \bibliography{references}

\end{document}